\documentstyle[amssymb,amsfonts]{amsart}

\newenvironment{proof}{\noindent {\bf Proof} }{\endprf\par}
\def \endprf{\hfill  {\vrule height6pt width6pt depth0pt}\medskip}
\def\emph#1{{\it #1}}
\def\textbf#1{{\bf #1}}

\newcommand{\bea}{\begin{eqnarray}}
\newcommand{\eea}{\end{eqnarray}}
\def\beaa{\begin{eqnarray*}}
\def\eeaa{\end{eqnarray*}}
\def\BB{{\cal B}}
\def\PP{{\cal P}}
\def\I{{\cal I}}
\def\J{{\cal J}}
\def\Jp{\J^{+}}
\def\Jm{\J^{-}}

\def\MM{{\cal M}}
\def\ba{\begin{array}}
\def\ea{\end{array}}
\def\be#1{\begin{equation} \label{#1}}
\def \eeq{\end{equation}}

\newcommand{\nn}{\nonumber}

\def\rr{{\bf R}}

\def\nn{\nonumber}

\def\NN{{\cal N}}
\def\cga{\overset\circ{\ga}}
\def\cnab{\overset\circ{\nab}}

\def\tr{\mbox{tr}}

\def\L{{\cal L}}
\def\Lb{\underline{L}}

\def\gg{{\bf g}}

\def\gat{{\gamma_{t}}}

\def\RR{{\cal R}}

\renewcommand{\div}{\mbox{div }}

\def\nnab{\overline{\nab}}

\def\a{\alpha}
\def\b{\beta}
\def\ga{\gamma}
\def\Ga{\Gamma}
\def\de{\delta}
\def\De{\Delta}
\def\ep{\epsilon}
\def\eps{\epsilon}

\def\La{\Lambda}
\def\si{\sigma}
\def\Si{\Sigma}
\def\om{\omega}

\def\th{\theta}
\def\ze{\zeta}
\def\nab{\nabla}

\def\bb{\underline{\b}}
\def\lap{\Delta}

\newcommand{\trchb}{\tr \chib}
\def\trchav{\overline{\trch}}
\newcommand{\chih}{\hat{\chi}}
\newcommand{\chib}{\underline{\chi}}

\newcommand{\chibh}{\underline{\hat{\chi}}\,}
\def\f14{\frac{1}{4}}
\def\f12{{\frac{1}{2}}}

\def\c{\cdot}

\newcommand{\les}{\lesssim}

\newcommand{\HH}{{\mathcal H}}

\def\Gd{A}

\def\pr{\partial}

\def\chih{\hat{\chi}}
\def\trch{\mbox{tr}\chi}
\newcommand{\ddd}{\nab}

\newcommand{\dddd}{\overline{\ddd}}
\newcommand{\N}{{\cal N}}

\begin{document}
\theoremstyle{plain}
  \newtheorem{theorem}[subsection]{Theorem}
  \newtheorem{conjecture}[subsection]{Conjecture}
  \newtheorem{proposition}[subsection]{Proposition}
  \newtheorem{lemma}[subsection]{Lemma}
  \newtheorem{corollary}[subsection]{Corollary}

\theoremstyle{remark}
  \newtheorem{remark}[subsection]{Remark}
  \newtheorem{remarks}[subsection]{Remarks}

\theoremstyle{definition}
  \newtheorem{definition}[subsection]{Definition}

\include{psfig}
\title[rough Einstein metrics ]{ Sharp trace theorems  for null hypersurfaces on
Einstein metrics with finite curvature flux}
\author{Sergiu Klainerman}
\address{Department of Mathematics, Princeton University,
 Princeton NJ 08544}
\email{ seri@@math.princeton.edu}

\author{Igor Rodnianski}
\address{Department of Mathematics, Princeton University, 
Princeton NJ 08544}
\email{ irod@@math.princeton.edu}
\subjclass{35J10\newline\newline
The first author is partially supported by NSF grant 
DMS-0070696. The second author is a Clay Prize Fellow and is partially 
supported by NSF grant DMS-01007791}

\vspace{-0.3in}
\begin{abstract}The main objective of the paper is to prove  a geometric version
of sharp trace and product estimates 
 on null hypersurfaces  
with finite curvature flux. 
These estimates play a crucial role 
to control the geometry of such 
null hypersurfaces. The  paper 
 is based on  an invariant  version of the classical
Littlewood -Paley  theory, in a noncommutative setting,  defined via 
 heat flow on  surfaces.
\end{abstract}
\maketitle
\section{Introduction} To motivate the problems studied in this paper 
we start with the simplest example of a sharp trace theorem.
This applies to smooth  functions $f=f(t,x^1, x^2)$ on $I\times{\Bbb R}^2\subset{\Bbb R}^3$
 with $I$ an interval in
${\Bbb R}$, for simplicity $I=[0,1]$. We denote by  $\|\,\,\|_{L^2}$
the standard $L^2$ norm on $I\times{\Bbb R}^2$ and by   $\|\,\,\|_{H^k}$, $k$ positive integer, the usual 
  norm of the  Sobolev space $H^k(I\times{\Bbb R}^2)$ . Thus,
\beaa
\|f\|_{L^2}=\|f\|_{L^2(I\times{\Bbb R}^2)}&=&\big(\int_0^1\int_{{\Bbb R}^2}|f(t,x)|^2\, dt\,
dx\big)^\f12\\
\|f\|_{H^k}=\|f\|_{H^k(I\times{\Bbb R}^2)}&=&\big(\sum_{i+j\le k}\int_0^1\int_{{\Bbb
R}^2}|\pr_t^i\,\,\nab^j f(t,x)|^2 dt\, dx\big)^\f12
\eeaa
with $\nab^j$ the partial derivatives of order $j$ with respect to the $x$ variables.
We shall also use the mixed norm notation,
\beaa
\|f\|_{L_t^qL_x^p}&=&\big(\int_0^1\|f(t\,,\,\c)\|_{L^p({\Bbb R}^2)}^q dt\big)^{\frac{1}{q}}\\
\|f\|_{L_x^pL_t^q}&=&\big(\int_{{\Bbb R}^2}\|f(\c\,,\, x)\|_{L_t^q(I)}^p dx\big)^{\frac{1}{p}}
\eeaa
with the obvious modifications when either $p=\infty$ or $q=\infty$.
\begin{proposition}
The following  inequality holds for  an arbitrary, smooth,  
scalar function $f=f(t, x^1, x^2)$    in
${\Bbb R}^3$:
\be{eq:intr-1}
\|\pr_t f\|_{L_x^\infty L_t^2(I\times{\Bbb R}^2)}  \les \|f\|_{H^2(I\times{\Bbb R}^2)}
\end{equation}
\label{prop:intr-1}
\end{proposition}
The inequality can be easily derived with the help of the
 $W^{2,1}({\Bbb R}^2)\subset L^\infty$ imbedding and a standard 
integration by parts. Observe that the estimate is false if one replaces
$\pr_t$ with the other partial  derivatives $\pr_{x^1},\pr_{x^2}$.

Using the a standard Littlewood-Paley theory it 
is not too difficult to prove a stronger version
of \eqref{eq:intr-1} in Besov spaces\footnote{This is not just a  minor technical improvement.
It turns out that this type of Besov space  improvement of  the sharp
trace estimate plays a fundamental role in \cite{KR1}. }.
\be{eq:intr-2}
 \|\int_I |\pr_t f|^2 dt\|_{B^1_{2,1}({\Bbb R}^2)} \les \|f\|_{H^2(I\times{\Bbb R}^2)}^2.
\end{equation}
Here, for a function $g=g(x^1, x^2)$,
$$\|g\|_{B^\th_{2,1}({\Bbb R}^2)}=\sum_{k\ge 0} 2^{k\th}\|P_kg\|_{L^2({\Bbb R}^2)}+\|P_{<0} g\|_{L^2({\Bbb
R}^2)}$$ denotes  the standard, inhomogeneous,  Besov norm in ${\Bbb R}^2$ with   $\th\ge 0$
and  $P_k$
the usual Littlewood Paley (LP) projections, see  for example \cite{S1}, \cite{S2} and\cite{B}
for applications to paradifferential calculus. Also,
$P_{<0}=\sum_{k<0} P_k$. Observe that,
$$\|g\|_{B^1_{2,1}}\les \|\nab g\|_{B^0_{2,1}}+\|g\|_{L^2}$$
and therefore \eqref{eq:intr-2} follows easily
from  its following  bilinear version:
\begin{proposition}[ Sharp bilinear trace] The following inequality holds   for  an arbitrary, smooth,  
scalar function $g,h$    on $I\times{\Bbb R}^2$:
\be{eq:intr-3}
\|\int_I \pr_t g\c h\|_{B^0_{2,1}({\Bbb R}^2)}\les \|g\|_{
H^1(I\times{\Bbb R}^2)}\c \|h\|_{H^1(I\times{\Bbb R}^2)}
\end{equation}
\label{prop:intr-2}
\end{proposition}
In addition to the  bilinear sharp trace estimate \eqref{eq:intr-3}
 we  also signal the following related   estimate,
\begin{proposition}[ Sharp integrated product]
The following inequality holds   for  an arbitrary, smooth,  
scalar function $g,h$    on $I\times{\Bbb R}^2$:
\be{eq:intr-prod-trace}
\|\int_I  g\c h\|_{B^0_{2,1}({\Bbb R}^2)}\les \big(\|g\|_{
H^1(I\times{\Bbb R}^2)}+\|g\|_{L_x^\infty L_t^2}\big)\c \|h\|_{L_t^2 B_x^0}
\end{equation}
where,
$$ \|h\|_{L_t^q B_x^0}=\big(\int_0^1\|h(t,\c)\|_{B^0_{2,1}({\Bbb R}^2)}^qdt\big)^{1/q}.$$
\label{prop:intr-producttrace}
Also,
\be{eq:intr-prod-trace2}
\|  g\c  \int_0^t h\|_{L_t^2 B_x^0}\les \big(\|g\|_{
H^1(I\times{\Bbb R}^2)}+\|g\|_{L_x^\infty L_t^2}\big)\c \|h\|_{L_t^1 B_x^0}
\end{equation}
\end{proposition}
\begin{remark}Proposition \ref{prop:intr-producttrace}
is intimately related with the following well known
estimate for functions in  ${\Bbb R}^2$:
\be{eq:}
\|f\c g\|_{B^0_{2,1}({\Bbb R}^2)}\les \big(\|g\|_{ H^1{(\Bbb R}^2)}+\|g\|_{L^\infty ({\Bbb R}^2)}\big)\c
\|h\|_{B^0_{2,1}({\Bbb R}^2)}
\end{equation}
\end{remark}

Propositions \ref{prop:intr-1} and \ref{prop:intr-producttrace} can be easily reformulated
in terms of functions defined on null hypersurfaces in Minkowski space
 ${\Bbb R}^{3+1}$. For simplicity consider the standard null 
hypersurface  defined by
the equation $u=-1$ where $u$ is the optical function
$u=t-r$, $r=\sqrt{\sum_{i=1}^3( x^i)^2}$.  Let $\HH$ denote the portion
of this null hypersurface contained between $t=0$ and $t=1$. An arbitrary null
geodesic $\Ga=\Ga_{\om}$  along $\HH$ can be parametrized by $(t, t\om)$
where $\om$ is a unit vector in ${\Bbb R}^3$. Given a scalar function
$f$ on $\HH$ we denote by $\int_{\Ga_{\om}} f=\int_0^1f(t,t\om) dt.$ 
We   denote by    $\ddd_L f$ its  derivative along the null
geodesic, i.e. $\ddd_L f=\frac{d}{dt}f(t,t\om)$. We also denote by $\nab f$ the angular
derivatives of $f$ and by $\dddd f$ all tangential derivatives of $f$ along
$\HH$, i.e. $\dddd f=(\ddd_L f, \nab f)$. To adapt 
proposition  \ref{prop:intr-2} to the case of the  null hypersurface $\HH$
we need to define a Besov  space $B^0_{2,1}({\Bbb S}^2)$ 
 analogous to $B^0_{2,1}({\Bbb R}^2)$ with $ {\Bbb S}^2$ the standard 
unit sphere in ${\Bbb R}^3$.  LP projections can, of course,  be easily
defined locally, in coordinate charts, and extended to
 all of ${\Bbb S}^2$ by a partition of unity. Besov spaces 
on ${\Bbb S}^2$ can then be formally introduced as before.
 A more intrinsic way to define such spaces would be
based on spherical harmonic decomposition. Yet another 
way to achieve the same result is to introduce a definition of LP
projections based on a heat flow for the corresponding Laplace-Beltrami
operator  on the leaves of the geodesic sphere
foliation  of $\HH$ given by
the level surfaces of the standard time function $t$, see \eqref{eq:LP-introd}. This is in fact
the approach we  develop
 here  to deal with null hypersurfaces in non flat spacetimes.

The propositions below  are straightforward   adaptations
of propositions  \ref{prop:intr-2},  \ref{prop:intr-producttrace}
to the case of the  null hypersurface $\HH$. The  norms used
in the proposition are,
\beaa
\|f\|_{L^2(\HH)}&=&\big(\int_0^t \int|g(t,t\om)|^2dt d\om\big)^\f12\\
\NN_1(g)&=&\|\ddd_L f\|_{L^2(\HH)}+\|\nab f\|_{L^2(\HH)}+\| f\|_{L^2(\HH)}\\
\|f\|_{ L_x^\infty L_t^2}&=&\sup_{\om\in{\Bbb S}^2}\big(\int_{\Ga_\om}|f|^2\big)^\f12\\
\|f\|_{\PP^0}&=&\big(\int\|f(t, t\,\c)\|_{B^0_{2,1}({\Bbb S}^2)}^2dt\big)^\f12
\eeaa

\begin{proposition} Let $g,h$ be arbitrary smooth functions on $\HH$.
then the following estimates hold true, uniformly in $\om\in {\Bbb S}^2$,
\be{eq:intr-bilineartrace-cone}
\|\int_{\Ga_\om} (\ddd_L g\c h)\|_{B^0_{2,1}({\Bbb S}^2)}\les \NN_1(g)\c\NN_1(h)
\end{equation}
\be{eq:intr-producttrace-cone}
\|\int_{\Ga_\om}  g\c h\|_{B^0_{2,1}({\Bbb S}^2)}\les\big( \NN_1(g)+\|g\|_{L_x^\infty
L_t^2}\big)\c\|h\|_{\PP^0}
\end{equation}
\be{eq:intr-producttrace-cone2}
\| g\c \int_{\Ga_\om}  h\|_{\PP^0}\les\big( \NN_1(g)+\|g\|_{L_x^\infty
L_t^2}\big)\c\|h\|_{\PP^0}
\end{equation}
As a corollary of   \eqref{eq:intr-bilineartrace-cone}   
we have the standard sharp trace theorem,
\be{eq:intr-6}
\sup_{\om\in {\Bbb S}^2}\int_{\Ga_\om} |\ddd_L f|^2\les \NN_2(f)^2
\end{equation}
with,
$$
\NN_2(f)=\|\ddd_L^2 f\|_{L^2(\HH)}+\|\nab\ddd_L f\|_{L^2(\HH)}+\|\nab^2 f\|_{L^2(\HH)}+\NN_1(f)
$$
\label{prop:intr-bilintrace-cone}
\end{proposition}

The goal of this paper is to adapt the estimates 
\eqref{eq:intr-bilineartrace-cone},  \eqref{eq:intr-producttrace-cone}  and  \eqref{eq:intr-6}   to null
hypersurfaces in curved backgrounds, verifying the Einstein-Vacuum equations,
and verifying the {\sl bounded  curvature flux}(BCF) condition of \cite{KR1}.
These results, which formed the content of the main lemma
in \cite{KR1}, played a crucial role in that paper. We recall
that the main result of \cite{KR1} was to prove that the
  (BCF)  condition  suffices to control
 the local  geometry of null hypersurfaces. The (BCF) condition
comes naturally in connection with  the {\sl bounded $L^2$ curvature
} conjecture and the result of \cite{KR1} is a crucial ingredient
in the resolution of that conjecture.

To illustrate the techniques needed in our work
we give, in section 2,  a quick  proof of propositions \ref{prop:intr-2}
and \ref{prop:intr-producttrace}. They are  based on the 
properties of the standard, euclidean,  Littlewood -Paley(LP)  projections $P_k$
which we recall below. 

The LP-projections $(P_k)_{k\in \Bbb Z}$, acting 
on functions $g(x)$, $x\in {\Bbb R}^2$,   are defined  as Fourier
multipliers according to the formula,
$$(P_k f )\hat{\,\,}(\xi)=\chi(2^{-k}\xi) f\hat{\,\,}(\xi)$$
with $f\hat{\,\,}(\xi)$ denoting the Fourier transform of $f$
and $\chi(\xi)=\chi(|\xi|)$ a real  smooth test function supported in $\f12\le |\xi|\le 2$.
Moreover, for all $\xi\in {\Bbb R}^2\setminus{0}$,\,\, $\sum_k  \chi(2^{-k}\xi)=1$.
We denote $P_J=\sum_{k\in J} P_k$ for all intervals $J\subset {\Bbb Z}$.

The following properties are at the heart of the classical LP theory:

{\bf LP 1.}\quad  {\sl Almost Orthogonality:} \quad The operators $P_k$
are selfadjoint and verify
 $P_{k_1}P_{k_2}=0$   for all pairs of integers  such that  $|k_1-k_2|\ge 2$.
In particular,
$$\|f\|_{L^2}\approx\sum_k\|P_k f\|_{L^2}$$

{\bf LP 2.} \quad {\sl $L^p$-boundedness:} \quad For any $1\le
p\le \infty$, and any interval $J\subset \Bbb Z$,
\be{eq:pdf1'}
\|P_Jf\|_{L^p}\les \|f\|_{L^p}
\end{equation}

{\bf LP 3.} \quad {\sl Finite band property:}\quad  We can write any partial derivative
$\nab P_k f$ in the form  $\nab P_k f=2^k \tilde{P}_k f$ where $\tilde{P}_k$
are the  LP-projections associated with a slightly  different test function $\tilde{\chi}$ 
and  verify the property {\bf LP2}. Thus, in particular,
for any $1\le p\le \infty$
\beaa
\|\nab P_k f\|_{L^p}&\les& 2^{k}\|  f\|_{L^p}\\
2^k\|P_k f\|_{L^p}&\les& \|\nab f\|_{L^p}
\eeaa

 {\bf LP 4.}  \quad{\sl Bernstein inequalities}. \quad For any $2\le p\le \infty$
 we have the Bernstein inequality and its dual,
$$  \|P_k f\|_{L^\infty}\les 2^{k} \|f\|_{L^2},\qquad \|P_k f\|_{L^2}\les 2^{k} \|f\|_{L^1}$$

{\bf LP 5.}\quad {\sl Commutation properties}\quad  Given functions $f(t,x)$,
 if we denote by $P_kf$  the action of the LP projections in $x$,  we have, trivially
$$\pr_t P_k f =P_k \pr_t f, \qquad P_k\int_0^1 f dt =\int_0^1 P_k fdt$$

The proof of proposition  \ref{prop:intr-2}
is  a typical illustration of the power of paradifferential calculus. Clearly
one needs, somehow, to integrate by parts, but because of the symmetry of the estimate
 \eqref{eq:intr-3} with respect to $f, g$ one does not achieve anything by a direct 
integration by parts.
 The idea is to  decompose both functions $f=\sum_{k'} P_{k'} f$ and $g=\sum_{k''} P_{k''}
 g$  in  the bilinear expression $\int_0^1 \ddd_L f\c g$ and integrate by parts only those integral 
terms $\int_0^t(\pr_t f_{k'}\c g_{k''})$ where $k'\ge k''$. This illustrates a general
 philosophy; 
the paradifferential calculus  gives us the flexibility to deal differently with
various parts of nonlinear expressions and thus  allows us to  separate and focus on  various
difficulties of the problem at hand.  It is truly a {\sl divide
and conquer strategy.}

In order to extend the results mentioned below to
 general null hypersurfaces $\HH$  we need to replace the LP
theory based on Fourier transform with a more  intrinsic geometric 
definition. Given a Riemannian  manifold $M$  we can define LP
projections $P_k$
according to the formula,
 \be{eq:LP-introd}P_k F=\int_0^\infty m_k(\tau) U(\tau) F d\tau
\end{equation}
where $m_k(\tau)=2^{2k}m(2^{2k}\tau)$ and $m(\tau)$ is a Schwartz function 
with a finite number of  vanishing moments.
 The operator $U(\tau)F$ denotes the unique solution
of the heat flow on $S$ with initial data provided by $F$,
$$\pr_\tau U(\tau)F -\lap U(\tau)F=0, \,\quad U(0)F=F.$$
where  $\lap$ denotes the standard Laplace-Beltrami operator for tensors,
$\lap =\ga^{ij}\nab_i\nab_j .$
We apply this definition to the $2$-dimensional leaves of the geodesic foliation on 
our null hypersurface $\HH$. Under some simple assumptions on the geometry
of these leaves we  prove, in \cite{KR2},  a sequence of properties of the LP 
projections similar to {\bf LP1}, {\bf LP4}. Some of our results are, of
course, weaker. For example the pointwise version of the almost orthogonality property
{\bf LP5} cannot possible be true. We can replace it however by a sufficiently
robust $L^p$  analogue of it. We also  find satisfactory  analogues
for {\bf LP2}-{\bf LP3}, though we have to be quite careful about 
what we can in fact prove with  our very limited regularity
assumptions. For example, we can prove    a version of the Bernstein inequalities
of {\bf LP4} for scalars $f$  but not for   tensorfields. Of course,
{\bf LP5} plays a fundamental role in the proof of propositions
\ref{prop:intr-2}, \ref{prop:intr-producttrace}. Such a property, however,
does not  hold for the  nonflat backgrounds we deal with in our 
work.  This lack of commutativity compounded  by the weak regularity
properties of the foliation, consistent with the (BCF) assumption,  leads to 
considerable conceptual and technical difficulties.
   
Once we have set up a satisfactory geometric LP theory
we can formulate and prove results, on non flat backgrounds,
similar to those 
of proposition \ref{prop:intr-bilintrace-cone}. 
We prefer to state these results 
only after a thorough discussion of 
the geometric framework and 
the properties of the intrinsic LP
projections.

We give complete proofs of both  propositions \ref{prop:intr-2}, \ref{prop:intr-producttrace}
in section 2 in order to prepare the reader for the methods used in the non flat situation.

In section 3 we discuss the main geometric notions  concerning null
hypersurfaces in a curved background. We also introduce our main assumptions 
{\bf BA1}, {\bf BA2}, {\bf WS},  {\bf K1} and {\bf K2}. All these assumptions
are consistent with the bounded curvature flux (BCF) condition  of \cite{KR1}.
The main results of this paper depend only on these   assumptions.
We also introduce our geometric  LP projections and state
their main properties, proved in \cite{KR2}. Finally we define our main Besov spaces
 and recall some of the properties proved in \cite{KR2}.

In section 4 we state a sequence of theorems
which  extend the results of \eqref{prop:intr-bilintrace-cone} to nonflat
 backgrounds.
These results were stated in section 5 of \cite{KR1},
and played an essential role in the proof of the main result there.

The remaining sections of the paper contain the proofs of these results.
\section{Proof of propositions \ref{prop:intr-2} and \ref{prop:intr-producttrace}} 
We shall  make use of the properties {\bf P1}--{\bf P5} of the classical
LP projections mentioned in the introduction.
\subsection{Proof of proposition \ref{prop:intr-2}}
By definition, 
\beaa
\|\int_0^1\pr_t g\c h dt\|_{B_{2,1}^0}&=&\sum_{k\ge 0} 
\|P_k \int_0^1\pr_t g\c h dt\|_{L_x^2}+\|P_{<0} \int_0^1\pr_t g\c h dt\|_{L_x^2}\\
&\les&\sum_{k} 
\|P_k \int_0^1\pr_t g\c h dt\|_{L_x^2}
\eeaa
We decompose, with respect to the $x$ variables,
  $f=\sum_k f_k$, $g=\sum_k g_k$ with $f_k=P_k f,\,\, g_k=P_k g$,   
 and write,
\beaa
 P_k\int_0^1(\pr_t g\c h)=A_k+B_k+C_k+D_k
\eeaa
\bea
A_k&=& P_k\int_0^1  (\pr_t g)_{<k}\c h_{\ge k},\qquad\quad
B_k=P_k\int_0^1  (\pr_t g)_{\ge k}\c h_{< k}\label{eq:ABCDkflat}\\
C_k&=&P_k\int_0^1  (\pr_t g)_{<k}\c h_{< k},\qquad\quad
D_k=P_k  \int_0^1(\pr_t g)_{\ge k}\c h_{\ge  k}\nn
\eea
Observe that  $P_k$ commute with the integral $\int_0^1$
and that $C_k$ is essentially\footnote{With the possible exception of 
a finite, $\,\les 8\,$,  number of terms  which
can be made part of  either $A_k, B_k, D_k$. 
This corresponds to the classical 
trichotomy formula and is due to {\bf LP1.)}} zero. 
Thus we only have to estimate $A_k, B_k, D_k$.

{\bf 1.)}\quad  {\sl Estimates for $A_k=P_k\int_0^t (\pr_t g)_{<k}\c h_{\ge k}$:} 
\qquad 
We use first {\bf LP2}  followed by the 
the direct Bernstein inequality  {\bf LP4}, and then {\bf LP3},
{\bf LP2}. We also use {\bf LP5} freely.
\beaa
\|A_k\|_{ L^2_x}&\les&\sum_{k'<k\le k''}\int_0^{1}
\|(\pr_t g)_{k'}\c h_{k''}\|_{L^2_x}\, dt
\les\sum_{k'<k\le k''}\int_0^{1}
\|(\pr_t g)_{k'}\|_{L^\infty_x} \c\|h_{k''}\|_{L^2_x}\,dt\\
&\les&\sum_{k'<k\le k''}2^{k'}\int_0^{1}
\|(\pr_t g)_{k'}\|_{L^2_x}\c\| h_{k''}\|_{L^2_x}\,dt\\
&\les& \sum_{k'<k\le k''}2^{k'-k''}
\|(\pr_t g)_{k'}\|_{L^2_t L^2_x}\c\| \nab h_{k''}\|_{L^2_t L^2_x}
\eeaa
Therefore, writing $2^{k'-k''}=2^{(k'-k)/2+(k'-k'')/2}\les 2^{(k'-k)/2+(k-k'')/2}$,
summing over $k$ and using {\bf LP1},
\beaa
\sum_{k\ge 0} \|A_k\|_{L_t^\infty L^2_x}&\les&\sum_k \sum_{k'<k\le k''}2^{(k'-k'')/2}\,\,2^{(k-k'')/2}
\|(\pr_t g)_{k'}\|_{L^2_t L^2_x}\| \nab h_{k''}\|_{L^2_t L^2_x}\\
&\les & \sum_{k'< k''}2^{(k'-k'')/2}
\|(\pr_t g)_{k'}\|_{L^2_t L^2_x}\| \nab h_{k''}\|_{L^2_t L^2_x}\\&\les &
\|\pr_t g\|_{L_t^2L_x^2}\c  \|\nab h\|_{L_t^2L_x^2}  .
\eeaa
{\bf 2.)}\quad   {\sl Estimates for $D_k=P_k\int_0^t (\pr_t F)_{\ge k}\c G_{\ge  k}$:}
\qquad 
We write, $D_k=D_k^1+D_k^2$  where,
\begin{align*}
D_{k}^{1}=\sum_{k\le k'\le k''} P_{k} \int_0^1 (\pr_t g)_{k'}
\c h_{k''},\qquad
D_{k}^{2}= \sum_{k\le k'<k''} P_{k} \int_0^1(\pr_t g)_{k''}\c h_{k'}
\end{align*}
The term  $D_k^1$ can be treated in a straightforward manner,
without integration by parts. We start by using the dual Bernstein inequality
{\bf LP4}, followed by the finite band property {\bf LP3},
\beaa
\|D_k^1\|_{L_x^2}&\les& 2^k\sum_{k\le k'\le k''}\| \int_0^1 (\pr_t g)_{k'}
\c h_{k''}\|_{L_x^1}\les 2^k\sum_{k\le k'\le k''}
\|\pr_t g_{k'}\|_{L_t^2L_x^2}\c\|h_{k''}\|_{L_t^2L_x^2}\\
&\les&2^{k-k''}\|\pr_t g_{k'}\|_{L_t^2L_x^2}\c\|\nab h_{k''}\|_{L_t^2L_x^2}
\eeaa
Thus, summing in $k$  and using the $L^2$ orthogonality property  {\bf LP1}
\beaa
\sum_{k\ge 0}\|D_k^1\|_{L_x^2}&\les& \sum_k \sum_{k\le k'<k''} 2^{k-k''}\|\pr_t
g_{k'}\|_{L_t^2L_x^2}\c\|\nab h_{k''}\|_{L_t^2L_x^2}\\
&\les&   \sum_{ k'<k''} 2^{k'-k''}\|\pr_t
g_{k'}\|_{L_t^2L_x^2}\c\|\nab h_{k''}\|_{L_t^2L_x^2}\les  \|\pr_t g\|_{L_t^2L_x^2}\c \|\nab
h\|_{L_t^2L_x^2}
\eeaa
To treat $D_{k}^{2}$ we need to transfer the $\pr_t$ derivative 
from the low frequency term $ g_{k'}$ to the high frequency 
term $h_{k''}$. After integration by parts  we can treat the resulting
integral exactly as $D_k^1$, the only terms
we need to worry about are the boundary terms 
 $\|I_k(1)-I_k(0)\|_{L_x^2}\les \sup_{0\le t\le 1}\|I_k(t,\c)\|_{L_x^2}=\|I_k\|_{L_t^\infty L_x^2}$,
where
 $$I_k =\sum_{k\le k'<k''}  P_k \big( g_{k''}\c h_{k'}\big).$$
To treat them we need the following 
\begin{lemma} For any $k, k', k''$, we have
\beaa
 \|P_k \big(g_{k'}\c h_{k''}\big)\|_{L_t^\infty L_x^2}\les 
2^{-\frac{1}{4}\big(|k'-k|+|k''-k|\big)}\|g_{k'}\|_{H^1}\c \|
h_{k''}\|_{H^1}
\eeaa
\label{le:integr-parts-k}
\end{lemma}
Using the lemma
we derive,
\beaa
\sum_k\|D_k^2\|_{L_x^2}&\les& \sum_k \sum_{k\le k'<k''} 2^{k-k''}\|\pr_t g_{k''}\|_{L_t^2L_x^2}\c\|\nab
h_{k'}\|_{L_t^2L_x^2}\\
&+& \sum_k \sum_{k\le k'<k''}2^{-\f12\big(|k'-k|+|k''-k|\big)}\| g_{k''}\|_{H^1}\c \|
h_{k'}\|_{H^1}\\
&\le&\sum_{ k'<k''}2^{-\f12|k'-k''|}\| g_{k''}\|_{H^1}\c \|
h_{k'}\|_{H^1}\les \| g\|_{H^1}\c\| h\|_{H^1}
\eeaa
as desired.

{3.)}\quad {\sl Estimates for 
$B_k=P_k \int_0^t (\pr_t g)_{\ge k}\c h_{< k} $:}\qquad 
We start by decomposing,
\beaa
B_k=\sum_{k'<k\le k''}P_k\int_0^t (\pr_t g)_{k''}\c h_{k'}
\eeaa
Integrating by parts and using lemma \ref{le:integr-parts-k}
to estimate the boundary terms, we derive,
we obtain, with $J_k=\sum_{k'< k\le k''}  P_k \big( g_{k''}\c h_{k'}\big)$,
\beaa
\|B_k\|_{ L^2_x}&\les& 
\sum_{k'<k\le k''}
\|\int_0^1   g_{k''}\c (\pr_t h)_{k'}\|_{ L^{2}_x} +\|J_k\|_{L_t^\infty L_x^2}\\
\|J_k\|_{L_t^\infty L_x^2}&\les&\sum_{k'<k\le k''}
2^{-\f12\big( |k'-k|+|k''-k|\big)}\| g_{k''}\|_{H^1}\c \|
h_{k'}\|_{H^1}
\eeaa
Now,
\beaa
\sum_{k'<k\le k''}\|\int_0^1   g_{k''}\c (\pr_t h)_{k'}\|_{ L^{2}_x}&\les &
\sum_{k'<k\le k''}\| g_{k''}\|_{L_t^2L_x^2}\c \|(\pr_t h)_{k'}\|_{L_t^2 L_x^\infty}\\
&\les&\sum_{k'<k\le k''}2^{k'-k''}\|\nab g_{k''}\|_{L_t^2L_x^2}\c \|(\pr_t h)_{k'}\|_{L_t^2 L_x^2}
\eeaa
Thus, summing  as before,

\beaa
\sum_k\|B_k\|_{ L^2_x}&\les&\sum_k\sum_{k'<k\le k''}2^{k'-k''}\|\nab g_{k''}\|_{L_t^2L_x^2}\c \|(\pr_t
h)_{k'}\|_{L_t^2 L_x^2}\\
&+&\sum_k\sum_{k'<k\le k''}
2^{-\f12\big( |k'-k|+|k''-k|\big)}\| g_{k''}\|_{H^1}\c \|
h_{k'}\|_{H^1}
\\
&\les& \sum_{k'< k''}2^{-\f12 |k'-k''|}\| g_{k''}\|_{H^1}\c \|
h_{k'}\|_{H^1}
\les \| g\|_{H^1}\c \|
h\|_{H^1}
\eeaa
as desired.

\begin{proof} {\bf of lemma \ref{le:integr-parts-k}}:\quad 
By symmetry it suffices to consider the following cases:
$$k'\ge k''\ge k,  \quad k'\ge k>k'', \quad k>k'\ge k'' $$
Observe that the last case, which we call ``low-low"  interaction, cannot occur. 
We are thus left with
only two cases:

{\bf Case 1.}\quad $k'\ge k''\ge k$
\beaa
\|P_k \big ( g_{k'}\c h_{k''}\big )\|_{L_t^\infty L_x^2} &\les & 
2^k \|g_{k'}\c h_{k''}\|_{L_t^\infty L_x^1}
\les  2^k \|g_{k'}\|_{L^\infty_t L^2_x}\| h_{k''}\|_{L_t^\infty L_x^2}
\eeaa
We now make use 
of the simple calculus inequality,
$$\|f\|_{L_t^\infty L_x^2}\les \|\pr_t f\|_{L_t^2L_x^2}^\f12\c  \|
f\|_{L_t^2L_x^2}^\f12+\|f\|_{L_t^2 L_x^2}.$$ Thus, 
\beaa
\|g_{k'}\|_{L^\infty_t L^2_x}&\les& \|\pr_t g_{k'}\|_{L_t^2L_x^2}^\f12\c 
 \| g_{k'}\|_{L_t^2L_x^2}^\f12+\|g_{k'}\|_{L_t^2L_x^2}\\
&\les&
 2^{k'/2}\big( \|\pr_t g_{k'}\|_{L_t^2L_x^2}^\f12\c 
 \| \nab g_{k'}\|_{L_t^2L_x^2}^\f12+\| g_{k'}\|_{L_t^2L_x^2}^\f12\c 
 \| \nab g_{k'}\|_{L_t^2L_x^2}^\f12\big)\\
&\les&2^{k'/2}\| \nab g_{k'}\|_{L_t^2L_x^2}^\f12\big(\|\pr_t g_{k'}\|_{L_t^2L_x^2}
+\| g_{k'}\|_{L_t^2L_x^2}\big)^\f12\\
&\les&2^{k'/2}\| g_{k'}\|_{H^1}^\f12\c \| g_{k'}\|_{H^1}^\f12=2^{k'/2}\| g_{k'}\|_{H^1}
\eeaa
Thus,
\bea
   \|g_{k'}\|_{L^\infty_t L^2_x}  &\les&    2^{k'/2}\| g_{k'}\|_{H^1},\qquad
\|h_{k''}\|_{L^\infty_t L^2_x}\les 2^{k''/2}\| h_{k''}\|_{H^1}\label{eq:funnyBernstein}
\eea
Therefore, since $k'\ge k''\ge k$,
\beaa
\|P_k \big ( g_{k'}\c h_{k''}\big )\|_{L_t^\infty L_x^2} &\les & 2^{k-\frac {k'}2-\frac {k''}2}
\| g_{k'}\|_{H^1}\c \| h_{k''}\|_{H^1}\\
&=& 2^{-\f12\big(|k'-k|+|k''-k|\big)}\| g_{k'}\|_{H^1}\c \| h_{k''}\|_{H^1}
\eeaa
as desired.

{\bf Case 2.}\quad $k'\ge k>k''$\quad Using once more the estimates \eqref{eq:funnyBernstein},
\beaa
\|P_k \big(g_{k'}\c h_{k''}\big)\|_{L_t^\infty L_x^2}&\les&
 \|g_{k'}\|_{L_t^\infty L_x^2}\c\|h_{k''}\|_{L_t^\infty L_x^\infty}
\les 2^{k''}\|g_{k'}\|_{L_t^\infty L_x^2}\c \|h_{k''}\|_{L_t^\infty L_x^2}\\
&\les& 2^{k''} 2^{-\frac{k'}{2}-\frac{k''}{2}}\| g_{k'}\|_{H^1}\c \| h_{k''}\|_{H^1}= 
2^{(k''-k')/2}\| g_{k'}\|_{H^1}\c \| h_{k''}\|_{H^1}\\
&=&2^{-\f12\big(|k'-k|+|k''-k|\big)}\| g_{k'}\|_{H^1}\c \| h_{k''}\|_{H^1}
\eeaa
as desired.
\end{proof}
\subsection{Proof of proposition \ref{prop:intr-producttrace}}\quad

\begin{proof} {\bf of}  \eqref{eq:intr-prod-trace}:\quad 
As before we have to estimate the sum,
\beaa
\sum_{k} 
\|P_k \int_0^1 g\c h \,dt\|_{L_x^2}
\eeaa
For each integer $k$ we  decompose, $h=h_{<k}+h_{\ge k}$
with $h_{<k}=\sum_{k'<k} h_{k'}$, and  $\,\, h_{\ge k}=\sum_{k'\ge k} h_{k'}$.
Thus,
\beaa
P_k \int_0^1 g\c h \,dt&=&P_k\big(\int_0^1 g\c h_{\ge k}\big)+P_k\big(\int_0^1 g\c h_{< k}\big)=
A_k+ B_k
\eeaa

 {\bf   1.)}\quad {\sl Estimates for $B_k$.}\quad
Observe that\footnote{In fact $ P_k(g\c h_{<k})=P_k(g_{> k+1}\c h_{< k})$, neglecting
a finite number of terms does not matter. }   $ P_k(g\c h_{< k})=P_k(g_{\ge k }\c h_{< k})$.
Thus,  using freely {\bf LP5}, 
$B_k=P_k\big(\int_0^1 g_{\ge k}\c h_{< k}\big)$. 
We now rely on  the  Bernstein and finite band  inequalities {\bf LP3 }, {\bf LP4},
\beaa
\|B_k\|_{L_x^2}&\les&  \|\int_0^1 g_{\ge k}\c h_{< k}\|_{L_x^2}
\le \| g_{\ge k}\|_{L_t^2L_x^2}\c\| h_{< k}\|_{L_t^2L_x^2}\\
&\les&2^k\sum_{k''<k\le k'}\| g_{ k'}\|_{L_t^2L_x^2}\c\| h_{k''}\|_{L_t^2L_x^\infty}\\
&\les&\sum_{k''<k\le k'}2^{k''-k'}\|\nab g_{ k'}\|_{L_t^2L_x^2}\c\| h_{k''}\|_{L_t^2L_x^2}
\eeaa
Therefore,
\beaa
\sum_k \|B_k\|_{L_x^2}&\les&\sum_k \sum_{k''<k\le k'}2^{k''-k'}\|\nab g_{ k'}\|_{L_t^2L_x^2}
\c\|h_{k''}\|_{L_t^2L_x^2}\\
&\les&\sum_{k''\le k'}2^{\frac{k''-k'}{2}}\|\nab g_{
k'}\|_{L_t^2L_x^2}\c\|h_{k''}\|_{L_t^2L_x^2}
\les\|\nab g\|_{L_t^2L_x^2}\c \| h\|_{L_t^2L_x^2}.
\eeaa

 {\bf   2.)}\quad {\sl Estimates for $A_k$.}\quad
To estimate $A_k$ we have to be more careful. According \footnote{In fact $\lap P_k f= 2^{2k}\tilde {P}_k
f$  with a slightly modified LP -projection.} to {\bf LP3},
$\lap P_k f= 2^{2k}P_k f$
\beaa
\|A_k\|_{L_x^2}&=&\|P_k\big(\int_0^1 g\c\lap h_{\ge k}\big)\|_{L_x^2}\les \sum_{k'\ge k}
2^{-2k'}\|P_k\big(\int_0^1 g\c \lap h_{ k'}\big)\|_{L_x^2}\\
&\les& \sum_{k'\ge k}
2^{-2k'}\|P_k\big(\nab\int_0^1  (g\c \nab h_{ k'})\big)\|_{L_x^2}
+\sum_{k'\ge k}2^{-2k'}\|P_k\big(\int_0^1  \nab g\c
\nab h_{ k'}\big)\|_{L_x^2}\\
&=& J_1+J_2
\eeaa
Now, using {\bf LP3},
\beaa
J_1&=& \sum_{k'\ge k}
2^{-2k'}\|P_k\big(\nab\int_0^1  (g\c \nab h_{ k'})\big)\|_{L_x^2}
\les \sum_{k'\ge k}
2^{-2k'+k}\|\int_0^1  g\c \nab h_{ k'}\|_{L_x^2}\\
&\les&\sum_{k'\ge k}
2^{-2k'+k}\|\big(\int_0^1|g|^2\big)^\f12\c\big(\int_0^1|\nab h_{ k'}|^2\big)^\f12\|_{L_x^2}\\
&\les&\sum_{k'\ge k}
2^{-2k'+k}\|g\|_{L_x^\infty L_t^2}\c \|\nab h_{ k'}\|_{L_t^2L_x^2}
\les \sum_{k'\ge k} 2^{-k'+k}\|g\|_{L_x^\infty L_t^2}\c \| h_{ k'}\|_{L_t^2L_x^2}
\eeaa
On the other hand, using the dual Bernstein inequality
\beaa
J_2&=& \sum_{k'\ge k}2^{-2k'}\|P_k\big(\int_0^1  \nab g\c
\nab h_{ k'}\big)\|_{L_x^2}
\les\sum_{k'\ge k}2^{-2k'+k}\|\big(\int_0^1  \nab g\c
\nab h_{ k'}\big)\|_{L_x^1}\\
&\les&\sum_{k'\ge k}2^{-2k'+k}\|\nab g\|_{L_t^2L_x^2}\c \|\nab  h_{ k'}\|_{L_t^2L_x^2}
\les\sum_{k'\ge k}2^{-k'+k}\|\nab g\|_{L_t^2L_x^2}\c \|  h_{ k'}\|_{L_t^2L_x^2}
\eeaa
Hence,
\beaa
\|A_k\|_{L_x^2}\les \sum_{k'\ge k}2^{-k'+k} \|  h_{ k'}\|_{L_t^2L_x^2}\c
\big(\|\nab g\|_{L_t^2L_x^2}+\|g\|_{L_x^\infty L_t^2}\big)
\eeaa
and therefore,
\beaa
\sum_k\|A_k\|_{L_x^2}&\les&\sum_{k'} \|  h_{ k'}\|_{L_t^2L_x^2} \c\big(\|\nab
g\|_{L_t^2L_x^2}+\|g\|_{L_x^\infty L_t^2}\big)\\
&\les&\big(\|\nab
g\|_{L_t^2L_x^2}+\|g\|_{L_x^\infty L_t^2}\big)\c \|h\|_{L_t^2 B_x^0}
\eeaa
 as desired.
\end{proof}
\begin{proof} {\bf of } 
 \eqref{eq:intr-prod-trace2}:\quad 
 We have to estimate the sum,
\beaa
\sum_{k} 
\|P_k\big(  g\c \int_0^t h \big)\,\|_{L_t^2 L_x^2}
\eeaa
For each integer $k$ we  decompose, $h=h_{<k}+h_{\ge k}$
with $h_{<k}=\sum_{k'<k} h_{k'}$, and  $\,\, h_{\ge k}=\sum_{k'\ge k} h_{k'}$.
Thus,
\beaa
P_k \big( g\c \int_0^t h\big) \,dt&=&
P_k\big(g\c \int_0^t h_{\ge k}\big)+P_k\big( g\c
\int_0^t h_{< k}\big)= A_k+ B_k
\eeaa

 {\bf   1.)}\quad {\sl Estimates for $B_k$.}\quad
Observe that\footnote{In fact $ P_k(g\c h_{<k})=P_k(g_{> k+1}\c h_{< k})$, neglecting
a finite number of terms does not matter. }   $ P_k(g\c \int_0^t h_{< k})=P_k(g_{\ge k }\c
\int_0^t h_{< k})$. Thus,  using freely {\bf LP5}, 
$B_k=P_k\big( g_{\ge k}\c \int_0^t h_{< k}\big)$. 
We  rely on   {\bf LP3 }, {\bf LP4},
\beaa
\|B_k\|_{L_t^2 L_x^2}&\les&  \| g_{\ge k}\c \int_0^t h_{< k}\|_{L_t^2L_x^2}
\le \| g_{\ge k}\|_{L_t^2L_x^2}\c\|\int_0^t h_{< k}\|_{L_t^\infty L_x^2}\\
&\les&2^k\sum_{k''<k\le k'}\| g_{ k'}\|_{L_t^2L_x^2}\c\| h_{k''}\|_{L_t^1L_x^\infty}\\
&\les&\sum_{k''<k\le k'}2^{k''-k'}\|\nab g_{ k'}\|_{L_t^2L_x^2}\c
\| h_{k''}\|_{L_t^1L_x^2}
\eeaa
Therefore,
\beaa
\sum_k \|B_k\|_{ L_t^2 L_x^2}&\les&\sum_k \sum_{k''<k\le k'}2^{k''-k'}\|\nab g_{
k'}\|_{L_t^2L_x^2}
\c\|h_{k''}\|_{L_t^1L_x^2}\\
&\les&\sum_{k''\le k'}2^{\frac{k''-k'}{2}}\|\nab g_{
k'}\|_{L_t^2L_x^2}\c\|h_{k''}\|_{L_t^1L_x^2}
\les\|\nab g\|_{L_t^2L_x^2}\c \| h\|_{L_t^1L_x^2}.
\eeaa

 {\bf   2.)}\quad {\sl Estimates for $A_k$.}\quad
To estimate $A_k$ we  proceed as in the proof of 
 \eqref{eq:intr-prod-trace}.
\beaa
\|A_k\|_{L_t^2L_x^2}&=&\|P_k\big( g\c
\int_0^t h_{\ge k}\big)\|_{L_t^2L_x^2}\les
\sum_{k'\ge k} 2^{-2k'}\|P_k\big( g\c \int_0^t\lap h_{ k'}\big)\|_{L_t^2 L_x^2}\\
&\les& \sum_{k'\ge k}
2^{-2k'}\|P_k\big(\nab  (g\c  \int_0^t \nab h_{ k'})\big)\|_{L_t^2 L_x^2}\\
&
+&\sum_{k'\ge k}2^{-2k'}\|P_k\big(  \nab g\c
 \int_0^t  \nab  h_{ k'}\big)\|_{L_t^2L_x^2}\\
&=& J_1+J_2
\eeaa
Now, using {\bf LP3},
\beaa
J_1&=& \sum_{k'\ge k}
2^{-2k'}\|P_k\big(\nab  (g\c \int_0^t \nab h_{ k'})\big)\|_{L_t^2 L_x^2}
\les \sum_{k'\ge k}
2^{-2k'+k}\| g\c \int_0^t \nab h_{ k'}\|_{L_t^2 L_x^2}
\\ &\les&\sum_{k'\ge k}
2^{-2k'+k}\|g\|_{L_x^\infty L_t^2}\c \|\int_0^t \nab h_{ k'}\|_{ L_x^2 L_t^\infty}
\les \sum_{k'\ge k} 2^{-k'+k}\|g\|_{L_x^\infty L_t^2}\c \| h_{ k'}\|_{L_t^1 L_x^2}
\eeaa
On the other hand, using the dual Bernstein inequality
\beaa
J_2&=& \sum_{k'\ge k}2^{-2k'}\|P_k\big( \nab g\c
\int_0^t \nab h_{ k'}\big)\|_{L_t^2 L_x^2}
\les\sum_{k'\ge k}2^{-2k'+k}\|  \nab g\c
\int_0^t \nab h_{ k'}\|_{L_t^2 L_x^1}\\
&\les&\sum_{k'\ge k}2^{-2k'+k}\|\nab g\|_{L_t^2L_x^2}\c \|\int_0^t \nab  h_{
k'}\|_{L_t^\infty L_x^2}\\
&\les&\sum_{k'\ge k}2^{-k'+k}\|\nab g\|_{L_t^2L_x^2}\c \|  h_{ k'}\|_{L_t^1 L_x^2}
\eeaa
Hence,
\beaa
\|A_k\|_{L_t^2 L_x^2}\les \sum_{k'\ge k}2^{-k'+k} \|  h_{ k'}\|_{L_t^1L_x^2}\c
\big(\|\nab g\|_{L_t^2L_x^2}+\|g\|_{L_x^\infty L_t^2}\big)
\eeaa
and therefore,
\beaa
\sum_k\|A_k\|_{L_t^2 L_x^2}&\les&\sum_{k'} \|  h_{ k'}\|_{L_t^1L_x^2} \c\big(\|\nab
g\|_{L_t^2L_x^2}+\|g\|_{L_x^\infty L_t^2}\big)\\
&\les&\big(\|\nab
g\|_{L_t^2L_x^2}+\|g\|_{L_x^\infty L_t^2}\big)\c \|h\|_{L_t^1 B_x^0}
\eeaa
 as desired.
\end{proof}

\section{Geometric set-up. Geometric LP- projections }
We assume given an Einstein spacetime $(\MM,\gg)$ a space -like hypersurface $\Si$ and
 an outgoing  null
hypersurface
$\HH$, initiating on a compact $2$ surface $S_0\subset\Si$ diffeomorphic to ${\Bbb S}^2$,   given by the level hypersurfaces 
of an optical function $u$, i.e. solution to the Eikonal equation 
\be{eq:I1}
\gg^{\a\b}\pr_\a u\pr_\b u=0.
\end{equation}
We briefly recall the main geometric definitions, see section 2 of \cite{KR1},  associated with
$\HH$.

1.) {\sl Geodesic foliation:}\quad Let $L=-\gg^{\a\b}\pr_a u\, \pr_b $ be the corresponding null generator vectorfield and $s$ 
its affine parameter, i.e. $L(s)=1,\quad s|_{S_0}=0$. The level surfaces $S_s$ of $s$
generates the geodesic foliation on $\HH$. We shall denote
by $\nab$ the covariant differentiation
on $S_s$ and by $\ddd_L$ the projection to $S_s$ of the
covariant derivative with respect to $L$, see section 2. We also denote by $r$ the
function on $\HH$ defined by  $r=r(s)=\sqrt{(4\pi)^{-1}|S_s|}   $, with $|S_s|$ the area of $S_s$. Let
$\HH_t$ be the portion of $\HH$ between $s=0$ and $s=t$ and, for simplicity, assume $\HH=\HH_1$.

2.) {\sl Null pair:} \quad  With our choice of $L$ we have 
 $$<L,L>=0,\qquad D_LL=0$$
where $<\,,\,>= <\,,\,>_g$ denoting the metric  of $\MM$.
At any point $P\in S_s\subset \HH$  we denote by 
$\Lb$ the null vector conjugate to $L$ relative
to the geodesic  foliation, i.e. 
$$<L\,,\,\Lb>=-2\,,\qquad <\Lb, X>=0 \quad \mbox{for all } \, X\in T_p(S_s).$$
We shall say that $L,\Lb$ form the canonical null  pair associated
to the foliation. We  denote by $\ga$ the induced metric on $S_v$, by $\nab$ the induced
covariant derivative and $K$ the Gauss curvature. An arbitrary orthonormal frame
on $S_v$ will be denoted by $(e_a)_{a=1,2}$. Clearly,
$$<e_a\,,\, L>=<e_a\,,\,\Lb>=0, \quad <e_a,e_b>=\de_{ab}.$$
A null pair together with an orthonormal frame $(e_a)_{a=1,2}$
as above is called a null 
frame associated to the foliation.

3.) {\sl Total curvature flux:}\quad  We introduce the  total curvature flux along
$\HH$ to be the integral\footnote{ Here $
\|F\|_{L^2(\HH)}=\big(\int_0^{1} ds\int_{S_s}  |F|^2 \big)^\f12$
where $ \int_{S_s}  |F|^2$
 denotes the integral with respect to  the
volume element $d\mu_s$ of  $S_s$. \label{footnote1}}, see precise definition 
 of the null curvature components $\a,\b,\rho,\si,\bb$ in section 2 of
\cite{KR1}, 
\be{eq:I2}
\RR_0=\bigg(\|\a\|_{L^2(\HH)}^2+\|\b\|_{L^2(\HH)}^2+
\|\rho\|_{L^2(\HH)}^2+\|\si\|_{L^2(\HH)}^2+\|\bb\|_{L^2(\HH)}^2\bigg)^\f12
\end{equation} 
with $\a,\b,\rho,\si,\bb$ null   components of the curvature tensor $\, \rr$
 of the spacetime, background metric $\bf g$.

In \cite{KR1} we worked under the assumption the $\RR_0$ is sufficiently small.
In this paper the only curvature components we shall need are $\b$ and the Gauss
curvature $K$ of the $S_s$ surfaces. We will  make specific assumptions about this,
consistent with the small curvature flux condition,  i.e. $\RR_0$ sufficiently small.

4.) {\sl Null connection coefficients:}\quad The null second fundamental forms $\chi,\chib$
 of the foliation $S_v$ are given by 
\be{eq:nullsforms}
\chi_{ab}=<D_aL\,,\, e_b>,\qquad \chib_{ab}=<D_a\Lb\,,\, e_b>
\end{equation} 
The torsion is given by,
\be{eq:torsion}
\ze_a=\f12 <D_a L\,,\,\Lb>
\end{equation}
We also denote  $\trch=\de^{ab}\chi_{ab}$ and $\trchb=\de^{ab}\chibh_{ab}$ and
$\chih=\chih-\f12\trch\de$, $\chibh=\chib-\f12\trchb\de$. Recall
 the definition of the mass aspect function $\mu=-\div\ze+\f12\chih\c\chibh -\rho+|\ze|^2$.

5.)\quad {\sl Commutator formulas}\quad Commutation formulas between 
$\ddd_L$ and $\nab $ play an important role in the paper. We recall, see section 2.15  of 
\cite{KR1},
\begin{proposition}
Consider an arbitrary k-covariant,
 S-tangent vectorfield $F_{\underline{a}}=F_{a_1\ldots a_k}$
Then,
\bea
\label{eq:commutationLnab}
\ddd_L\nab_b F_{\underline{a}}-\nab_b\ddd_L F_{\underline{a}}& =&-\chi_{bc}\nab_c
F_{\underline{a}}
+\sum_i(\chi_{a_i b}\,\ze_c-\chi_{bc}\,\ze_{a_i}+\in_{a_ic}\,^\star
\b_b)F_{a_1\ldots c\ldots a_k}\nn
\eea
In particular for scalars f,
\be{eq:commscalarf}
\ddd_L\nab_b f-\nab_b\ddd_L f=-\chi_{bc}\nab_c f
\end{equation}
Also, for a one form $F$,
\be{eq:commutationdiv}
L(\div F)-\div(\ddd_LF)=-\chi\c \nab F  +\big(\f12 \trch\c\ze
 +\chih\c\ze -\b\big)\c F\nn
\end{equation}
and again for scalars,
\bea
L(\De f)-\De(Lf)&=&- \trch\lap f-2\chih\c\nab^2 f
+\bigg( \nab\c\trch+2\chih\c\ze+\trch\c\ze\bigg)\c\nab f \nn
\eea
\label{prop:commutationL}
\end{proposition}
6.)\quad {\sl Bochner identity}\quad
Bochner identity holds for scalars and tensors on 
surfaces $S=S_s$, $0\le s\le 1$,
\begin{proposition}\quad 

i)\,\, For a scalar function $f$,
\begin{equation}
\label{eq:scalBoch}
\int_{S} |\nab^{2} f|^{2} = \int_{S} |\lap f|^{2} - 
\int_{S} K |\nab f|^{2}
\end{equation}

\noindent
ii)\,\, For a tensorfield  $F$
\begin{equation}
\label{eq:vectorBoch}
\int_{S} |\nab^{2} F|^{2} = \int_{S} |\lap F|^{2} -
\int_{S} K (2|\nab F|^{2}-|\div F|^{2}) + \int_{S} K^{2} |F|^{2}
\end{equation}
\label{prop:Bochner}
\end{proposition}
\begin{remark} The difference between scalars and tensors is 
substantial as terms quadratic in the curvature are much more difficult to control,
see the properties  {\bf K1}-{\bf K2} we state in the next section.
\end{remark}
\subsection{Main geometric properties of $\HH$}\quad
The proof of the Main Theorem of \cite{KR1} was based on
  the  bootstrap 
assumptions {\bf BA1}--{\bf BA4}  concerning  the geometric quantities $\trch$, $\chih$, $\ze$, $
\trchb$, $\chibh$. In this paper we   make
a consistent but  somewhat different
set of  assumptions, {\bf A1}, {\bf A2}, {\bf WS}, {\bf K1}, {\bf K2} 
 concerning $\trch$, $\chih$,
$\ze$  and the  Gauss curvature $K$.
In the following theorem we stress the fact that these assumptions 
follow from the bootstrap assumptions {\bf BA1}--{\bf BA3},
the small curvature flux and initial conditions
$\RR_0$, $\I_0$ of \cite{KR1}.
\begin{proposition}
 The geometric properties of the geodesic foliation of ${\cal H}$ described
below in {\bf A1}, {\bf A2}, {\bf WS}, {\bf K1}, {\bf K2} follow from 
the bootstrap assumptions {\bf BA1}--{\bf BA3},
the small curvature flux  and the initial conditions 
$\RR_0$, $\I_0$ of \cite{KR1}
\end{proposition}
\begin{proof}:\quad It will become obvious that 
{\bf A1}, {\bf A2} are contained in assumptions {\bf BA1}, {\bf BA2}.  
Properties {\bf WS}, {\bf K1}, {\bf K2}
have been carefully derived in section 4 of \cite{KR1} as  consequences
{\bf BA1}--{\bf BA3},
the small curvature flux  and the initial conditions 
$\RR_0$, $\I_0$. 
\end{proof}
We now describe properties {\bf A1}, {\bf A2}, {\bf WS}, {\bf K1}, {\bf K2} and 
their immediate consequences.

The most primitive
assumption is, as in \cite{KR1}:

{\bf A1.}\qquad \qquad 
$
 \sup_{\HH} r\,|\trchav-\frac{2}{r}|\le \De_0,\qquad   \sup_{\HH} r\,|\trch-\trchav|\le
\De_0.\qquad         
$
where $0<\De_0<\f12$ is a sufficiently small constant.

 Based on the assumption  {\bf A1} we could easily deduced, see section 3.7  in \cite{KR1},
 $$r_0+\f12 s\le r\le r_0+\frac{3}{2}s.$$
Moreover,
\be{eq:controlvol}
1  \le \frac{\sqrt|\ga_s|}{\sqrt|\ga_0|}\le 2\big(\frac{3}{2}\big)^6,\qquad \mbox{
for all}\quad 0\le s\le t\le 1
\end{equation}
  i.e.  the volume elements  of $S_s$ and $S_0$ remain
comparable in the interval  $0\le s\le 1$.
As a consequence of  \eqref{eq:controlvol}   we also  infer  that the $L^2(\HH)$ norm,
 defined in the footnote \ref{footnote1} below,  is equivalent to  the product 
 norm on $[0,t]\times S_0$,
\be{eq:L2normHequiv}
\|F\|_{L^2}=\|F\|_{L_t^2L_x^2}=\big(\int_0^t\int_{S_0}|F|^2 \, ds \, d\mu_0\big)^\f12=
\bigg(\int_0^t ds \int_{S_0}|F(s,\om)|^2 \,  \,\sqrt{|\ga_0|} d\om \bigg)^\f12
\end{equation}
We shall also make use of the following norms,
\be{eq:LxLt2norm}
\|F\|_{L_x^\infty L_t^2}=\sup_{\om\in S_0}\big(\int_0^t ds \,|F(s,\om)|^2\big)^\f12
\end{equation}

\be{eq:Lx2Ltinftynorm}
\|F\|_{L_x^2 L_t^\infty}=\| \,\sup_{0\le s\le 1
} |\,F(s,\om)|\,\,\|_{L^2(S_0)}
\end{equation}
as well as, for $1\le p\le \infty$,
\bea
\|F\|_{L_t^\infty L_x^p }&=& \,\sup_{0\le s\le t}\|\,F(s)\,\|_{L^p(S_0)}
\label{eq:LtinftynormLxp}\\
\|F\|_{L_t^2 L_x^p }&=&\big( \,\int_0^t \|\,F(s)\,\|_{L^p(S_0)}^2 ds\big)^{\f12}
\eea
Observe that $\|F\|_{L_t^\infty L_x^2 }\le \|F\|_{L_x^2 L_t^\infty}$.
We recall the following transport lemma, see section 3.7 in \cite{KR1}.
\begin{lemma}
Consider the equation
$\ddd_L F+k\trch F=G$
for S-tangent tensors  $F, G$ on $\HH$. 
Then, for any $p\ge 1$,
\be{eq:LtinftyLx2estimfortransport}
\|F\|_{L_x^pL_t^\infty}\les \|F(0)\|_{L^p(S_0)}+\|F\|_{L_x^pL_t^1}
\end{equation}
\label{le:transportformulaL2}
\end{lemma}
We shall also make us of the following notations:
\begin{definition} Given an arbitrary $S$-tangent tensor $H$ on $\HH=\HH_1$  we denote 
$$\|\nnab F\|_{L^2}=\|\nab F\|_{L^2}+\|\ddd_L F\|_{L^2}.$$
We also introduce the following  norms ,
\beaa
\NN_1(F)&=&\|F\|_{L^2}+\|\nnab F\|_{L^2}\\
&=&\|F\|_{L^2}+\|\nab F\|_{L^2}+\|\ddd_L F\|_{L^2}\\
\NN_2(F)&=&\|F\|_{L^2}+\|\nnab F\|_{L^2}+\|\nab\nnab F\|_{L^2}\\
&=&\|F\|_{L^2}+\|\nab F\|_{L^2}+\|\nab^2 F\|_{L^2}\\
&+&\|\ddd_L
F\|_{L^2}+\|\nab\c\ddd_L F\|_{L^2}\\
\eeaa
\label{def:NNnorms}
where $L^2$ here stands for $L_t^2L_x^2$.
\end{definition}

The second set
of assumptions we need is:

{\bf A2.}\quad 
\beaa
\|\chih\|_{L_x^\infty L_t^2}\,\,,\,\,\,\, \|\ze\|_{L_x^\infty L_t^2}, &\le& \De_0,\\
\|\nab\trch\|_{ L_x^2L_t^\infty}\,\,,\quad \NN_1(\chih)\,\,,
\quad  \NN_1(\ze)  &\le&\De_0
\eeaa
\begin{remark} The assumptions {\bf A2} are essentially the same as 
{\bf BA2} of \cite{KR1} except for the bound on $\mu$ which is not
needed here.
\end{remark}
As in \cite{KR1} we can simplify our various calculations by 
introducing 
  the following  symbolic notations for 
 connection coefficients.
\begin{definition}  We denote by $\Gd$ the collection formed by
 the connection coefficients:  $\Gd =\trch-\frac{2}{r},\chih,\ze$ 
\label{def:symbolic-notation}
\end{definition}
With these notation the assumptions,
{\bf A1} and {\bf A2} take
the form,
\be{eq:BA1BA2}
\|\trch-\frac{2}{r}\|_{L_t^\infty L_x^\infty},\quad  \|\nab\trch\|_{L_x^2
L_t^\infty},\quad\|\Gd\|_{L_x^\infty L_t^2},\quad 
\NN_1[\Gd]\les
\De_0
\end{equation}
The following  inequalities are straightforward consequences of {\bf A1} and {\bf A2},
 see\cite{KR1}:

\begin{lemma}
\label{lem:Basic-est}
The following estimates hold
for  an arbitrary, smooth, $S$-tangent tensorfield  $F$:
\beaa
\|F\|_{L_t^\infty L_x^2},\,\,\, \|F\|_{L_t^\infty L_x^4},\quad\|F\|_{L_t^6
L_x^6}&\les&\NN_1[F]\\
\|F\|_{L_t^\infty L_x^\infty}&\les&\NN_2[F]
\eeaa
\label{corr:calcineq}
\end{lemma}
\begin{lemma}  \label{lem:Direct}
Let  $w$ be   a solution of the scalar transport equation 
\be{eq:transp-u}
\nab_L w = f,\qquad w|_{S_0}=0,
\end{equation}
For  any $p\ge 1$,
  \be{eq:nabcom1}
  \|\nab  w \, \|_{L^{p}_{x}L^\infty_t }\les
  \|\nab f\|_{L^p_x L^1_t} .
 \end{equation} 
\end{lemma}
\begin{proof}: We only need to differentiate 
according to the commutator formula of proposition \ref{prop:commutationL},
$[\nab_{L},\nab]w=-\frac{1}{2}\trch\c \nab w- \chih\c\nab w $
and  then apply lemma \ref{le:transportformulaL2} to the transport
equation,
\beaa
\ddd_L\nab w+\f12\trch\c\nab w=\nab f-\chih\c\nab w
\eeaa
\end{proof}
The notation introduced  in definition \ref{def:symbolic-notation} 
allows us to express in a compact form
the commutator formulas of proposition \ref{prop:commutationL}
More precisely,
\begin{proposition} In what follows we denote arbitrary $S$ tangent tensorfields by 
capital letters $F$ and scalars by low case letters $f$.
\bea
[\ddd_L,\nab] f&=&-\f12\trch \c\nab f+\Gd \c \nab f\label{eq:comm-with-nab-scalar}
\eea
\bea
[\ddd_L,\lap] f&=&\nab\big( (\frac{1}{r}+\Gd)\c\nab f\big)+\nab\Gd\c\nab f+\Gd\c\Gd\c \nab f
\label{eq:comm-with-lap-scalar}
\eea
\beaa
[\ddd_L,\nab] F&=&-\f12\trch \c\nab F+\Gd \c \nab F+\b\c F+(\Gd+\frac{1}{r} )\c \Gd\c F
\eeaa
\beaa
[\ddd_L,\lap] F&=&-\trch\lap f+\Gd\c\nab^2 F+\nab A\c\nab F+(\Gd+\frac{1}{r})\c\Gd\c\nab F\\
&+&\b \c \nab F +\nab\big(\b\c\nab F+(\Gd+\frac{1}{r} )\c \Gd\c F\big)
\eeaa
\label{prop:comm-nabL-symbolic}
\end{proposition}

\begin{lemma}\label{lem:Reverse}
For a given 1-form $F$ let  $w$ be   a solution of the scalar transport equation 
\be{eq:transp-w}
\nab_L w = \div F,\qquad w|_{S_0}=0,
\end{equation}
 and let 1-form  $W$ be a solution of the equation
\be{eq:transp-W}
\nab_L W - \chi\c W =  F,\qquad W|_{S_0}=0.
\end{equation}
Then for  any $1\le p\le 2$,
  \be{eq:nabcom10}
  \|\div W- w \, \|_{L^{p}_{x}L^\infty_t }\les
  \De_0 \|F\|_{L^{\frac {2p}{2-p}}_x L^1_t} .
 \end{equation} 

\end{lemma}
\begin{proof}:\quad 
We commute the  equation \eqref{eq:transp-W} with $\div$, 
using the commutation formula of
proposition \ref{prop:comm-nabL-symbolic},  and subtract the 
transport equation for $w$ we obtain
$$
\nab_L (\div W-w) = \nab \chi\c W +\b\c W+(\Gd+\frac{1}{r} )\c \Gd\c W
$$
Applying the estimate \eqref{eq:LtinftyLx2estimfortransport} of 
lemma \ref{le:transportformulaL2} we infer that 
\beaa
\|\,\div W-w\|_{L^p_x L^\infty_t} &\les & 
\|\nab\chi \c W \|_{L^p_x L^1_t} +
\|\b\c W\|_{L^p_x L^1_t}+\|(\Gd+\frac{1}{r} )\c \Gd\c W\|_{L^p_x L^1_t}\\
&\les &  
\|\Gd\|^2_{L^\infty_x L^2_t} \c\|W\|_{L^p_x L^\infty_t} +
(\|\b\|_{L^2_t L^2_x} + \|\nab\Gd\|_{L^2_t L^2_x})\c
\|W\|_{L^{\frac {2p}{2-p}}_x L^2_t}
\eeaa
Again applying the estimate \eqref{eq:LtinftyLx2estimfortransport} of 
lemma \ref{le:transportformulaL2}  to the transport equation for
$W$ and using the conditions {\bf A1}, {\bf A2}, we derive
\be{eq:almost}
\|\,\div W-w\|_{L^p_x L^\infty_t} \les  
\De_0 \,\big ( \|F\|_{L^{\frac {2p}{2-p}}_x L^1_t}+ \De_0 \,\|F\|_{L^p_x L^1_t} \big )
\les \De_0 \,\|F\|_{L^{\frac {2p}{2-p}}_x L^1_t}
\end{equation}
as desired.
\end{proof}

In addition to {\bf A1}--{\bf A2} we need two other type of  assumptions.

{\bf WS.}\quad  The initial surface $S_0$ can be covered with
a finite number of coordinate charts $(w^1, w^2)$ such that relative to 
the transported coordinates $(s,\om^a)$ on $\HH$, with $s$ the afine parameter,
 the metric $\ga$ 
 and  its partial derivatives $\pr\ga$, relative to the coordinates $(s,\om)$,
verify the estimates 
\be{eq:weak-regular}
\|\ga\|_{L_t^\infty L_x^\infty}, \qquad\|\ga^{-1}\|_{L_t^\infty L_x^\infty}\les 1,\qquad
\|\pr(\ga-\cga_s)\|_{L_x^2L_t^\infty}\les
\De_0
\end{equation}
where $\cga_s=(1+s)^2\cga$  and
$\cga$ denotes the standard  metric on
$S$, isometric to that of 
$\,\,{\Bbb S}^2$. We shall also make
assumptions on   the Gauss curvature $K$
of the surfaces $S_s$.

{\bf K1.}\quad The Gauss curvature $K$  of the $S_s$ surfaces and the null  curvature component
 $\b$, see \eqref{eq:I2},  verify:
\beaa
 \|K-\frac{1}{r^2}\|_{L_t^2L_x^2}, \quad    \|\b\|_{L_t^2L_x^2} &\les& \De_0\\
\eeaa

{\bf K2.}\quad The Gauss curvature $K$ of the $S_s$ surfaces satisfies 
$$
\|\La^{-\ga} \big(K-r^{-2}\big)\|_{L^2_x L^\infty_t}\les \De_0
$$
with $\La^{-\ga}=(1-\Delta)^{-\ga/2}$, for any $\ga>\f12$.

The properties {\bf WS}, {\bf K1}, {\bf K2} allow us to apply all   the results of  
 \cite{KR2}. In what follows we shall present a summary of
the results proved in \cite{KR2}   which shall be needed in this paper.
\subsection{Calculus inequalities on surfaces}
\begin{proposition} The following calculus inequalities hold true for 
our surfaces $S=S_s$ for any tensorfield $F$.
\be{eq:GNirenberg}
\|F\|_{L^p(S)}\les \|\nab F\|_{L^2(S)}^{1-\frac{2}{p}}\c\|F\|_{L^2(S)}^{\frac 2 p}
+\|F\|_{L^2(S)},\qquad 2\le p<\infty
\end{equation}
Also, for every $2\le p<\infty$,
\be{eq:LinftyL2tensor}
\|F\|_{L^\infty(S)}\les \|\nab^2 F\|^{\frac 1p}_{L^2(S)} 
\|\nab F\|_{L^2(S)}^{\frac {p-2}{p}} \|F\|_{L^2(S)}^{\frac 1{p}} + \|\nab F\|_{L^2(S)}
\end{equation}
As a consequence of the B\"ochner identity for tensors, 
see proposition \ref{prop:Bochner}
\bea
\|\nab^2 F\|_{L^2(S)} &\les &\|\Delta F\|_{L^2(S)} + 
\|K\|_{L^2(S)}^{\frac p{p-1}}
 \|\nab F\|_{L^2(S)}^{\frac {p-2}{p-1}} \|F\|_{L^2(S)}^{\frac 1{p-1}}\nn\\
& + &
 \|K\|_{L^2(S)}\|\nab F\|_{L^2(S)}\label{eq:conseq-Bochner-ineq}
\eea
while for scalars $f$,
\be{eq:conseq-Bochner-ineq-scalar}
\|\nab^2 f\|_{L^2(S)}+\|\nab f\|_{L^2(S)}\les \|\Delta f\|_{L^2(S)}
\end{equation}
\end{proposition}
\begin{proof}:\quad For the first three inequalities see  \cite{KR2}.
The proof of \eqref{eq:conseq-Bochner-ineq-scalar}  can be found in  \cite{KR1} section 4.
\end{proof}
\subsection{Properties of the heat flow} Given a tensor $F$ on $S=S_s$
we define the corresponding heat flow $U(\tau) F$ 
to be   the unique solution
of the equation,
$$\pr_\tau U(\tau)F -\lap U(\tau)F=0, \,\,\quad U(0)F=F.$$
Here $\lap$ denotes the standard Laplace-Beltrami operator on tensors,
$$\lap G=\ga^{ij}\nab_i\nab_j G.$$
\begin{proposition}\label{prop:heat-flow-estimates}
The heat flow $U(\tau)F$ 
verifies the following properties:
\bea
\|U(\tau) F\|_{L^p(S)}&\les &\|F\|_{L^p(S)},\qquad 1\le p\le \infty\label{eq:lpheat1-tensor}\\
\|\nab U(\tau) F\|_{L^2(S)}&\les& \|\nab F\|_{L^2(S)}\label{eq:l2heat2-nab-tensor}
\\ \|\nab U(\tau) F\|_{L^2(S)}&\les& \tau^{-\f12}\|F\|_{L^2(S)}\label{eq:l2heat2-tensor}\\
\| U(\tau)\nab F\|_{L^2(S)}&\les& \tau^{-\f12}\|F\|_{L^2(S)}
\label{le:L2heat-tensor}\\ 
\|\lap U(\tau) F\|_{L^2(S)}&\les& \tau^{-1}\|F\|_{L^2(S)}
\label{eq:l2heat3-tensor}
\eea
Also, for $2\le p<\infty$,
\bea
\|U(\tau) F\|_{L^p(S)}\les \big(1+\tau^{-(1-2/p)}\big)\|F\|_{L^2(S)}\label{eq:heat-GN}
\eea
and the dual estimate, for $1<q\le 2$,
\bea
\|U(\tau) F\|_{L^2(S)}\les \big(1+\tau^{(1-2/q)}\big)\|F\|_{L^q(S)}\label{eq:dual-heat-GN}
\eea
In addition, if $f$ is a scalar function\footnote{We do  not know if such estimate
holds in the tensor case. This failure is also connected with the absence of strong
tensor Bernstein inequality, to be discussed in the next subsection.}
\be{eq:strong-scalar-heat}
\|U(\tau) f\|_{L^\infty(S)}\les \big(1+\tau^{-1}\big)\|f\|_{L^2(S)}
\end{equation}
and its dual 
\be{eq:strong-scalar-heat-dual}
\|U(\tau) f\|_{L^2(S)}\les \big(1+\tau^{-1}\big)\|f\|_{L^1(S)}
\end{equation}
\end{proposition}
\begin{proof}:\quad
See \cite{KR2}.
\end{proof}
\subsection{Geometric LP-projections}
Finally we recall  below the definition  and main properties of the Littlwood-Paley(LP)
 projections introduced in \cite{KR2}. 

\begin{definition}
Consider  the class $\cal M$ of smooth functions $m$ on $[0,\infty)$,
vanishing sufficiently fast at $\infty$,
verifying the  vanishing  moments property:
\be{eq:moments}
\int_0^\infty \tau^{k_1}\pr_\tau^{k_2} m(\tau) d\tau=0, \,\,\,\,
|k_1|+|k_2|\le N 
\end{equation}
 We set,
  $$m_k(\tau)=2^{2k}m(2^{2k}\tau)$$  and 
 define the geometric Littlewood -Paley (LP) 
projections $P_k$, associated to the LP- representative 
function $m\in \MM$, for   arbitrary tensorfields  $F$ on a given surface  $S=S_s$, $0\le s\le 1$,
to be 
\be{eq:LP}P_k F=\int_0^\infty m_k(\tau) U(\tau) F d\tau
\end{equation}
where $U(\tau)F$ is the heat flow on $S$. 

Given an interval $I\subset \Bbb Z$ we define $$P_I=\sum_{k\in I} P_k f.$$
In particular we shall use the notation $P_{<k}, P_{\le k}, P_{>k}, P_{\ge k}$.
\end{definition}
Observe that $P_k$ are selfadjoint. They verify the following properties:
\begin{proposition} The following properties of
the LP 
projections depend only on the conditions {\bf WS}, {\bf K1}, {\bf K2}.

i)\quad {\sl $L^p$-boundedness} \quad For any $1\le
p\le \infty$, and any interval $I\subset \Bbb Z$,
\be{eq:pdf1}
\|P_IF\|_{L^p(S)}\les \|F\|_{L^p(S)}
\end{equation}

ii) \quad {\sl $L^p$- almost orthogonality}\quad  Consider two families
of LP-projections $P_k, \tilde P_k$ associated to $m$ and  respectively 
$\tilde m$, both in ${\cal M}$. For any  
$1\le p\le
\infty$:
\be{eq:pdf2}
\|P_k\tilde P_{k'}F\|_{L^p(S)}\les 2^{-2|k-k'|} \|F\|_{L^p(S)}
\end{equation}
iii) \quad  {\sl Bessel inequality} 
$$\sum_k\|P_k F\|_{L^2(S)}^2\les \|F\|_{L^2(S)}^2$$
iv)\quad {\sl Reproducing property\footnote{see precise statement in\cite{KR1}.}} \quad 
 Given an appropriately defined   $\bar m\in \MM$ there exists $m\in \MM$ such that 
 such that $\bar m=  m\star m$. Thus,
$$^{(\bar m)}P_k =^{(m)}P_k\c ^{(m)}P_k.$$
Whenever there is no danger of confusion we shall simply write $P_k=P_k\c P_k$.

v)\quad {\sl Finite band property}\quad For any $1\le p\le \infty$, $k\ge 0$,
\beaa
\|\lap P_k F\|_{L^p(S)}&\les& 2^{2k}
\|F\|_{L^p(S)}\qquad\qquad\qquad\qquad\qquad[\bf{\lap\text{FB}}]\\
\|P_kF\|_{L^p(S)} &\les& 2^{-2k} \|\lap F 
\|_{L^p}\qquad\qquad\qquad\qquad[\bf{\lap\text{FB}}^{-1}]
\eeaa
In addition, the $L^2$ estimates
\beaa
\|\nab P_k F\|_{L^2(S)}&\les& 2^{k}
\|F\|_{L^2(S)}\qquad\qquad\qquad\qquad\qquad[\bf{\nab\text{FB}}]\\
\|P_kF\|_{L^2(S)} &\les& 2^{-k} \|\nab F 
\|_{L^2}\qquad\qquad\qquad\qquad[\bf{\nab\text{FB}}^{-1}]
\eeaa
hold together with the dual estimate
$$\| P_k \nab F\|_{L^2(S)}\les 2^k \|F(S)\|_{L^2}$$

vi) \quad{\sl Weak Bernstein inequality}. \quad For any $2\le p<\infty$
\begin{align*}
&\|P_k F\|_{L^p(S)}\les (2^{(1-\frac 2p)k}+1) \|F\|_{L^2(S)},\qquad\qquad\qquad\qquad
\, [\text{{\bf wB}}]\\
&\|P_{<0} F\|_{L^p(S)}\les \|F\|_{L^2(S)}
\end{align*}
together with the dual estimates 
\begin{align*}
&\|P_k F\|_{L^2(S)}\les (2^{(1-\frac 2p)k}+1) \|F\|_{L^{p'}(S)},\qquad\qquad\qquad
\qquad [\text{{\bf wB}}^*]\\
&\|P_{<0} F\|_{L^2(S)}\les \|F\|_{L^{p'}(S)}
\end{align*}

vii)\quad {\sl Strong Scalar Bernstein Inequality}\quad For any scalar 
function $f$ and $k\ge 0$
\begin{align*}
&\|P_k f\|_{L^\infty(S)}\les 2^k \|f\|_{L^2(S)}\qquad\qquad\qquad\qquad
\qquad [\text{{\bf ssB}}],\\
&\|P_{<0} f\|_{L^\infty(S)}\les \|f\|_{L^2(S)}
\end{align*}
and the dual estimates,
\beaa
\|P_k f\|_{L^1(S)}&\les& 2^k \|f\|_{L^2(S)}\qquad\qquad\qquad\qquad
\qquad [\text{{\bf ssB}}^*]\\
\|P_{<0}f\|_{L^1(S)}&\les&  \|f\|_{L^2(S)}
\eeaa

In addition we have the following curvature dependent estimates.

viii)\quad{\sl Strong Tensor Bernstein Inequality} \quad 
For any tensor-field $F$, $k\ge 0$
\begin{align*}
&\|P_k F\|_{L^\infty(S)} \les \big (2^k + 2^{k\frac {p-2}{p-1}} 
\|K\|_{L^2(S)}^{\frac
1{p-1}}\big )\c
\|F\|_{L^2(S)},\quad\qquad [\text{{\bf stB}}]\\
& \|P_{<0} F\|_{L^\infty(S)} \les \big (1+
\|K\|_{L^2(S)}^{\frac
1{p-1}}\big )\c\|F\|_{L^2(S)}
\end{align*}

ix)\quad{\sl Dyadic B\"ochner inequality}\quad 
For any tensor-field $F$ and  $2\le p<\infty$, $k\ge 0$,
\begin{align*}
&\|\nab^2 P_k F\|_{L^2(S)} \les \big (2^{2k} + 2^k \|K\|_{L^2(S)} + 2^{k\frac
{p-2}{p-1}}
\|K\|_{L^2(S)}^{\frac p{p-1}}\big )\c \|F\|_{L^2(S)},\\
&\|\nab^2 P_{<0} F\|_{L^2(S)} \les \big (1 +  \|K\|_{L^2(S)} + 
\|K\|_{L^2(S)}^{\frac p{p-1}}\big )\c \|F\|_{L^2(S)}
\end{align*}
x) \quad{\sl Dyadic $L^\infty$ inequality}\quad 
For any tensor-field $F$ and  $2\le p<\infty$, $k\ge 0$,
\begin{align*}
&\|P_k F\|_{L^\infty} \les 
\big (2^k + 2^{k\frac {p-1}{p}} \|K\|^{\frac 1p}_{L^2} + 2^{k\frac {p-2}{p-1}} 
\|K\|_{L^2}^{\frac 1{p-1}}\big ) \c\|F\|_{L^2},\\
&\|P_{<0} F\|_{L^\infty} \les 
\big (1 + \|K\|^{\frac 1p}_{L^2} + 
 \|K\|_{L^2}^{\frac 1{p-1}}\big )\c \|F\|_{L^2}
\end{align*}
\label{thm:geometric-LP}
\end{proposition}
\begin{proof}:\quad See \cite{KR2}.
\end{proof}
In what follows we outline some of the main differences between the properties of  geometric 
LP theory projections recorded in theorem \ref{thm:geometric-LP} and the properties
{\bf LP1}--{\bf LP5}. 
\begin{enumerate}
\item  The simple, pointwise,
orthogonality property  {\bf LP1} does not hold. The  replacement  by the almost 
orthogonality  \eqref{eq:pdf2} is not going to create major difficulties, however 
 the usual trichotomy properties of products are not longer valid. More precisely,
in the classical LP theory {\sl low-low} interactions of
  the type\footnote{In fact $P_k(f_{<k-2}\c
g_{<k-2})$ }
 $P_k(f_{<k}\c g_{<k})$ are forbidden. This is no longer valid 
 for our geometric LP theory.
\item The geometric LP projections commute with the geometric
laplacean $\lap$ but fail to commute with covariant derivatives.
Because of this one has to  be very careful when applying 
the finite band properties recorded in v).
\item In our applications to null hypersurfaces 
we don't have a  bound\footnote{According
to assumption {\bf K1} we only control $ \|K-\frac{1}{r^2}\|_{L^2(\HH)}$. 
The only bound for the Gauss curvature, on a fixed surface $S=S_s$,  we posses is  given by {\bf K2}.} for
 the quantity $\|K-\frac{1}{r^2}\|_{L^2(S)}$. Because of this we have to be very careful
when we  apply the strong Bernstein estimates (  $L^2-L^\infty$)
for tensorfields. However, we do have an unconditional strong Bernstein inequality
for scalars.
\item In flat Minkowski space both the classical
and  geometric   LP projections  commute\footnote{modulo a $\frac 1r$ term generated by the
mean curvature of the sphere foliation of the Minkowski null cone.} with $\ddd_L$
derivatives. This is no longer true for the geometric LP
projections for null hypersurfaces on curved backgrounds.
Moreover, due to our weak regularity assumptions ${\bf BA1}-{\bf BA2}$
as well as {\bf WS}, {\bf K1}, {\bf K2}, the commutators are often not
any better, in terms of their regularity properties, than the principal term.
\end{enumerate}

\subsection{Besov spaces on surfaces $S=S_s$}\quad

The following result was proved in \cite{KR1}
\begin{proposition} $\,$\newline {\bf i}.)\quad Consider the  LP projections
 $P_k$ associated
to an arbitrary $m\in{\cal M}_2$. Then,
\bea
\sum_k \|P_k f\|_{L^2(S)}^2&\les& \|f\|_{L^2(S)}^2\label{onesided1}\\
\sum_k 2^{2k}  \|P_k f\|_{L^2(S)}^2&\les& \|\nab f\|_{L^2(S)}^2\label{onesided2}
\eea
{\bf ii}.)\quad If in addition the LP-projections $P_k$ verify:
\be{eq:sumPk2}
\sum_kP_k^2=I
\end{equation}
Then, 
\bea
\|f\|_{L^2(S)}^2&=&\sum_k\|P_k f\|_{L^2(S)}^2\\
\|\nab f\|_{L^2(S)}^2&\les&\sum_k2^{2k}\|P_k f\|_{L^2(S)}^2
\eea
\label{prop:Sobequivalence}
\end{proposition}
 
Using a family  of LP projections  $P_k$ verifying  $\sum P_k^2=I $
we can now  define our main Besov type norms: 
\begin{definition} Given an arbitrary tensor $F$ on a fixed  $S=S_s$ we  
 define the Besov norm $B^{a}_{2,1}(S)$ for every $0\le \a <\infty$,
\be{eq:Besov}
\|F\|_{B^a_{2,1}(S)}=\big(\sum_{k\ge 0} 2^{a k}\|P_k F\|_{L^2(S)}\big)+
\|\sum_{k< 0} P_k F\|_{L^2(S)}.
\end{equation}
\end{definition}
We recall the following product estimates:
\begin{proposition}
\label{prop:Product1}
Let $a,a', b,b' >0$ such that $a + b =a'+b'\ge 1$.
Then for all tensorfields $F, G$ and 
 any $0\le c <1$,
\begin{equation}
\label{eq:prod'}
\|F\c G\|_{B^c_{2,1}(S)}\les 
\|\La^{a+c} F\|_{L^{2}(S)} \|\La^{b} G\|_{L^{2}(S)}+\|\La^{a'} G\|_{L^{2}(S)} 
\|\La^{b'+c}G\|_{L^{2}(S)}
\end{equation}
\end{proposition}
\begin{proof}:\quad See \cite{KR2}.
\end{proof}
We shall also need the following estimate connecting the norms 
$B^1_{2,1}(S_0)$ and $B^0_{2,1}(S_0)$ for scalars.
\begin{proposition} Given a scalar function $f$ on $S$  we have the inequality:
\be{eq:estimate-B1-B0}
\|f\|_{B^1_{2,1}(S)}\les \|f\|_{B^0_{2,1}(S)}+\|\nab f\|_{B^0_{2,1}(S)}
\end{equation}
\end{proposition}
\begin{proof}:\quad See \cite{KR2}.

\end{proof}
\subsection{Besov spaces on null hypersurfaces $\HH$}\quad

Using the geometric  LP projections 
we are ready to define our main Besov
type norms on $\HH$.
\begin{definition}
For $S$-tangent tensors $F$ on $\HH$  we   introduce the norms, for $0\le a \le 1$:
\bea
\|F\|_{\BB^a}&=&\sum_{k\ge 0} 2^{ak}\|
P_k F\|_{L_t^\infty L_x^2}+\|P_{<0}F\|_{L_t^\infty
L_x^2},\label{eq:Besovnorm}\\
\|F\|_{\PP^a}&=&\sum_{k\ge 0} 2^{ak} \| P_k F\|_{L_t^2 L_x^2}+\|P_{<0}F\|_{L_t^2
L_x^2}\label{eq:Penrosenorm}
\eea
\end{definition}
The following is a crucial result allowing us to pass from tensorial $\BB^0$ estimates 
to their $\BB^0$ scalar counterparts. 
\begin{proposition}
\label{prop:equiv-Besov} There exist a finite number of vectorfields
$X_1,\ldots X_l$ verifying the following properties\footnote{recall
that $\cnab$ represents the covariant derivative with
respect to the background metric $(1+s)^2 \cga$.},
\beaa
\|X,\cnab X\|_{L_t^\infty L_x^\infty}\les 1, \quad \|(\nab- \cnab)X\|_{L_x^2L_t^\infty}\les\De_0,\quad 
\|\nab( \cnab\, X)\|_{ L_x^2L_t^\infty}\les 1
\eeaa

An arbitrary $S$-tangent tensor  $F\in L_t^\infty L_x^2$ is in   $ \BB^0$  if and only if
$F\c X_i\in\BB^0$ for all $1\le i\le l$.
Moreover the  $X_i$'s can be chosen to be coordinate 
vectorfields with $\ddd_L X_i=0$.
\end{proposition} 
\begin{proof}:\quad
To prove  the proposition we 
first  choose the vectorfields $X$
to be the coordinate vectorfields  associated to the transported coordinates
$(t,\om)$. The conditions $\|X\|_{L_t^\infty L_x^\infty }\les 1$
 and $\|(\nab -\cnab) X\|_{L_x^2L_t^\infty}\les \De_0$  are equivalent
to our  {\bf WS} condition. Moreover $\cnab\,^2 X$ are clearly bounded
and the condition on  $\nab\cnab X$ follows from $\|\pr\ga\|_{L_X^2L_t^\infty}\les 1$.

With the help of these vectorfields we  note 
 the following characterization of the $L_t^\infty L_x^2$
norm for tensors\footnote{for simplicity of
notations we restrict to 1-forms.} $F$.
\be{eq:L2-equiv}
\|F\|_{L_t^\infty L_x^2}\approx \max_{1\le i\le l}\|X_i\c F\|_{L_t^\infty L_x^2}
\end{equation}
We now proceed as follows\footnote{We only show the $X_i\c F\in \BB^0$
implies $F\in \BB^0$. This is the implication  which will be needed later. Observe
that we also drop the term $\|P_{<0}F\|_{L_t^\infty L_x^2}\les\|F\|_{L_t^\infty L_x^2} $, which is
trivial, in the expression for $\|F\|_{\BB^0}$}:
\beaa
\|F\|_{\BB^0}&=&\sum_{k\ge 0}\|P_k F\|_{L_t^\infty L_x^2}\approx 
\max_{1\le i\le l} \sum_{k\ge 0}\| X_i\c P_k F\|_{L_t^\infty L_x^2}\\
&=&\max_{1\le i\le l} \sum_{k\ge 0}\|  P_k(X_i\c F)\|_{L_t^\infty L_x^2}
+O\big(\sum_{k\ge 0}\|  [P_k\,,\,X]\c F\|_{L_t^\infty L_x^2}\big)
\eeaa
It suffices to prove that,
\bea
\sum_{k\ge 0}\|  [P_k\,,\,X]\c F\|_{L_t^\infty L_x^2}&\les & \|(\nab-\cnab)
X\|_{L_x^2L_t^\infty}\|F\|_{\BB^0} +\| F\|_{L_t^\infty L_x^2} \label{eq:PkX-comm}\\
&\les&
\De_0 \|F\|_{\BB^0}  +\| F\|_{L_t^\infty L_x^2}.
\eea
To prove \eqref{eq:PkX-comm} we have to recall the formula
\beaa
 [P_k\,,\,X]\c F&=&\int_0^\infty m_k(\tau)\Phi(\tau)\label{eq:PkF-comm}\\
\Phi(\tau;k)&=& \int_0^\tau U(\tau-\tau')[\lap\,,\, X]U(\tau')F d\tau'=\Phi_1(\tau)+\Phi_2(\tau)
+\Phi_3(\tau)+\Phi_4(\tau)\nn\\
\Phi_1(\tau;k)&=&\int_0^\tau U(\tau-\tau')  (\nab-\cnab) X \c  \nab U(\tau')F d\tau'\nn\\
\Phi_2(\tau;k)&=&\int_0^\tau U(\tau-\tau')  \div \big( \nab-\cnab) X \c   U(\tau')F\big) d\tau'\nn\\
\Phi_3(\tau;k)&=&\int_0^\tau U(\tau-\tau') \cnab X\c  \nab U(\tau')F d\tau'\nn\\
\Phi_4(\tau;k)&=&\int_0^\tau U(\tau-\tau')  \div \big( \cnab X \c   U(\tau')F\big) d\tau'
\eeaa
since 
\beaa
[\lap, X]F&=&\nab^i(\nab_i X\c F)+\nab_i X\c\nab^i F
=\nab^i\big((\nab_i-\cnab_i) X\c F\big)+(\nab_i -\cnab_i) X\c\nab^i F\\
&+&\nab^i\big(\cnab_i X\c F\big)+\cnab_i X\c\nab^i F
\eeaa
To estimate $\Phi_1$ we observe that $U(\tau-\tau')$ corresponds to a scalar heat flow
and therefore we make use of the following  scalar heat flow estimate
\eqref{eq:strong-scalar-heat-dual}:
 $$\|U(\tau-\tau') f\|_{L_x^2}\les (1+|\tau-\tau'|^{-\f12})\| f\|_{L_x^1}. $$
\begin{remark}
In what follows, we shall systematically replace the above estimate, and all other
heat flow estimates like it,  with their  
slightly incorrect versions where we ignore the non-singular term $1$.
We thus write
$$\|U(\tau-\tau') f\|_{L_x^2}\les |\tau-\tau'|^{-\f12}\| f\|_{L_x^1}. $$
\end{remark}
Therefore,
\beaa
\|\Phi_1(\tau;k)\|_{L_t^\infty L_x^2} &\les&\int_0^\tau 
|\tau-\tau'|^{-\f12}
\| (\nab-\cnab) X\c \nab U(\tau') F\|_{L_t^\infty L_x^1}\\
&\les&  \int_0^\tau  |\tau-\tau'|^{-\f12} 
 \|( \nab-\cnab) X\|_{ L_x^2L_t^\infty}\c\| \nab U(\tau') F\|_{L_t^\infty L_x^2}\\
&\les&\De_0\c \int_0^\tau |\tau-\tau'|^{-\f12}
\| \nab U(\tau') F\|_{L_t^\infty L_x^2}
\eeaa
To estimate the integral
$$J(\tau)=\int_0^\tau |\tau-\tau'|^{-\f12}  
\| \nab U(\tau') F\|_{L_t^\infty L_x^2}$$
we decompose it as follows:
\beaa
J(\tau)\les \sum_m\int_0^\tau  |\tau-\tau'|^{-\f12}
\| \nab P_m U(\tau') F\|_{L_t^\infty L_x^2}=\sum_mJ_m(\tau)
\eeaa
Next we make use of
\begin{lemma}\label{le:Jmtau}
\be{eq:Jm-tau}
J_m(\tau)\les \min \big(2^m\tau^\f12, 2^{-m}\tau^{-\f12}\big)\|P_m F\|_{L_t^\infty L_x^2}
\end{equation}
\end{lemma}
and proceed  as follows:
$$
\int_0^\infty m_k(\tau)\|\Phi_1(\tau;k)\|_{L_t^\infty L_x^2}\les
\sum_m \|P_m F\|_{L_t^\infty L_x^2}\int_0^\infty m_k(\tau)
\min \big(2^m\tau^\f12, 2^{-m}\tau^{-\f12}\big)
$$
Now,
\beaa
\int_0^\infty m_k(\tau)
\min \big(2^m\tau^\f12, 2^{-m}\tau^{-\f12}\big)&= & 2^m\int_0^{2^{-2m}} m_k(\tau)
\tau^\f12+ 2^{-m}\int_{2^{-2m}}^\infty m_k(\tau)\tau^{-\f12}\\ &=&
2^{m-k}\int_0^{2^{2(k-m)}} \tilde m(\tau)
+ 2^{k-m}\int_{2^{2(k-m)}}^\infty \hat m(\tau)\\ &\les & 2^{-|k-m|}
\eeaa
Here $\tilde m(\tau) = \tau^\f12 m(\tau)$ and $\hat m(\tau) =\tau^{-\f12} m(\tau)$.
Moreover, to arrive at the inequality above we used the following bounds:
\beaa
\int_0^\infty \tilde m(\tau),\quad \int_0^\infty \hat m(\tau) \les 1,\quad
\int_0^a \tilde m(\tau) \les a, \quad\int_A^\infty \hat m(\tau) \les A^{-1}
\eeaa 
which hold for all sufficiently small $a$ and all sufficiently large $A$.
Thus,
\beaa
\sum_{k\ge 0} \int_0^\infty m_k(\tau) \|\Phi_1(\tau;k)\|_{L^2_t L^2_x} \les 
\De_0\sum_{k\ge 0} \sum_m 2^{-|m-k|} \|P_m F\|_{L^2_t L^2_x} \les \De_0\|F\|_{\PP^0}
\eeaa
To estimate $\Phi_2$ we first observe that the following estimate
for the scalar heat flow holds for any $1<p\le 2$:
\be{eq:heat-deriv-p}
\|U(\tau-\tau') \div g\|_{L^2_x} \les (\tau-\tau')^{-\frac 1p} \|g\|_{L^p_x}
\end{equation}
Using this we obtain for some $p<2$ sufficiently close to $p=2$,
\beaa
\|\Phi_2(\tau;k)\|_{L^\infty_t L^2_x} &\les &\int_0^\tau (\tau-\tau')^{-\frac 1p} 
\|(\nab-\cnab) X \c U(\tau') F\|_{L^\infty_t L^p_x} \\
&\les & \int_0^\tau (\tau-\tau')^{-\frac 1p} 
\|(\nab-\cnab) X\|_{L^2_x L^\infty_t} \| U(\tau') F\|_{L^\infty_t L^{\frac {2p}{2-p}}_x} \\
&\les & \De_0 \sum_m \int_0^\tau (\tau-\tau')^{-\frac 1p} 
 \| P_m U(\tau') F\|_{L^\infty_t L^{\frac {2p}{2-p}}_x} \\
&\les & \De_0 \sum_m \int_0^\tau (\tau-\tau')^{-\frac 1p} 
 \| \nab P_m U(\tau') F\|_{L^\infty_t L^2_x}^{2-\frac 2p} \|P_m U(\tau') F\|^{\frac 2p-1}_{L^\infty_t
L^2_x} 
\eeaa
The remaining argument now is a straightforward modification of the proof for
$\Phi_1$. We infer that 
$$
\sum_{k\ge 0} \int_0^\infty m_k(\tau) \|\Phi_2(\tau;k)\|_{L^\infty_t L^2_x} \les \De_0\sum_{k\ge 0} \sum_m
2^{-(2-\frac 2p)|k-m|} \|P_m F\|_{L^\infty_t L^2_x} \les \De_0 \|F\|_{\BB^0}
$$
It only remains to estimate the easier terms $\Phi_3, \Phi_4$.
\beaa
\Phi_3(\tau;k)&=&\int_0^\tau U(\tau-\tau') \cnab X\c  \nab U(\tau')F d\tau' 
=\Phi_{31}(\tau;k)+\Phi_{32}(\tau;k)\\
&=&\int_0^\tau U(\tau-\tau') (\nab\cnab X)\c   U(\tau')F d\tau'+
\int_0^\tau  U(\tau-\tau')\nab \big( \cnab X\c   U(\tau')F\big) d\tau'
\eeaa
Now, starting  as for $\Phi_{1}$, $\Phi_2$,
\beaa
\|\Phi_{31}(\tau;k)\|_{L_t^\infty L_x^2} &\les&\int_0^\tau (\tau-\tau')^{-\f12}
\|\nab \cnab X\c  U(\tau') F\|_{L_t^\infty L_x^1}\\
&\les&  \int_0^\tau (\tau-\tau')^{-\f12}  
 \| \nab\cnab X\|_{ L_x^2L_t^\infty}\c\|  U(\tau') F\|_{L_t^\infty L_x^2}
\les\tau^{\f12}  
\| F\|_{L_t^\infty L_x^2}
\eeaa
\beaa
\|\Phi_{32}(\tau;k)\|_{L_t^\infty L_x^2} &\les&\int_0^\tau (\tau-\tau')^{-\f12}
\|\cnab X\c  U(\tau') F\|_{L_t^\infty L_x^2}\les\tau^\f12\|F\|_{L_t^\infty L_x^2}
\eeaa
Hence, 
\beaa
\sum_{k\ge 0} \int_0^\infty m_k(\tau) \|\Phi_3(\tau;k)\|_{L^\infty_t L^2_x}\les 
\sum_{k\ge 0} \int_0^\infty \tau^\f12 m_k(\tau) d\tau\les \sum_k 2^{-k}
\eeaa
and similarly for $\Phi_4$.
Thus, going back to \eqref{eq:PkF-comm} we have 
$$
\sum_k \|[P_k\,,\,X] F\|_{L^\infty_t L^2_x} \les \De_0 \|F\|_{\BB^0}+\|F\|_{ L_t^\infty  L_x^2}
$$
\end{proof}
\begin{proof}{\bf of lemma \ref{le:Jmtau}.}\quad 
Recall that our goal is to prove the estimate 
$$
\int_0^\tau (\tau-\tau')^{-\f12}  
\| \nab P_m U(\tau') F\|_{L_t^\infty L_x^2}\les \min \big(2^m\tau^\f12, 2^{-m}\tau^{-\f12}\big)\|P_m
F\|_{L_t^\infty L_x^2}.
$$
We first observe that,
\be{eq:2mmin}
\| \nab P_m U(\tau') F\|_{L_t^\infty L_x^2}\les 2^m \min \big(1, 2^{-4m}(\tau')^{-2}\big)\|P_m
F\|_{L_t^\infty L_x^2}
\end{equation}
Indeed we have both
$$\| \nab P_m U(\tau') F\|_{L_t^\infty L_x^2}\les 2^m  \|P_m
F\|_{L_t^\infty L_x^2}$$
and 
\beaa
\| \nab P_m U(\tau') F\|_{L_t^\infty L_x^2}&\les& 2^m  \|P_m
 U(\tau') F\|_{L_t^\infty L_x^2}\\
&\les&  2^{m} 2^{-4m} \| P_m\lap^2
 U(\tau') F\|_{L_t^\infty L_x^2}\\
&\les&  2^{m} 2^{-4m}(\tau')^{-2} \| P_m F\|_{L_t^\infty L_x^2}
\eeaa
To show that, 
\beaa
I_m(\tau)&=&\|P_m F\|_{L_t^\infty L_x^2}^{-1}\c\int_0^\tau(\tau-\tau')^{-\f12}\| \nab P_m U(\tau')
F\|_{L_t^\infty L_x^2}\\
&\les& \min\big( 2^m\tau^\f12\,,\,2^{-m}\tau^{-\f12}\big)
\eeaa
it suffices to prove 
\beaa
I_m(\tau)&\les& 2^m\tau^\f12\qquad\mbox{if}\quad\, \tau\le 2^{-2m}\\
I_m(\tau)&\les& 2^{-m}\tau^{-\f12}\quad\mbox{if}\quad \tau\ge 2^{-2m}
\eeaa
Using 
 \eqref{eq:2mmin} we infer that,
\beaa
I_m(\tau)&\les& 2^m
\int_0^\tau(\tau-\tau')^{-\f12}\min \big(1, 2^{-4m}(\tau')^{-2}\big)d\tau'
\eeaa
For $\tau\le 2^{-2m}$ we  have,
\beaa 
I_m(\tau)&\les& 2^m
\int_0^\tau(\tau-\tau')^{-\f12}d\tau'\les 2^m \tau^\f12
\eeaa
For $\tau\ge 2^{-2m}$,
\beaa
I_m(\tau)&\les& 2^{-3m}
\int_{2^{-2m}}^\tau(\tau-\tau')^{-\f12}\c(\tau')^{-2}d\tau'
+     
2^m \int_0^{2^{-2m}}(\tau-\tau')^{-\f12}d\tau'\\
&\les& 2^{-m}\tau^{-\f12}
\eeaa
as desired.
\end{proof}

\section{Main results}
We are now ready to state our main results. They can be viewed as 
extensions of  results mentioned in the introduction, see
proposition \ref{prop:intr-bilintrace-cone}, to null hypersurfaces
$\HH$ verifying the assumptions {\bf A1}, {\bf A2}, {\bf WS}
and {\bf K1}, {\bf K2}.
 The first
is a generalized sharp bilinear trace theorem.
\begin{theorem}[Bilinear-Trace-Transport]
\label{thm:maintrace}
 Consider the transport equation along $\HH$:
$$
\ddd_L W= \nab_LF\c G
$$
 for $S$-tangent tensors  $W$, $F$, $G$.  Then,
\be{eq:mainlemma1}
\|W\|_{\BB^0}\les \|W|_{S_0}\|_{B^0_{2,1}(S_0)}+\NN_1[F]\c\big (\NN_1[G]+
\|G\|_{L^\infty_x L^2_t}\big )
\end{equation}
\end{theorem}
\begin{remark}
We have a stronger estimate in the case of a transport 
equation for a {\it scalar} function $w$.
\be{eq:mainlemma10}
\|w\|_{\BB^0}\les \|w|_{S_0}\|_{B^0_{2,1}(S_0)}+\NN_1[F]\c\NN_1[G].
\end{equation}
\end{remark}
The next two  theorems are the noncommutative versions of the 
 sharp product estimates  of proposition  \ref{prop:intr-producttrace}.

\begin{theorem}[Product-Transport I]
\label{thm:mainBesov}
 Consider the transport equation along $\HH$:
$$
\ddd_L W= F\c G
$$
 for $S$-tangent tensors  $W$, $F$, $G$.
We have the estimate,
\be{eq:mainlemma2}
\|W\|_{\BB^0}\les  \|W|_{S_0}\|_{B^0_{2,1}(S_0)}+ \|F\|_{\PP^0}\c\big(\NN_1[G]+\|
G\|_{L_x^\infty L_t^2}\big)
\end{equation}
\end{theorem}

\begin{theorem}[Product-Transport II]\label{thm:mainBesov2}
 Given any    pair of $S$-tangent tensors $G, W$,
of same type such that $W$ satisfies a transport equation,
of the form,
$$\ddd_L W=  F.$$
Then,
\be{eq:bilin2}
\| G\c W\|_{\PP^0} \les \big(\|F\|_{\PP^0}+\|W|_{S_0}\|_{B^0_{2,1}(S_0)}\big)\c
\big(\N_1[G]+\|G\|_{L_x^\infty L_t^2} \big)
\end{equation}
\end{theorem}
As a consequence of theorems \ref{thm:maintrace}, \ref{thm:mainBesov} we derive the following.
\begin{theorem}[Sharp-Trace]
For any $S$-tangent tensor $F$, which allows a
 decomposition of the form  $\nab F=\nab_L \check F + G$,
we have
\be{eq:maintheorem3}
\|F\|_{L^\infty_x L^2_t} \les \N_1[F]+ \N_1[\check F] + \|G\|_{\PP^0}
\end{equation}
\end{theorem}
\begin{proof}:
\quad The proof can be found in section 5 of \cite{KR1}.
The idea is to introduce the scalar function  $f(t)=\int_0^t|F|^2$
and observe that it  verifies the transport
equation,
$$\ddd_L f=|F|^2,\qquad U(0)=0.$$
Differentiating it one 
derives,
\beaa
\ddd_L (\nab f)+\f12 \trch (\nab f)&=&2 F \c\nab F-\chih \c(\nab f)\\
&=&2 F\c\ddd_L\check F+ 2F\c G-\chih \c(\nab f)
\eeaa
and  apply\footnote{One has to take some care to eliminate the term 
$\trch \c \nab f$ first, as it is done in \cite{KR1}.} propositions  \ref{thm:maintrace},
\ref{thm:mainBesov}. 
\end{proof}
In particular, we have the following noncommutative  version
of the classical sharp trace theorem.
\begin{corollary}
For any $S$-tangent tensor $F$
\be{eq:maintheorem4}
\int_0^t |\nab_L F|^2\les \int_{\cal H} \Big (|\nab^2 F|^2 + |\nab^2_L F|^2 + |F|^2\Big )
\end{equation}
\end{corollary}
\subsection{Reduction to scalar estimates}\quad

The first two transport theorems can be reduced to the case of a scalar transport equation.
More precisely, we have the following 
\begin{proposition}\label{thm:reduction}
Assume that the conclusions of theorems \ref{thm:maintrace} and \ref{thm:mainBesov}
have been verified for scalar transport equations. Then they also hold true in the tensor
case. Moreover , in the particular case of theorems \ref{thm:mainBesov},
 and \ref{thm:mainBesov2} we can reduce the corresponding estimates 
to a fully scalar situation, i.e both $F$ and $G$  are scalar functions. 
\end{proposition}
\begin{proof}:\quad In view of the scalar characterization of
the space $\BB^0$ stated in proposition \ref{prop:equiv-Besov} it suffices to 
do the following. We multiply the transport equation for a tensor\footnote{Assume for simplicity that $W$
is a 1-form}
$W$, in either of the theorems,  by the vectorfields $X$ to derive a scalar equation
$$
\nab_L (X\c W) = \nab_L F\c (G\widehat{\otimes} X )\qquad \bigg(\mbox {or}\,\,\,F\c (G\widehat{\otimes}
X)\bigg)
$$
where $\widehat{\otimes}$ denotes either a tensor product
or a contraction. It only remains 
to  observe that the norm   $\N_1[G\widehat{\otimes} X]+ \|G\widehat{\otimes} X\|_{L^\infty_x L^2_t}$  
is  invariant with respect to a tensor multiplication by a vectorfield $X$ with the properties guaranteed
by proposition
\ref{prop:equiv-Besov}.
To prove the second part of the proposition  consider the case of theorem 
\ref{thm:mainBesov} where we have already reduced to the case  $\ddd_L w=F\c G$  
with $w$ scalar and $F\c G$ denotes the scalar product between two tensorfields.
Clearly $F\c G$ can be expressed, at every point, as a product between various
scalar components of $F$ and $G$ inner product with the vectorfields $X_i$
 and the components of the metric $\ga$. Therefore we can apply the proposition 
\ref{prop:equiv-Besov} to  each of the scalar 
components\footnote{The components of the metric $\ga$
can be 
combined with $G$. }  and derive our result in view of the invariance of the norms involved.
\end{proof}
The reduction to scalar equations is a very 
important simplification in so far as it allows
us  to work with integral
estimates. We state below a result which, in view of the reduction made above,
implies  
 theorems \ref{thm:maintrace}, \ref{thm:mainBesov}.
\begin{theorem} The following statements hold true
for arbitrary $S$ tangent tensorfields $F, G$ of same order\footnote{such that $F\c G$
denotes a scalar.}:
\be{eq:maintrace-scalar}
\|\int_0^t \ddd_L F\c G\|_{\BB^0}\les \NN_1[F]\c\NN_1[G]
\end{equation}
\be{eq:mainBesov-scalar}
\|\int_0^t  F\c G\|_{\BB^0}\les\|F\|_{\PP^0}\c\big(\NN_1[G]+\|
G\|_{L_x^\infty L_t^2}\big)
\end{equation}
Moreover for solutions of the homogeneous scalar transport equation
 $\ddd_Lw=0$  we have,
\be{eq:homog-transport}
\|w\|_{\BB^0}\les \|w|_{S_0}\|_{\BB^0}
\end{equation}
\label{thm:main-scalar}
\end{theorem}
Indeed, once we are in the scalar case, for example
in the case of  theorem \ref{thm:maintrace},
$\ddd_L w=F\c G$, we can integrate and therefore reduce 
the statement of  the theorem  to,
$$\sum_k\|P_k\int_0^t F\c G\|_{L_t^\infty L_x^2}\les \NN_1[F]\c\big (\NN_1[G]+
\|G\|_{L^\infty_x L^2_t}\big )$$

\begin{remark}\label{rem:mainBesov-dyadic}
We can also prove a more precise dyadic
version of the estimate \eqref{eq:mainBesov-scalar}: for any $k\ge 0$ and  some
$\si >0$,
\bea
\|P_k \int_0^t F\c G\|_{L^\infty_t L^2_x}&\les & 
\Big (\sum_{k'} 2^{-\si{|k-k'|}}
\|P_{k'}F\|_{L^2_t L^2_x} + 2^{-\si k} \|F\|_{L^2_t L^2_x}\Big )\c\nn\\ &\c& 
\Big (\N_1[G] + \|G\|\Big )\label{eq:mainlemma2-dyadic}
\eea
\end{remark}
\section{ Some dyadic estimates} 

In the proof of the theorems we shall need the notion of  an $\NN_1$  envelope
of a tensorfield. It plays the role of an LP -localized  version of the $\NN_1$ norm
of definition \ref{def:NNnorms}.

\begin{definition}\label{def:Envelope}
For a given smooth $S$-tangent tensor-field $F$ and a sufficiently
small $\epsilon >0$  we define its
$\N_1$ -envelope (of order\footnote{Unless otherwise specified 
we shall assume that $\epsilon$  is a fixed, sufficiently small constant. } $\epsilon$) to be
{\it any} sequence 
 of positive real numbers $\N_1[F_k]$ satisfying 
the following properties:  
 \,\,\,
\beaa
\N_1[F_k]&\les& 2^{\epsilon |k-k'|} \N_1[F_{k'}],\qquad\mbox{for any $k, k'$,}
\\
\sum_k \N_1[F_k]^2 &\approx& \N_1[F]^2,
\eeaa
\end{definition}
The existence of an envelope follows from the following elementary
construction. Let $\bar{\N_1}[F_k]$,  be defined as follows:  
\begin{align}
&\bar{\N_1}[F_k] = \|F_k\|_{L^2_t L^2_x}+ \|\nab F_k\|_{L^2_t L^2_x}
+ \|(\nab_L F)_k\|_{L^2_t L^2_x},\label{eq:N1-dyadic}
\end{align}
Note that, in view of proposition \ref{prop:Sobequivalence}, 
$$
\sum_k\bar{ N_1}[F_k]^2 \approx \N_1[F]^2.\qquad 
$$
We now easily check that 
the sequences 
\beaa
\N_1[F_{k}]&=&\sum_{k'} 2^{-\epsilon |k-k'|} \bar{\N_1}[F_{k'}]
\eeaa
are  desired  envelopes.
The following simple result provides us with  a useful tool in handling various 
error terms.
\begin{lemma}
Let $\{\N_1[F_k]\}$ be an envelope (of order $\ep$) for a tensor-field $F$. Then 
for any $0<\a\le \epsilon$ the sequence $\{\N_1[F_k] + 2^{-\a k}\N_1[F]\}$
is also an envelope.   Moreover, for any $\a >\epsilon$ the sequence 
$\{\N_1[F_k] + 2^{-\a k}\N_1[F]\}$ is dominated by an $\NN_1$-envelope for
$F$.
\end{lemma}
\begin{proof}:\quad The first part of the lemma is obvious from the definition. 
From a purely technical point of view $\{\N_1[F_k] + 2^{-\a k}\N_1[F]\}$
is not an an envelope of order $\ep$ if $\a\ge \ep$ yet it can be clearly dominated by one.
\end{proof}
We now formulate and prove a number of results which will be routinely
used in the proof of our main theorem. The next result allows us to treat the
multitude of  commutator terms with $\ddd_L$  which will   appear  throughout
the next sections.
\begin{lemma}\label{lem:Comm-nabL}
For any smooth $S$-tangent tensor field $F$ and all $q<2$ sufficiently close
 to $q=2$\footnote{By    $2^{-\frac k2+}$ we mean $2^{-ak}$ with $a<1/2$ arbitrarily
close to $1/2$.  }, 
\be{eq:corq-L}
\|[P_k,\nab_L] F\|_{L^q_t L^2_x} + 2^{-k} \|\nab [P_k,\nab_L] F\|_{L^q_t L^2_x} 
\les  2^{-\frac k2+} \N_1[F] 
\end{equation}
while for $q=1$,
\be{eq:corq-L1}
\|[P_k,\nab_L] F\|_{L^1_t L^2_x} + 2^{-k} \|\nab [P_k,\nab_L] F\|_{L^1_t L^2_x} 
\les  2^{-k+} \N_1[F] 
\end{equation}
\end{lemma}
\begin{proof}:\quad See  section \ref{sec:Commutator}.
\end{proof}
As a corollary of this we see
 how to control the $\nab_L$-derivative of the LP pieces
of a tensor-field $f$ in terms of its envelope.
\begin{lemma}\label{lem:Equiv-L}
For all $q, \,1\le q <2$ sufficiently close to $q=2$ 
and any smooth $S$-tangent tensor-field $F$,
\be{eq:equiv-L}
\|\nab_L F_k\|_{L^q_t L^2_x} \les \N_1[F_k]
\end{equation}
for some $\N_1$-envelope for $f$.
\end{lemma}
\begin{proof}:\quad
We write
$$
\nab_L F_k = P_k \nab_L F + [P_k,\nab_L]  F
$$
The commutator estimate of Lemma \ref{lem:Comm-nabL} implies that 
$$
\| [P_k,\nab_L] F\|_{L^q_t L^2_x} \les 2^{-\frac k2} \N_1[F]
$$
Thus,
\beaa
\|\nab_L F_k\|_{L^q_t L^2_x}\les  
\|P_k\nab_L F\|_{L^q_t L^2_x} +  2^{-\frac k2}
\N_1[F] \les\N_1[F_k]
\eeaa
\end{proof}
The above result can be complemented by a 
dyadic Gagliardo-Nirenberg estimate with respect to the time
variable $t$.
\begin{lemma}\label{lem:Gagl-Nir-L}
For any smooth $S$-tangent tensor field $f$, any $2\le q\le \infty$, 
we have the following dyadic  Gagliardo-Nirenberg inequality
$$
\qquad\qquad\qquad\|F_k\|_{L^q_t L^2_x} \les
  2^{-k(\frac 12+\frac 1q)}\N_1[F_k]\qquad\qquad\qquad 
\qquad\qquad\qquad{\bf
[GN_k]}
$$

\end{lemma}
\begin{proof}:\quad
First, we trivially estimate 
$$
\|F_k\|^q_{L^q_t L^2_x}\les \|F_k\|_{L_t^2L_x^2}^2 
\|F_k\|_{ L^2_xL^\infty_t}^{q-2}  
$$
On the other hand,
\beaa
\|F_k\|^2_{L^\infty_t} &\les &\|\nab_L F_k\c f_k\|_{L^1_t} + \|F_k\|^2_{L^2_t}\\
&\les &  \|(\nab_L F)_k\c F_k\|_{L^1_t} + \|[P_k,\nab_L]F\c F_k\|_{L^1_t} + \|F_k\|^2_{L^2_t}
\eeaa
Thus,
$$
\|F_k\|_{ L^2_xL^\infty_t}^2 \les \|(\nab_L F)_k\|_{L_t^2L_x^2} 
\|F_k\|_{L_t^2L_x^2} + \|[P_k,\nab_L]F\c F_k\|_{L^1_t L^1_x}
+ \|F_k\|_{L_t^2L_x^2}^2
$$
We choose an exponent $r<2$ sufficiently close to $r=2$ and estimate the
commutator term as follows.
\beaa
\|[P_k,\nab_L]F\c F_k\|_{L^1_t L^1_x}&\les & 
\|[P_k,\nab_L]F\|_{L^r_t L^2_x} \|F_k\|_{L^{r'}_t L^2_x}\\
&\les & 2^{-\frac k2} \N_1[F] \|F_k\|_{L^{r'}_t L^2_x}
\eeaa
Combining all the estimates we obtain that 
\beaa
\|F_k\|^q_{L^q_t L^2_x}\les \|F_k\|_{L_t^2L_x^2}^2 
\Big (\|(\nab_L F)_k\|_{L_t^2L_x^2}
\|F_k\|_{L_t^2L_x^2} +
2^{-\frac k2} \N_1[F] \|F_k\|_{L^{r'}_t L^2_x} + 
\|F_k\|_{L_t^2L_x^2}\Big )^{\frac q2-1}
\eeaa
It follows that for all $r'>2$ sufficiently close to $r'=2$,
$$
\|F_k\|_{L^q_t L^2_x}\les 
\|(\nab_L F)_k\|_{L_t^2L_x^2}^{\frac 12-\frac 1q}
\|F_k\|_{L_t^2L_x^2} ^{\frac 12+\frac 1q} + \|F_k\|_{L_t^2L_x^2}
+ 2^{-\frac k2(\frac 12-\frac 1q)} (\N_1[F] )^{\frac 12-\frac 1q}
\|F_k\|_{L^{r'}_t L^2_x} ^{\frac 12-\frac 1q}
\|F_k\|_{L_t^2L_x^2}^{\frac 2q}
$$
Using  the above estimate with $q=r'$ and
then plugging the result into the above estimate for a given $q$, 
we obtain that for any $\a<\frac 14 -\frac 1{2q}$,
$$
\|F_k\|_{L^q_t L^2_x} \les 2^{-k(\frac 12+\frac 1q)} 
\Big (\N_1[F_k] + 2^{-\a k} \N_1[F]\Big ).
$$
The desired result now easily follows.
\end{proof}

\begin{lemma}\label{lem:Gagl-Nir-nab4}
For any smooth $S$-tangent tensor-field $F$ and any Lebesque
exponent $2\le q <4$, 
\be{eq:trickyBernstein}
\|\nab  \,F_k\|_{L^q_t L^4_x}\les 2^{k(1-\frac 1q)} \N_1[F_k]
\end{equation}
\end{lemma}
\begin{proof}
By Gagliardo-Nirenberg \eqref{eq:GNirenberg},
\beaa
\|\nab F_k\|_{L^q_t L^4_x}\les \big {\|} \|\nab^2 F_k\|_{L^2_x}^\f12 
\|\nab F_k\|_{L^2_x}^\f12 \big {\|} _{L^q_t}
\eeaa
On the other hand, according to the Bochner inequality\footnote{Simplified a bit and
ignoring lower order terms.}
\eqref{eq:conseq-Bochner-ineq},
$$
\|\nab^2 F_k\|_{L^2_x} \les \|\Delta F_k\|_{L^2_x}  + 
\|K\|_{L^2_x}^{\frac p{p-1}} \|\nab F_k\|_{L^2_x}
$$
for any $2\le p<\infty$. Choose the exponent $p$ such that 
$$
2=\frac {qp}{2(p-1)}
$$
The existence of such $p$ is guaranteed by the condition that $q<4$.
Thus, using H\"older inequality and condition {\bf K1},
\beaa
\|\nab F_k\|_{L^q_t L^4_x}&\les &\|\Delta F_k\|_{L^q_t L^2_x}^\f12 
\|\nab F_k\|_{L^q_t L^2_x}^\f12  +  \|\nab F_k\|_{L^\infty_t L^2_x} 
\\ &\les & 2^{\frac {3k}2} \|F_k\|_{L^q_t L^2_x } + 2^k 
\|F_k\|_{L^\infty_t L^2_x} 
\eeaa
 It remains to apply the Gagliardo-Nirenberg inequality $[{\bf GN_k}]$ of Lemma \ref{lem:Gagl-Nir-L}.
\end{proof}
Our next result is the integrated version of the strong Bernstein inequality.
\begin{lemma}\label{lem:Int-SB-q}
For any  
$S$-tangent tensor-field $F$ and exponent $2\le q <\infty$,
\be{eq:int-strB-q}
\|F_k\|_{L^q_t L^\infty_x} \les 2^{k(\frac 12-\frac 1q)}\,\N_1[F_k]
\end{equation}
\end{lemma}
\begin{proof}:\quad 
Observe that the dyadic $L^\infty$ inequality of x) of proposition \ref{thm:geometric-LP}
implies that for all sufficiently large $p$
$$
\|F_k\|_{L^\infty_x} \les 2^k(1+ 2^{-\frac 1{p}k} \|K\|_{L^2_x}^{\frac 1{p}}+ 
2^{-\frac 1{p-1}k} \|K\|_{L^2_x}^{\frac 1{p-1}})
\|F_k\|_{L^2_x}
$$
Taking the $L^q_t$ norm we and using the condition {\bf K1} we obtain
$$
\|F_k\|_{L^\infty_x} \les 2^k \|F_k\|_{L^q_t L^2_x} + 2^{-\frac 1{p}k} 
\|F_k\|_{L^{\frac {2qp}{2p-q}}_t L^2_x}+
2^{-\frac 1{p-1}k}
\|F_k\|_{L^{\frac {2q(p-1)}{2(p-1)-q}}_t L^2_x}
$$
Thus, applying  the dyadic Gagliardo-Nirenberg estimate of Lemma
\ref{lem:Gagl-Nir-L} we derive
\beaa
\|F_k\|_{L^\infty_x} &\les&\big ( 2^{k(\f12 -\frac 1q)}  + 2^{-\frac 1{p}k}\c 
2^{k(\frac 12 - \frac {2p-q}{2pq})}  + 
2^{-\frac 1{p-1}k}\c 2^{k(\frac 12 - \frac {2(p-1)-q}{2(p-1)q})}\big )\, \N_1[F_k]\,\\
&\les & 2^{k(\f12 -\frac 1q)}\c\big (1 + 2^{-\frac 1{2p}k}   + 
2^{-\frac 1{2(p-1)}k}\big )\,
\N_1[F_k]\\ &\les &  2^{k(\f12 -\frac 1q)}\, \N_1[F_k]
\eeaa
as desired.
\end{proof}
\section{Notations. Outline of the remainder of the paper}
In the proof we shall often refer to various properties 
of the LP calculus and other analytical tools developed above. To help 
the reader we give below a glossary of our main notations:
\begin{enumerate}
\item[]{ [{\bf H\"o}]}\quad -\quad H\"older inequality
\item[] {[{\bf Leib}]}\quad -\quad Leibnitz rule
\item[] {[{\bf GN}]}\quad--\quad Gagliardo-Nirenberg estimates in
\eqref{eq:GNirenberg}
\item[]{[{\bf Env}]}\quad -\quad envelope properties of definition 
\ref{def:Envelope}
\item[] {[ $\nab${\bf FB}]}\quad -\quad derivative finite band condition of v) of proposition
\ref{thm:geometric-LP}
 \item[] {[{\bf FB}$\nab$]}\quad -\quad dual derivative finite band condition of v) of proposition
\ref{thm:geometric-LP}
\item[] {[$\nab${\bf FB}$^{-1}$]}\quad -\quad inverse derivative finite band condition of v) of proposition
\ref{thm:geometric-LP}
\item[] {[$\De${\bf FB}]}\quad -\quad laplacean finite band condition of v) of proposition
\ref{thm:geometric-LP}
\item[] {[$\De${\bf FB}$^{-1}$]}\quad -\quad dual laplacean finite band condition of v) of proposition \ref{thm:geometric-LP}
\item[] {[{\bf wB}]}\quad -\quad weak Bernstein inequality of vi) of proposition
\ref{thm:geometric-LP}
\item[] {[{\bf wB}$^*$]}\quad -\quad dual weak Bernstein inequality of vi) of proposition
\ref{thm:geometric-LP}
\item[] {[{\bf ssB}]}\quad -\quad strong scalar Bernstein inequality of vii) of proposition
\ref{thm:geometric-LP}
\item[] {[{\bf ssB}$^*$]}\quad -\quad dual strong scalar Bernstein inequality of vii) of proposition
\ref{thm:geometric-LP}
\item[] {[{\bf GN}$_k$]}\quad -\quad dyadic Gagliardo-Nirenberg inequality of lemma
\ref{lem:Gagl-Nir-L}
\item[] {[\,,\,]}\quad -\quad commutator estimates of lemma \ref{lem:Comm-nabL}
\end{enumerate}

The remainder of the paper is focused on the proof of theorems \ref{thm:mainBesov2}
and \ref{thm:main-scalar}.

We start with section 7 where we prove an integration by parts lemma which is the 
noncommutative counterpart of lemma \ref{le:integr-parts-k}.

Section 8 contains the proof of the sharp bilinear trace estimate of 
theorem \ref{thm:main-scalar}. The proof follows the outline of the 
corresponding flat statement with modifications taking into the account
non-commutativity and absence of the strong tensor Bernstein inequality.

Sections 9 and 10 provide the proofs of the integrated sharp product estimates
I and II of theorems \ref{thm:main-scalar} and \ref{thm:mainBesov2}.

In section 11 we finish the proof of theorem  \ref{thm:main-scalar}
by establishing the $\BB^0$ estimates for solutions of a homogeneous
scalar transport equation. The proof is quickly reduced to an estimate 
for the commutator $[\nab_L, P_k]$ applied to a scalar $\BB^0$ function.

Finally, section 12 contains the proof of non-sharp tensor commutator 
estimates for $[\nab_L, P_k]$ needed to control various error terms throughout
the paper.
\section{Dyadic integration by parts}
In this section we shall deal with an integration by parts lemma which is needed 
in the proof of the sharp bilinear trace theorem of theorem \ref{thm:main-scalar}.
The lemma is the non-commutative analogue of the integration by parts argument in the 
proof of the flat sharp bilinear trace estimate of proposition \ref{prop:intr-2}.

We shall estimate the time integral of the expression 
$$
``\,\nab_{L}" \big (F_{k'}\c G_{k''}\big ):= (\nab_{L} F)_{k'}\c G_{k''}+F_{k'}\c (\nab_{L}
G)_{k''},
$$
where $F, G$ are tensors of the same order and $\cdot $ denotes the scalar product.
Observe that the expression above differs from the perfect derivative
$\frac d{ds} \big (F_{k'}\c G_{k''}\big )$ by the commutator terms 
$$
[P_{k'},\nab_{L}] F\c G_{k''},\qquad 
 F_{k'}\c [P_{k''},\nab_{L}] G
$$
\begin{proposition} Let   $F, G$  be tensors of same order. Then for any  $k,k', k''$,
\beaa
\|P_k\int_0^t ``\,\nab_{L}" \big (F_{k'}\c G_{k''}\big )\|_{L_t^\infty L_x^2}\les 2^{-\si(
|k'-k''|+|k'-k|)}\,\NN_1[F_{k'}]\c
\NN_1[G_{k''}]
\eeaa
with a strictly positive  $\si$ independent of $k,k', k''$.
\label{prop:integrbypartsPk}
\end{proposition}
\begin{proof}:\quad By symmetry it suffices to consider 
the following three cases:\qquad \qquad
\,\,\qquad \qquad {\bf a)}\,\, $k'\ge k''\ge k$, \qquad {\bf b)}\,\, $k'\ge k > k''$,
\qquad {\bf c)}\,\, $k>k'\ge k''$\newline
In all of the cases the proof of the proposition reduces to the estimate
for the commutators
\bea
\|P_k\int_0^t [P_{k'},\nab_{L}] F\c G_{k''}\|_{L_t^\infty L_x^2}  
&\les &2^{-\si( |k'-k''|+|k'-k|)}\,\NN_1[F_{k'}]\c
\NN_1[G_{k''}],
\label{eq:Commut1-Pk}\\ 
\|P_k\int_0^t F_{k'}\c [P_{k''},\nab_{L}] G\|_{L_t^\infty L_x^2}&\les & 
2^{-\si( |k'-k''|+|k'-k|)}\,\NN_1[F_{k'}]\c
\NN_1[G_{k''}]\label{eq:Commut2-Pk}
\eea
as well as the estimate for the boundary terms
\be{eq:Bound-Pk}
\|P_k \big ( F_{k'}(t)\c G_{k''}(t) - F_{k'}(0)\c G_{k''}(0)\big )\|_{L_t^\infty L_x^2} 
\les 2^{-\si( |k'-k''|+|k'-k|)}\,\NN_1[F_{k'}]\c \NN_1[G_{k''}]
\end{equation}
Since $k''\le k'$ the estimate \eqref{eq:Commut2-Pk} is
more sensitive  than \eqref{eq:Commut1-Pk}  and we shall only
prove it and \eqref{eq:Bound-Pk}
 in what follows.

We start with the easier of the three case.

{\bf a)}\,\, We start by applying the dual strong Bernstein inequality for scalars
followed by Cauchy-Schwartz with $q<2$, the dyadic Gagliardo-Nirenberg
inequality of Lemma \ref{lem:Gagl-Nir-L} 
and the commutator estimate of Lemma \ref{lem:Comm-nabL},
to infer that 
\beaa
\|P_k\int_0^t F_{k'}\c [P_{k''},\nab_{L}] G\|_{L_t^\infty L_x^2}&\les & 
2^k\|F_{k'}\c [P_{k''},\nab_{L}] G\|_{L_t^1 L_x^1}\qquad\qquad\qquad  [\mbox{{\bf ssB}}^*]\\
&\les &  2^k\|F_{k'}\|_{L^{q'}_t L^2_x} \| [P_{k''},\nab_{L}] G\|_{L_t^q L_x^2}
\qquad\qquad [\mbox{{\bf H\"o}}]\\
\text{{[\bf GN$_k$]}}\, \&\,[\,,\,]\qquad\qquad \qquad &\les & 2^k2^{-k'(\frac 12 +\frac 1{q'})}2^{-\frac {k''}2}
\N_1[F_{k'}]\c
\N_1[G]\\ &\les & 2^{-(\frac 12+\frac 1{q'}) (k'-k'')} 2^{-(k''-k)} \N_1[F_{k'}]\c 2^{-\frac
1{q'}k''}\N_1[G]\\ \text{{\bf [Env]}}\qquad\qquad\qquad&\les & 2^{-(\frac 12+\frac 1{q'})(k'-k'')} 2^{-(k''-k)}
\N_1[F_{k'}]\c
\N_1[G_{k''}]
\eeaa
The estimate for the boundary terms proceeds as follows.
Using the dual strong Bernstein inequality for scalars followed by
the Cauchy-Schwartz and the dyadic Gagliardo-Nirenberg estimate 
of Lemma \ref{lem:Gagl-Nir-L} with $q=\infty$, we obtain 
\beaa
\|P_k \big ( F_{k'}\c G_{k''}\big )\|_{L_t^\infty L_x^2} &\les & 
2^k \|F_{k'}\c G_{k''}\|_{L_t^\infty L_x^1}\qquad \qquad\qquad\qquad\,\,  [\mbox{{\bf ssB}}^*]\\
&\les &
 2^k \|F_{k'}\|_{L^\infty_t L^2_x}\| G_{k''}\|_{L_t^\infty L_x^2}\qquad \qquad\qquad
[\mbox{{\bf H\"o}}]\\ &\les & 2^{k-\frac {k'}2-\frac {k''}2} \N_1[F_{k'}]\c
\N_1[G_{k''}]\quad\qquad
\qquad\text{{[\bf GN$_k$]}}\, \\ &\les &  2^{-\frac 12(k'-k'')} 2^{-(k''-k)}
\N_1[F_{k'}]\c \N_1[G_{k''}]
\eeaa
{\bf b)}\quad We apply the H\"older inequality with $q<2$ followed by
the Gagliardo-Nirenebrg and weak Bernstein inequalities in the $x$ variable,
the dyadic Gagliardo-Nirenberg inequality of Lemma \ref{lem:Gagl-Nir-L},
and finally the commutator estimates of Lemma \ref{lem:Comm-nabL},
to infer that
\beaa
\|P_k\int_0^t F_{k'}\c [P_{k''},\nab_{L}] G\|_{L_t^\infty L_x^2}&\les & 
\|F_{k'}\|_{L^{q'}_t L^4_x}\| [P_{k''},\nab_{L}] G\|_{L_t^q L_x^4}\qquad \qquad\qquad
[\mbox{{\bf H\"o}}]\\
\text{{[\bf wB]}} \,\,\&\,\,  \text{{[\bf GN]}} \qquad\qquad
\qquad&\les & 
2^{\frac {k'}2}\|F_{k'}\|_{L^{q'}_t L^2_x} 
\| \nab [P_{k''},\nab_{L}] G\|_{L_t^q L_x^2}^\f12 
\| [P_{k''},\nab_{L}] G\|_{L_t^q L_x^2}^\f12\\
\text{{[\bf GN$_k$]}} \,\&\,\,[\,,\,]\qquad\qquad
\qquad &\les & 2^{-\frac 1{q'}k'}  \N_1[F_{k'}]\c \N_1[G]\\
&\les & 2^{-\frac 1{2q'} (k'-k'')} 2^{-\frac 1{2q'}(k-k'')} \N_1[F_{k'}]\c 2^{-\frac
1{q'}k''}\N_1[G]\\ 
\text{{\bf [Env]}}\qquad\qquad\qquad&\les & 2^{-\frac 1{2q'} (k'-k'')} 2^{-\frac
1{2q'}(k-k'')}\N_1[F_{k'}]\c
\N_1[G_{k''}]
\eeaa
The boundary terms are estimated with the help of the H\"older inequality
followed by the weak Bernstein inequality and the dyadic Gagliardo-Nirenberg 
estimate of Lemma \ref{lem:Gagl-Nir-L}.
\beaa
\|P_k \big ( F_{k'}\c G_{k''}\big )\|_{L_t^\infty L_x^2} &\les & 
\|F_{k'}\|_{L^\infty_t L^{\frac 83}_x}\| G_{k''}\|_{L_t^\infty L_x^8} \qquad \qquad\qquad
[\mbox{{\bf H\"o}}]\\
\text{{[\bf wB]}} \,\,\,\,   \qquad\qquad
\qquad &\les &
 2^{\frac 14 k'+\frac 34 k''} \|F_{k'}\|_{L^\infty_t L^2_x}\| G_{k''}\|_{L_t^\infty L_x^2} \\
\text{{[\bf GN$_k$]}} \,\qquad\qquad
\qquad &\les &
2^{-\frac 14 ({k'}-{k''})} \N_1[F_{k'}]\c \N_1[G_{k''}]\\ &\les & 
2^{-\frac 18(k'-k'')} 2^{-\frac 18(k-k'')} \N_1[F_{k'}]\c \N_1[G_{k''}]
\eeaa

{\bf c)}\,\,\, We start by applying the inverse finite band condition
followed by Lemma \ref{lem:Direct} and H\"older inequality 
with an exponent $1\le q<2$ chosen to be sufficiently close to $q=2$, 
\beaa
\|P_k \int_0^t F_{k'}\c [P_{k''},\nab_L] G\|_{L^\infty_t L^2_x} &\les& 
 2^{-k} \|\nab\int_0^t  F_{k'} \c [P_{k''},\ddd_L ]  G \|_{L^\infty_t L^2_x}\qquad 
\quad [{{\bf \nab FB^{-1}}}]\\
&\les &  2^{-k}\|\nab \big (F_{k'} \c [P_{k''},\ddd_L ] \,G\big )\|_{L^{1}_t L^2_x}\\
\,[\mbox{\bf {H\"o}}]\, \&\, [\mbox{\bf {Leib}}] \qquad\qquad\qquad 
&\les &  2^{-k}\|\nab F_{k'}\|_{L^{q'}_t L^4_x}  
\| [P_{k''},\ddd_L ] \,G\|_{L^{q}_t L^4_x} \\ &+&
 2^{-k}\|F_{k'}\|_{L^{q'}_t L^\infty_x}  \|\nab [P_{k''},\ddd_L ] \,G\|_{L^{q}_t L^2_x}
\\ &=& C^{(1)}_{kk'k''} +  C^{(2)}_{kk'k''}
\eeaa
To estimate the first term above we use the derivative Bernstein inequality of
Lemma \ref{lem:Gagl-Nir-nab4} together with the Gagliardo-Nirenberg estimate
followed by the commutator estimates of Lemma \ref{lem:Comm-nabL},
\beaa
C^{(1)}_{kk'k''}&\les & 
2^{-k} 2^{k'(\frac 12 -\frac 1{q'})} \N_1[F_{k'}] \c\| \nab [P_{k''},\ddd_L ] \,G\|_{L^{q}_t L^2_x}^\f12 \| [P_{k''},\ddd_L ] \,G\|_{L^{q}_t L^2_x}^\f12 \\
&\les &   2^{-k} 2^{k'(1-\frac 1{q'}} \N_1[F_{k'}] \c \N_1[G]\\
&\les & 2^{-\frac 1{2q'}(k'-k'')} 2^{-\frac 1{2q'} (k-k'')}\N_1[F_{k'}]  \c 2^{-\frac 1{q'}k''}\N_1[G]\\ &\les &  2^{-\frac 1{2q'}(k'-k'')} 2^{-\frac 1{2q'} (k-k'')
}\N_1[F_{k'}]\c  \N_1[G_{k''}]
\eeaa
To estimate the term $C^{(2)}_{kk'k''}$ we apply the integrated version of the
strong Bernstein inequality of Lemma \ref{lem:Int-SB-q} together with the 
commutator estimate of Lemma \ref{lem:Comm-nabL}. 
\beaa
C^{(2)}_{kk'k''}&\les & 2^{-k} 2^{k'(\frac 12 -\frac 1{q'})}
2^{-\frac {k''}2} \N_1[F_{k'}]\c \N_1[G] \\ &\les & 
2^{-\frac 1{q'}(k'-k'')} 2^{-\frac 14 (k-k'')} \N_1[F_{k'}]\c 2^{-k''(\frac 34 +\frac 1{q'})}
\N_1[G]\\ &\les &2^{-\frac 1{q'}(k'-k'')} 2^{-\frac 14 (k-k'')} \N_1[F_{k'}]\c
\N_1[G_{k''}]
\eeaa
To estimate the boundary terms for  the low-low interaction ($k''\le k'< k$) 
we argue as follows.
We start by applying the finite band property,
\beaa
 \|P_k(F_{k'}\c G_{k''})\|_{L^\infty_t L^2_x} &\les& 
  2^{-2 k} \|\Delta (F_{k'}\c G_{k''})\|_{L^\infty_t L^2_x}
\qquad \qquad\quad [{\bf \Delta FB^{-1}}]\\\,[\mbox{\bf {Leib}}]\qquad \qquad \qquad &\les &
  2^{-2k} \|\Delta F_{k'}\c G_{k''}\|_{L^\infty_t L^2_x} 
  + 2^{-2k} \|F_{k'}\c \Delta G_{k''}\|_{L^\infty_t L^2_x} \\ &+ &2^{-2k} 
  \|P_k(\nab F_{k'}\c\nab G_{k''})\|_{L^\infty_t L^2_x}\\ 
\, [{\bf \Delta FB}]\qquad \qquad \qquad &\les & 
2^{-2(k-k')} \|F_{k'}\c G_{k''}\|_{L^\infty_t L^2_x} 
  + 2^{-2(k-k'')} \|F_{k'}\c G_{k''}\|_{L^\infty_t L^2_x} \\ &+ &2^{-2k} 
  \|P_k(\nab F_{k'}\c\nab G_{k''})\|_{L^\infty_t L^2_x}
  \eeaa
  The most difficult is the last term. We apply the dual strong scalar Bernstein inequality
  followed by H\"older inequality, finite band property, and the dyadic Gagliardo-Nirenberg inequality of Lemma \ref{lem:Gagl-Nir-L},
  \beaa
   2^{-2k} 
  \|P_k(\nab F_{k'}\c\nab G_{k''})\|_{L^\infty_t L^2_x}& \les & 
  2^{-k} \|\nab F_{k'}\c\nab G_{k''}\|_{L^\infty_t L^1_x} \qquad \qquad 
[\mbox{\bf{ssB}}^*]\\ 
\, [\mbox{\bf{H\"o}}]\qquad \qquad \qquad \qquad &\les & 
   2^{-k} \|\nab F_{k'}\|_{L^\infty_t L^2_x}\|\nab G_{k''}\|_{L^\infty_t L^2_x}
   \\\, [\nab \mbox{\bf {FB}}]\qquad\qquad \qquad \qquad &\les & 2^{-k+k' +k''} 
   \|F_{k'}\|_{L^\infty_t L^2_x}\|G_{k''}\|_{L^\infty_t L^2_x}\\
\,[\mbox{\bf{GN}}_k]\qquad\qquad \qquad \qquad &\les & 
2^{-k+\frac{k'}2 +\frac{k''}2} \N_1[F_{k'}]\c \N_1[G_{k''}]\\
&\les & 2^{-\frac 14 (k'-k'')} 2^{-\frac 14(k-k'')}\N_1[F_{k'}]\c \N_1[G_{k''}]
  \eeaa

\end{proof}
\section{Proof of the sharp bilinear trace  theorem}
In this section we provide the proof of the sharp
 bilinear trace  estimate \eqref{eq:maintrace-scalar} of
theorem \ref{thm:main-scalar}:
\begin{equation}
\| \int_0^t \nab_L F\c G\|_{\BB^0}\les\N_1[F]\c
\N_1[G].\label{eq:transp-tr-inh-again}
\end{equation}
In fact it suffices to prove,
\begin{equation}
\sum_{k\ge 0}\|P_k \int_0^t \nab_L F\c G\|_{L_t^\infty L_x^2}\les\N_1[F]\c
\N_1[G].\label{eq:transp-tr-inh-k}
\end{equation}
Indeed for the low frequencies we can use the dual  strong  scalar Bernstein 
inequality [${\bf ssB}^*$],
\beaa
\| P_{<0} \int_0^t \nab_L F\c G\|_{L_t^\infty L_x^2}&\les&
 \|\int_0^t \nab_L F\c G\|_{L_t^\infty
L_x^1}\qquad\qquad \qquad [\text{\bf ssB}^*]\\
&\les&\|\nab_L F\|_{L_t^2L_x^2}\c \| G\|_{L_t^2L_x^2}\les \NN_1[F]\c\NN_1[G].
\eeaa
 We start with the LP -  decompositions of $\nab_L F\c  G$.
\beaa  (\ddd_LF\c G)
&=&(\ddd_L F)_{<k}\c G_{\ge k}+(\ddd_L F)_{\ge k}\c G_{<k}+(\ddd_L F)_{\ge k}\c
G_{\ge k}\\ &+&(\ddd_L F)_{<k}\c G_{<k}
\eeaa
Thus,
$$P_k \int_0^t \ddd_LF\c G =A_k+B_k+C_k+D_k$$
\bea
A_k&=& P_k \int_0^t (\ddd_L F)_{<k}\c G_{\ge k},\qquad\quad
B_k= P_k \int_0^t (\ddd_L F)_{\ge k}\c G_{< k}\label{eq:ABCDk}\\
C_k&=&P_k \int_0^t (\ddd_L F)_{<k}\c G_{< k},\qquad\quad
D_k=P_k \int_0^t (\ddd_L F)_{\ge k}\c G_{\ge  k}\nn
\eea
\subsection{{\bf {Estimates for $A_k=P_k\int_0^t (\ddd_L F)_{<k}\c G_{\ge k}$}}}
\qquad 

This is the easiest term. We start by giving a ``slightly wrong ''
proof based on the use of strong, tensorial form, of the Bernstein inequality 
``[{\bf stB}]":
\be{eq:fakeBernst}\| P_k F\|_{L_x^\infty}\les 2^k \|F\|_{L_x^2}
\end{equation}
which we don't in fact  possess. We shall indicate however how to circumvent this problem
in remark \ref{rem:correct-incorrection}.

 Indeed, using  \eqref{eq:fakeBernst}, 
\beaa
\|A_k\|_{L^\infty_t L^2_x}&\les&\sum_{k'<k\le k''}\int_0^{1}
\|(\ddd_L F)_{k'}\c G_{k''}\|_{L^2_x}\, dt\\\, [\text{\bf H\"o}]\qquad\qquad
&\les&\sum_{k'<k\le k''}\int_0^{1}
\|(\ddd_L F)_{k'}\|_{L^\infty_x} \|G_{k''}\|_{L^2_x}\,dt\\
\,\,``[\text{\bf stB}^*]"\qquad\qquad &\les&\sum_{k'<k\le k''}2^{k'}\int_0^{1}
\|(\ddd_L F)_{k'}\|_{L^2_x}\| G_{k''}\|_{L^2_x}\,dt\\\,
 [\nabla\text{\bf FB}^{-1}]\qquad\qquad &\les&  \,\sum_{k'<k\le k''}2^{k'-k''}
\|(\ddd_L F)_{k'}\|_{L^2_t L^2_x}\| \nab G_{k''}\|_{L^2_t L^2_x}
\eeaa
Therefore,
\beaa
\sum_{k} \|A_k\|_{L_t^\infty L^2_x}&\les& \sum_{k'<k\le k''}2^{k'-k''}
\|(\ddd_L F)_{k'}\|_{L^2_t L^2_x}\| \nab g_{k''}\|_{L^2_t L^2_x}\\
&\les & \sum_{k'< k''}2^{\frac {k'-k''}2}
\|(\ddd_L F)_{k'}\|_{L^2_t L^2_x}\| \nab G_{k''}\|_{L^2_t L^2_x}\les 
\N_1[F]\,\c \N_1[G].
\eeaa
\begin{remark}\label{rem:correct-incorrection}
It is easy to see that we don't need \eqref{eq:fakeBernst}; indeed we can replace it
with the weak Bernstein inequalities for $L^{p'}, L^{p''}$, 
 with ${p'}^{-1}+{p''}^{-1}=2^{-1}$ and $p''<<p'$,
 as follows 
\beaa
\|A_k\|_{L^\infty_t L^2_x}&\les&\sum_{k'<k\le k''}\int_0^{1}
\|(\ddd_L F)_{k'}\c G_{k''}\|_{L^2_x}\, dt\\
\, [\text{\bf H\"o}]\qquad\qquad &\les&\sum_{k'<k\le k''}\int_0^{1}
\|(\ddd_L F)_{k'}\|_{L^{p'}_x} \|G_{k''}\|_{L^{p''}_x}\,dt\\
\, [\text{\bf wB}]\qquad\qquad &\les&\sum_{k'<k\le k''}2^{k'(1-\frac{2}{p'})}2^{k''(1-\frac{2}{p''})}\int_0^{1}
\|(\ddd_L F)_{k'}\|_{L^2_x}\| G_{k''}\|_{L^2_x}\,dt\\
 \, [\nab\text{\bf FB}^{-1}]\qquad\qquad &\les& \sum_{k'<k\le k''}2^{(k'-k'')\c\frac{2}{p''}}
\|(\ddd_L F)_{k'}\|_{L^2_t L^2_x}\| \nab G_{k''}\|_{L^2_t L^2_x}
\eeaa
The proof then follows as before.
\end{remark}
\subsection{{\bf{Estimates for $D_k=P_k\int_0^t (\ddd_L F)_{\ge k}\c G_{\ge  k}$}}}
\qquad 

We write, $D_k=D_k^1+D_k^2$  where,
\begin{align*}
D_{k}^{1}=\sum_{k\le k'\le k''} P_{k} \int_0^t (\ddd_L F)_{k'}
\c G_{k''},\qquad
D_{k}^{2}= \sum_{k\le k'<k''} P_{k} \int_0^t(\ddd_L F)_{k''}\c G_{k'}
\end{align*}
 We only provide the proof for the term $D_{k}^{2}$ which
 is more difficult to treat since
the
$\nab_{L}$  derivative there falls on the factor with a higher frequency.
The term $D_k^1$ can be treated in the same way without the integration by parts.

We need to apply the integration by parts estimate of
proposition \ref{prop:integrbypartsPk} followed by the scalar dual strong 
 Bernstein inequality, H\"older, finite band  and the property of envelopes,
$$\NN_1[F_{m}]\les 2^{\ep|m-l|}\NN_1[F_{l}], \qquad \forall m, l\ge 0$$
 Thus, writing
$$(\nab_L F)_{k''}\c
G_{k'}=``\, \ddd_L "( F_{k''}\c G_{k'})  - F_{k''}\c
 (\nab_LG)_{k'}.  $$
we derive, for $\si>\ep$,
\beaa
\|D_{k}^{2}\|_{L_t^\infty L^{2}_x}&\les& 
 \sum_{k\le k'<k''} \|P_k\int_0^t F_{k''}\c
 (\nab_LG)_{k'}\|_{L_t^\infty L^{2}_x}
\qquad\qquad [\text{\bf prop. \ref{prop:integrbypartsPk}}]
\\ &+& \sum_{k\le k'<k''} 2^{-\si\big(
|k'-k''|+|k-k'|\big)}
\NN_1[F_{k''}]\c \NN_1[G_{k'}]\\
\,[\text{\bf ssB}] \, \text{\&}\,  [\text{\bf Env}]\qquad &\les& 2^k\sum_{k\les k'<k''} 
\|\int_0^t F_{k''}\c
(\nab_{L} G)_{k'}\|_{L_t^\infty L^{1}_x} +\NN_1[F_{k}]\c\NN_1[G_k]\\
\,[\text{\bf H\"o}]\qquad\qquad&\les&
 2^k\sum_{k\les k'<k''} \|F_{k''}\|_{L_t^2L_x^2}
 \|(\nab_LG)_{k'}\|_{L_t^2L_x^2}
 +\NN_1[F_{k}]\c\NN_1[G_k]\\
\, [\text{\bf Env}]\qquad\qquad &\les& \sum_{k\les k'<k''}2^{k-k''}\NN_1[F_{k''}]\c\NN_1[G_{k'}]
+\NN_1[F_{k}]\c\NN_1[G_k]\les \NN_1[F_{k}]\c\NN_1[G_k]
\eeaa
Thus, 
$$
\sum_{k} \|D_{k}^{2}\|_{L^\infty_tL_x^2} \les 
\sum_{k} \NN_1[F_{k}]\c\NN_1[G_k]
\les \N_1[F]\, \c\N_1[G]
$$

\subsection{{\bf{Estimates for 
$B_k=P_k \int_0^t (\ddd_L F)_{\ge k}\c G_{< k} $}}}\qquad 

We start by decomposing,
\beaa
B_k=\sum_{k'<k\le k''}P_k\int_0^t (\ddd_L F)_{k''}\c G_{k'}
\eeaa
Integrating by parts with the help of proposition \ref{prop:integrbypartsPk}
we obtain,
\beaa
\|B_k\|_{L_t^\infty L^2_x}&\les& 
\sum_{k'<k\le k''}\bigg(
\|\int_{\gat}   F_{k''}\c (\nab_{L} G)_{k'}\|_{L_t^\infty L^{2}_x}+
2^{-\si\big( |k'-k''|+|k-k'|\big)}
\NN_1[F_{k''}]\c \NN_1[G_{k'}]\bigg)\\
&\les&\sum_{k'<k\le k''}
\|\int_0^t  F_{k''}\c (\nab_{L} G)_{k'}\|_{L_t^\infty
L^{2}_x}+\NN_1[F_{k}]\c\NN_1[G_{k}]
\eeaa
To estimate the sum on the right hand side we
proceed exactly as in the proof for 
$A_{k}$.

{\bf Estimates for 
$C_k=P_k \int_0^t (\ddd_L F)_{<k}\c G_{< k}$}\qquad 
This term, which is absent 
in the classical paradifferential  calculus, 
 is by far the most difficult and requires a lot
more work than the previous ones.
We further decompose 
\begin{align*}
&C_k=C_k^{(1)} + C_k^{(2)},\\
&C_k^{(1)} = \sum_{k'\le k''<k} P_k \int_0^t (\ddd_L F)_{k'}\c G_{k''},\\
&C_k^{(2)} = \sum_{k'<k''<k} P_k \int_0^t(\ddd_L F)_{k''}\c G_{k'}
\end{align*}
As before we ignore the term $C_k^{(1)}$ since
 the term $C_k^{(2)}$ is  clearly the more difficult one
as the 
$\nab_L$ derivative  falls on a term with higher frequency
 and thus requires an integration by parts.

We first integrate by parts with the help 
of proposition \ref{prop:integrbypartsPk}.
\beaa
\|C_k^{(2)}\|_{L_t^\infty L^2_x} &\les & \sum_{k'<k''<k}
 \|P_k \int_0^t F_{k''}\c (\nab_L G)_{k'}\|_{L_t^\infty L^2_x} +
\NN_1[F_k]\c\NN_1[G_k]\\ 
\eeaa
To control the first term above we need 
to use first  the finite band condition $[\bf{\nab\text{FB}}^{-1}]$
followed by Lemma \ref{lem:Direct} and the Leibnitz rule
\beaa
 \|P_k \int_0^t  F_{k''}\c (\nab_L G)_{k'}\|_{L_t^\infty L^2_x}&\les &2^{-k}\|\nab
\int_0^t F_{k''}\c (\nab_L G)_{k'}\|_{L^\infty_t
L^2_x}\qquad\qquad\qquad[\bf{\nab\text{FB}}^{-1}]\\
\textbf{[Le \,\,\ref{lem:Direct}]}\qquad\qquad\qquad &\les & 2^{-k}\big( \|\nab\big (
F_{k''}\c (\nab_L G)_{k'}\big )\|_{L^1_t L^2_x}+\|
F_{k''}\c (\nab_L G)_{k'}\|_{L^1_t L^2_x}\big)\\ 
\,[\textbf{Leib}]\qquad\qquad\qquad&\les & 2^{-k} \|\nab F_{k''}\c (\nab_L G)_{k'}\|_{L^1_t L^2_x} + 2^{-k} \| F_{k''}\c \nab (\nab_L G)_{k'}\|_{L^1_t L^2_x}\\&+&2^{-k}
\| F_{k''}\c  (\nab_L G)_{k'}\|_{L^1_t L^2_x}
\eeaa
 To continue it suffices to consider only  the first term since in that case the derivative
falls on a term with higher frequency.
Using H\"older followed by the  inequality \eqref{eq:trickyBernstein} of Lemma
\ref{lem:Gagl-Nir-nab4} and the weak Bernstein inequality, 
we infer that,
\beaa
\|\nab F_{k''}\c (\nab_L G)_{k'}\|_{L^1_t L^2_x}&\les & 
\|\nab F_{k''}\|_{L^2_t L^4_x} \|(\nab_L G)_{k'}\|_{L^2_t L^4_x}
\qquad\qquad\, [\text{\bf H\"o}]\\
\, [\textbf{\eqref{eq:trickyBernstein}}]\qquad\qquad\qquad 
&\les & 2^{\frac {k''+k'}2}\N_1[F_{k''}] \c
\NN_1[ G_{k'}]
\eeaa
This leads to the estimate 
\beaa
 \sum_k \sum_{k'<k''<k}\|P_k \int_0^t  
F_{k''}\c (\nab_L G)_{k'}\|_{L^\infty_t L^2_x}&\les &
\sum_k \sum_{k'<k''<k}2^{\frac {k''+k'}2-k}\N_1[F_{k''}]\,\c
\N_1[G_{k'}]\\ &\les &\N_1[F]\,\c \N_1[G]
\eeaa
\section{ Integrated  sharp product estimate I}
In this section we prove
  the integrated sharp product estimate \eqref{eq:mainBesov-scalar}
of theorem \ref{thm:main-scalar}. We have already shown that it suffices
to prove the estimate only in a fully scalar case. Namely, that for two 
scalar functions $f$ and $g$,   
\be{eq:mainBesov-scalar-again}
\|\int_0^t  f\c g\|_{\BB^0}\les\|f\|_{\PP^0}\c\big(\NN_1[g]+\|
g\|_{L_x^\infty L_t^2}\big)
\end{equation}
Observe that it suffices to prove,
\be{eq:mainBesov-scalar-k}
\sum_{k\ge 0}\|P_k\int_0^t  f\c g\|_{L_t^\infty L_x^2}\les\|f\|_{\PP^0}\c\big(\NN_1[g]+\|
g\|_{L_x^\infty L_t^2}\big)
\end{equation}
\begin{remark}
The proof below will show that we have, in fact, a stronger dyadic estimate:
for all $k\ge 0$,
$$
\|P_k\int_0^t  f\c g\|_{L_t^\infty L_x^2}\les\big(\NN_1[g]+\|
g\|_{L_x^\infty L_t^2}\big)\c\Big (\sum_{k'}2^{-|k-k'|}\|P_{k'}f\|_{L^2_t L^2_x}
+ 2^{-k}\|f\|_{L^2_t L^2_x}\Big ),
$$
needed to establish the second part of theorem \ref{thm:mainBesov}.
\end{remark}
We  start by applying the LP-decomposition to $f\c g$.
$$ 
f\c g
= f_{<k}\c g + f_{\ge k}\c g
$$
Thus, 
$$
P_k\big(\int_0^t  f\c g\big)=A_k+B_k,
$$ where 
\beaa
A_k =P_k\int_0^t f_{<k}\c g,\qquad
B_k=  P_k\int_0^t  f_{\ge k}\c g
\eeaa
\subsection {Estimates for 
$A_k=P_k\int_0^t f_{<k}\c g$.}\quad

Applying the dual finite band condition
 followed by  \eqref{eq:nabcom1} of Lemma \ref{lem:Direct} and 
the Leibnitz rule we obtain
\beaa
\|A_k\|_{L^\infty_t L^2_x}&\les &2^{-k} 
\|\nab \int_0^t f_{<k}\c g\|_{L^\infty_t L^2_x}\qquad \qquad \qquad
[\nab\mbox {\bf {FB}}^{-1}]\\
&\les &
2^{-k}\|\nab (f_{<k}\c g)\|_{ L^2_xL^1_t} \\ \,[\mbox{\bf{Leib}}]\qquad
&\les & 2^{-k}\|\nab f_{<k}\c g\|_{L^2_x L^1_t} 
+ 2^{-k} \|f_{<k}\c \nab g\|_{L^2_x L^1_t} \\ 
\,[\mbox{\bf{H\"o}}]\qquad&\les & 2^{-k}\|\nab f_{<k}\|_{L^2_x L^2_t} \| g\|_{L^\infty_x L^2_t} 
+ 2^{-k}\|f_{<k}\|_{L^2_t L^\infty_x} \|\nab g\|_{L^2_t L^2_x }  \\
\,[\text{\bf FB}]\,\&\, [\mbox{\bf {ssB}}]
\qquad &\les &  
 \sum_{k'<k}  2^{k'-k}\|f_{k'}\|_{L^2_t L^2_x} \big (\| g\|_{L^\infty_x L^2_t} + 
\|\nab g\|_{L^2_t L^2_x } \big )
\eeaa
Therefore,
\beaa
\sum_k \|A_k\|_{L^\infty_t L^2_x}&\les & 
\big (\| g\|_{L^\infty_x L^2_t} + 
\|\nab g\|_{L^2_t L^2_x } \big )
\sum_k \sum_{k'<k}  2^{k'-k}\|f_{k'}\|_{L^2_t L^2_x} \\
&\les & \big (\| g\|_{L^\infty_x L^2_t} + 
\|\nab g\|_{L^2_t L^2_x } \big )\sum_{k'}  \|f_{k'}\|_{L^2_t L^2_x} 
\\&\les & \|f\|_{\PP^0}\c \big (\| g\|_{L^\infty_x L^2_t} + 
\|\nab g\|_{L^2_t L^2_x } \big )
\eeaa
\subsection {Estimates for $B_k=  P_k\int_0^t  f_{\ge k}\c g$}\quad
We decompose $f_{\ge k} = \sum_{k'\ge k} f_{k'}$ and on each $k''$
dyadic LP-piece apply the finite band condition 
$f_{k'} \approx 2^{-2k'} \Delta f_{k'}$.
\beaa
\|B_k\|_{L^\infty_t L^2_x} &\les & \| P_k \int_0^t f_{\ge k}\c g\|_{L^\infty_t L^2_x}\les 
 \sum_{k' \ge k}\| P_k \int_0^t f_{k'}\c g\|_{L^\infty_t L^2_x}\\
\,[\De{\bf FB}^{-1}]\qquad\qquad  &\les &  \sum_{k' \ge k} 2^{-2k'}
\| P_k \int_0^t  \Delta f_{k'}\c   g\|_{L^\infty_t L^2_x}\\
&\les &  \sum_{k' \ge k} 2^{-2k'}
\| P_k \int_0^t  \div\big (\nab f_{k'}\c  g\big )\|_{L^\infty_t L^2_x}
\\ &+& \sum_{k' \ge k} 2^{-2k'}\| P_k \int_0^t \nab f_{k'}\c  \nab g\|_{L^\infty_t L^2_x}=  B_k^{(1)} + B_k^{(2)}
\eeaa
To estimate $B_k^{(1)}$ we use a result of Lemma \ref{lem:Reverse} according 
to which for any vector-field $F$
\be{eq:represent}
\int_0^t \div F = \div W+ E,
\end{equation}
where $W$ is a solution of the transport equation 
$$
\nab_L W -\chi\c W= F,\quad W|_{S_0}=0
$$
and the error term $E$ satisfies the estimate 
\be{eq:E-q}
\|E\|_{L^{\frac {2p}{2-p}}_x L^\infty_t}\les \|F\|_{L^p_x L^1_t}
\end{equation}
for any $1\le p \le 2$.
Note that for any $r\ge 1$ we have the following estimate for $W$:
\be{eq:W-q}
\|W\|_{L^r_xL^\infty_t}\les \|F\|_{L^r_x L^1_t}
\end{equation}
Therefore we use the representation \eqref{eq:represent} with $F = \nab f_{k'}\c   g$,  
apply the dual finite band condition
together with strong scalar Bernstein inequality followed by the estimates
 \eqref{eq:E-q} and \eqref{eq:W-q} with $p=1$ and $r=2$, and finite band 
 condition. 
\beaa
B_k^{(1)} &\les& \sum_{k' \ge k} 2^{-2k'} \|P_k \div W\|_{L^\infty_t L^2_x}
+ \sum_{k' \ge k} 2^{-2k'} \|P_k E\|_{L^\infty_t L^2_x}\\
\, [\text{\bf FB}\nab]\, \&\, [\text{\bf ssB}^*]\qquad
&\les &\sum_{k' \ge k} 2^{k-2k'} \|W\|_{L^\infty_t L^2_x}
+ \sum_{k' \ge k} 2^{k-2k'} \|E\|_{L^\infty_t L^1_x}\\
&\les & \sum_{k' \ge k} 2^{k-2k'} \|\nab f_{k'}\c g\|_{ L^2_x L^1_t }
+ \sum_{k' \ge k} 2^{k-2k'} \|\nab f_{k'}\c g\|_{ L^2_x L^1_t }\\
&\les & \sum_{k' \ge k} 2^{k-2k'} \|\nab f_{k'}\|_{L^2_x L^2_t}
\| g\|_{ L^\infty_x L^2_t }\\ \, [\nab\mbox{\bf{FB}}]\qquad &\les & 
\sum_{k'' \ge k} 2^{k-k'} \| f_{k'}\|_{L^2_t L^2_x} \|g\|_{L^\infty_xL^2_t}
\eeaa
Thus,
$$
\sum_k B_k^{(1)} \les 
\sum_k \sum_{k' \ge k} 2^{k-k'} \| f_{k'}\|_{L^2_t L^2_x} 
\|g\|_{L^\infty_x L^2_t}\les  \| f\|_{\PP^0} \|g\|_{L^\infty_x L^2_t}
$$
We now  estimate the term $B_k^{(2)}$.
 Using the dual strong scalar Bernstein inequality followed by H\"older inequality and
 the finite band property, we obtain
\beaa
B_k^{(2)}&\les & \sum_{k' \ge k} 2^{k-2k'}\|\int_0^t \nab f_{k'}\c  \nab g\|_{L^\infty_t L^1_x}
\qquad \qquad [\mbox{\bf{ssB}}^*]\\ 
&\les & \sum_{k' \ge k}  2^{k-2k'} \|\nab f_{k'}\|_{L^2_t L^2_x} 
\|\nab g\|_{L^2_t L^2_x}\\ \,[\nab\mbox{\bf {FB}}]\qquad
&\les & \sum_{k' \ge k}  2^{k-k'} \|f_{k'}\|_{L^2_t L^2_x} 
\|g\|_{L^2_t L^2_x}
\eeaa
Therefore, as in the case of $B_k^{(1)}$,
$$
\sum_k B_k^{(2)} \les \| f\|_{\PP^0} \|g\|_{L^\infty_x L^2_t}
$$
  It then follows that
  $$
 \sum_k \|B_k\|_{L^\infty_t L^2_x}\les 
\| f\|_{\PP^0} \|g\|_{L^\infty_x L^2_t}
  $$

\section{Integrated sharp product estimate II}
In this section we indicate how to prove   the sharp  product estimate of   theorem
\ref{thm:mainBesov2}, which has already been reduced to the fully scalar case.
  For  scalar functions  $w, g$. where $w$ is a solution of the transport equation
  $$
  \nab_L w=f,\qquad f|_{S_0}=f_0,
  $$ we need to prove the estimate
\bea
\| w\c g\, \|_{P^{0}}\les  
\big (\|f_0\|_{B^{0}_{2,1}} + \| f\|_{\PP^0} \big )\big (\N_1[g] + \|g\|_{L^\infty_x L^2_t}\big
) 
\label{eq:main2-P}
\eea
\begin{proof}:\quad 
  We define the  functions $w_{(<k)}, w_{(k)}, w_{(\ge k)}$ as solutions 
  of the transport equations 
 \begin{align*}
& \nab_L w_{(<k)} = f_{<k},\qquad w_{(<k)}|_{S_0} = {f_0}_{<k},\\
 & \nab_L w_{(k)} = f_{<k},\qquad w_{(k)}|_{S_0} = {f_0}_{k},\\
 & \nab_L w_{(\ge k )} = f_{<k},\qquad w_{(\ge k)}|_{S_0} = {f_0}_{\ge k},
 \end{align*}
 Observe that according to the results of lemmas \ref{le:transportformulaL2} and \ref{lem:Direct}
 we have 
 \bea
 \|w_{([s]k)}\|_{L^p_x L^\infty_t} &\les &
\|{f_0}_{[s]k}\|_{L^p_x}+  \|f_{[s]k}\|_{L^p_x L^1_t}, \qquad \label{eq:LP-w}\\
 \|\nab w_{([s]k)}\|_{L^p_x L^\infty_t} &\les &
 \|\nab {f_0}_{[s]k}\|_{L^p_x}+ \|\nab f_{[s]k}\|_{L^p_x L^1_t} \label{eq:LP-w1}
 \eea
where $[s]k=<k, k, \ge k$.
We now consider 
the following LP-decomposition of $F\c G$:
$$ 
w\c g
=  w_{(<k)}\c g + w_{(\ge k)}\c g
$$
Thus,
$
P_k\big ( w\c g\big)=A_k+B_k,
$ where 
\beaa
A_k =P_k( w_{(<k)}\c g),\qquad
B_k=  P_k ( w_{(\ge k )}\c f).
\eeaa

 \subsection{Estimates for $A_k =P_k( w_{(<k)}\c g)$}\quad
 
 We have
 \beaa
\|A_k\|_{L^2_t L^2_x}&\les &2^{-k} 
\|\nab (w_{(<k)}\c g)\|_{L^2_t L^2_x}\qquad \qquad \qquad
[\nab\mbox {\bf {FB}}^{-1}]\\ \,[\mbox{\bf{Leib}}]\qquad
&\les & 2^{-k}\|\nab w_{(<k)}\c g\|_{L^2_x L^2_t} 
+ 2^{-k} \|w_{(<k)}\c \nab g\|_{L^2_x L^2_t} \\ 
\,[\mbox{\bf{H\"o}}]\qquad&\les & 2^{-k}\|\nab w_{(<k)}\|_{L^2_x L^\infty_t} \| g\|_{L^\infty_x L^2_t} 
+ 2^{-k}\|w_{(<k)}\|_{L^\infty_t L^\infty_x} \|\nab g\|_{L^2_t L^2_x }  \\
 &\les & 2^{-k}\sum_{k'<k}\|\nab w_{(k')}\|_{L^2_x L^\infty_t} \| g\|_{L^\infty_x L^2_t} 
+ 2^{-k}\|w_{(k')}\|_{L^\infty_t L^\infty_x} \|\nab g\|_{L^2_t L^2_x }\\\,
[\eqref{eq:LP-w}]\qquad  &\les &
2^{-k}\sum_{k'<k}\big (\|\nab {f_0}_{k'}\|_{L^2_x} + 
\|\nab f_{k'}\|_{L^2_x L^1_t}\big ) 
\| g\|_{L^\infty_x L^2_t} \\
&+& 2^{-k}
\big (\|{f_0}_{k'}\|_{L^\infty_x} + 
\|f_{k'}\|_{L^\infty_x L^1_t}\big )  \|\nab g\|_{L^2_t L^2_x }\\ 
\,[\nab\text{\bf FB}]\,\&\, [\mbox{\bf {ssB}}]
\qquad &\les &  
 \sum_{k'<k}  2^{k'-k}\big (\|{f_0}_{k'}\|_{L^2_x} +\|f_{k'}\|_{L^2_t L^2_x} \big )
\big (\| g\|_{L^\infty_x L^2_t} + 
\|\nab g\|_{L^2_t L^2_x } \big )
\eeaa
Therefore,
\beaa
\sum_k \les \|A_k\|_{L^\infty_t L^2_x}\les  
\big (\|f_0\|_{B^0_{2,1}}+ \|f\|_{\PP^0}\big )\c \big (\| g\|_{L^\infty_x L^2_t} + 
\|\nab g\|_{L^2_t L^2_x } \big )
\eeaa
\subsection {Estimates for $B_k=  P_k (w_{(\ge k)}\c g)$}\quad

We define the functions $[\Delta w]_{(k')}$ as solutions of the transport
equations 
$$
\nab_L [\Delta w]_{(k')}= \Delta f_{k'},\qquad 
 [\Delta w]_{(k')}|_{S_0} = \Delta {f_0}_{k'}
$$
Observe that $ \Delta f_{k'}\approx 2^{2k'} f_{k'}$
and therefore $[\Delta w]_{(k')}\approx  2^{2k'} w_{(k')}$.
In view of this,
\beaa
\|B_k\|_{L^2_t L^2_x} \les  \sum_{k'\ge k} 2^{-2k'}
\| P_k \big ( [\Delta w]_{ (k')}\c g\big )\|_{L^2_t L^2_x}
\eeaa
Now, according to lemma \ref{lem:Reverse} for $p=1$, we can write
$$[\Delta w]_{ (k')}=\div W_{(k')}+E_{(k')}$$
where,
\beaa
\ddd_L  W_{(k')}-\chi \c W_{(k')}&=&\nab f_{k'}, \quad W_{(k')}|_{S_0}=\nab (f_0)_{k'}\\
\|E_{(k')}\|_{L_x^2L_t^\infty}&\les&
 \|\nab (f_0)_{k'}\|_{L_2(S_0)}+\|\nab f_{k'}\|_{L_x^2 L_t^1}
\eeaa
The proof continues essentially as in the proof of the estimates for $B_k$
in the previous section.
\end{proof}

 \section{Sharp scalar commutator estimates}
The goal of this section is to
prove the last part of theorem \ref{thm:main-scalar}. In other
words
 we have to show that any solution $f$ of a homogeneous
scalar transport equation 
$$
\nab_L f = 0,\qquad f|_{S_0}=f_0
$$
verifies the estimate
\be{eq:transp-comm}
\|f\|_{\BB^0}\les \|f_0\|_{B^0_{2,1}(S_0)} 
\end{equation}
We first observe that it is trivial to take care of the low frequency
component $\|P_{<0} f\|_{L_t^\infty L_x^2}$ of $\|f\|_{\BB^0}$.
Indeed, $\|P_{<0} f\|_{L_t^\infty L_x^2}\les \| f\|_{L_t^\infty L_x^2}$
and clearly $\|f\|_{L_t^\infty L_x^2}\le \|f_0\|_{L^2}\les \|f_0\|_{B^0_{2,1}(S_0)}$.
Therefore we only have to show that,
\be{eq:transp-comm-high}
\sum_{k\ge 0}\|P_kf\|_{L_t^\infty L_x^2}\les \|f_0\|_{B^0_{2,1}(S_0)} 
\end{equation}
Commuting the equation with the LP-projection
we obtain
$$
\nab_L P_k f = [\nab_L , P_k] f,\qquad (P_k f)|_{S_0}=P_k f_0
$$
Thus,
$P_k f=P_k f_0+\int_0^t [\nab_L , P_k] f$ and, 
\beaa
\|P_k f(t)\|_{L_t^\infty L_x^2}\les \|P_k f_0\|_{L^2(S_0)}+\|\int_0^t[\nab_L , P_k]f\|_{L_t^\infty
L_x^2}
\eeaa
\begin{remark}During  the argument below we need  to keep track of the intervals
of integration in $t$. For this reason we denote by  
$\|f(t)\|_{L_t^qL_x^p}$
spacetime norms in the interval $[0,t]$. Often though, when
no confusion
is possible, we will drop the explicit dependence on $t$.
\end{remark}
Consequently,
$$\sum_{k\ge 0}\|P_k f(t)\|_{L_t^\infty L_x^2}\les \|f_0\|_{B^0_{2,1}(S_0)} +
\sum_{k\ge 0}\|\int_0^t[\nab_L , P_k]f(t)\|_{L_t^\infty L_x^2}$$
In the proposition below
we shall prove the following inequality:
\beaa
\sum_{k\ge 0}\|\int_0^t[\nab_L , P_k]f(t)\|_{L_t^\infty L_x^2}\les \De_0  \sum_{k\ge 0}\|P_k
f(t)\|_{L_t^\infty L_x^2}+\int_0^t\sum_{k\ge 0}\|P_k f(s)\|_{L_t^\infty L_x^2}ds.
\eeaa
Thus, for small $\De_0$,
$$\sum_{k\ge 0}\|P_k f(t)\|_{L_t^\infty L_x^2}\les \|f_0\|_{B^0_{2,1}(S_0)}+\int_0^t\sum_{k\ge 0}\|P_k
f(s)\|_{L_t^\infty L_x^2} ds.
$$
and the desired estimate follows by a straightforward Gronwall inequality.
Therefore to prove \eqref{eq:transp-comm-high} it  suffices to establish  the following:
\begin{proposition}\label{prop:Sharp-Comm}
Let $f$ be a scalar function on ${\cal H}$. Then 
\be{eq:sharp-comm}
\sum_{k\ge 0} \| \int_0^t [\nab_L, P_k] f\|_{L_t^\infty L^2_x} \le \epsilon
\|f\|_{\BB^0}+\int_0^t\|f(s)\|_{B^0}\, ds
\end{equation}
where $\|f(s)\|_{\BB^0}= \sum_{k\ge 0} \|P_k f(s)\|_{L_t^\infty L_x^2} +\|P_{<0} f(s)\|_{L_t^\infty
L_x^2} $ in the spirit of the remark above. 
\end{proposition}
\begin{proof}:\qquad
Since $P_k =\int_0^\infty m_k(\tau) U(\tau) f$
we have,
$$ \int_0^t [\nab_L, P_k] f= \int_0^\infty  m_k(\tau) d\tau\int_0^t [\nab_L, U(\tau)] fds$$
or, introducing the notation
\be{eq:def-Ltau}
\L(\tau) f =\int_0^t [\nab_L, U(\tau)] f,
\end{equation}
we have,
\bea
 \int_0^t [\nab_L, P_k] f&=& \int_0^\infty  m_k(\tau)\L(\tau)f d\tau\label{eq:maincomm-Ltau}\\
&=&\int_0^1  m_k(\tau)\L(\tau)f d\tau +\int_1^\infty  m_k(\tau)\L(\tau)f d\tau\nn
\eea
Observe that the integral $\int_1^\infty  m_k(\tau)\L(\tau)f d\tau$
is a lot easier to estimate. In fact all we need to treat it
is to show that $\L(\tau)f$ verifies the estimate,
\be{eq:tau1-infty}
\|\L(\tau)f\|_{L_t^\infty L_x^2}\les \tau^ K \|f\|_{L_t^2 L_x^2}, \qquad \mbox{for}\quad
\tau\ge 1
\end{equation}
for some positive value of $K$.
Indeed \eqref{eq:tau1-infty} easily implies the better estimate
$$\sum_k\|\int_1^\infty  m_k(\tau)
\L(\tau)f d\tau\|_{L_t^\infty L_x^2}\les \|f\|_{L_t^2L_x^2}.$$
In what follows we shall only concentrate
on the more difficult term $\int_0^1  m_k(\tau)\L(\tau)f d\tau$
and ignore the contribution of $\tau\ge 1$.

We  start with  a simple  property  of a $\BB^0$ function in terms of the heat 
semigroup $U(\tau)$.
\begin{lemma}\label{lem:Besov-heat}
For any smooth scalar function $f$ on ${\cal H}$
$$
\int_0^1 \|\nab^2 U(\tau) f\|_{L^\infty_t L^2_x}\,d\tau \les 
\|f\|_{\BB^0}
$$
\end{lemma} 
\begin{proof}:\qquad
First,  in view of the scalar  Bochner 
inequality \eqref{eq:conseq-Bochner-ineq-scalar},
$$
\|\nab^2 U(\tau) f\|_{L^2_x} \les \|\Delta U(\tau) f\|_{L^2_x}
$$
We now decompose
$$
\|\Delta U(\tau) f\|_{L^2_x}\les \sum_{k\ge 0} 
\|\Delta P_k U(\tau) f\|_{L^2_x}+\|f\|_{L^2_x}.
$$
We note  the following,
\be{eq:altern-bound-U}
\|\De P_k U(\tau) f\|_{L^2_x}\les \min(2^{2k}, 2^{-2k} \tau^{-2})\|P_k f\|_{L^2_x}
\end{equation}
Indeed 
we have both, 
\begin{align}
& \|\Delta P_k U(\tau) f\|_{L^2_x}\les 2^{2k} \|P_k  f\|_{L^2_x},\label{eq:compete1}\\
 &\|\Delta P_k U(\tau) f\|_{L^2_x}\les \tau^{-1} \|P_k  U(\tau/2)f\|_{L^2_x}\les
 2^{-2k} \tau^{-2} \|P_k f\|_{L^2_x}\label{eq:compete2}
\end{align}
Therefore,
\beaa
&&\int_0^1 \|\nab^2 U(\tau) f\|_{L^\infty_t L^2_x}\,d\tau \les  
\sum_{k\ge 0} \int_0^1\|\Delta P_k U(\tau) f\|_{L^\infty_t L^2_x}\, d\tau +\|f\|_{L^2_x}\\
&&\les 
\sum_{k\ge 0} \big(\int_0^{2^{-2k}}2^{2k}\| P_k f\|_{L^\infty_t L^2_x}\,d\tau +
 \int_{2^{-2k}}^1  2^{-2k} \tau^{-2}\|P_k f\|_{L^\infty_t L^2_x} d\tau\big)\\
&&\les  \sum_{k\ge 0} \| P_k f\|_{L^\infty_t L^2_x}+\|f\|_{L^2_x} \les\|f\|_{\BB^0}
\eeaa
as desired.
\end{proof}
The main step in proving
the proposition
is  the following:
\begin{lemma}The following estimates hold:
\be{eq:Ltau-bound-trivial}
\sup_{0\le \tau\le 1}\|\L(\tau ) f(t)\|_{L^\infty_t L^2_x} \les \|f\|_{\BB^0}
\end{equation}
Moreover for any $\tau\approx 2^{-2m_0}$,
\be{eq:Ltau-bound}
\|\L(\tau ) f(t)\|_{L^\infty_t L^2_x} \les  \sum_{k\ge 0} 2^{-\eps |m_0-k|}
\big( \De_0\c\|P_k f(t)\|_{L^\infty_t L^2_x} + 
 \int_0^t\|P_k f(s)\|_{ L^2_x}\big)
\end{equation}
\label{le:main-comm}
\end{lemma}
We postpone the proof
of the lemma and show now 
 how  \eqref{eq:Ltau-bound} 
 implies  proposition \ref{prop:Sharp-Comm}.
 Indeed  observe that  estimate \eqref{eq:Ltau-bound} is equivalent
to the bound,
$$
\|\L(\tau) f(t)\|_{L^\infty_t L^2_x} \les \sum_m
\min \big(2^m\tau^\f12, 2^{-m}\tau^{-\f12}\big)^\eps \big (\De_0\|P_m f(t)\|_{L^\infty_t L^2_x}
+ \int_0^t\|P_m f(s)\|_{L^2_x}ds\big )
$$
Therefore, in view of \eqref{eq:maincomm-Ltau},
\beaa
&&\|\int_0^t [\nab_L, P_k] \|_{L^\infty_t L^2_x}\les 
\int_0^1 m_k(\tau)\|\L(\tau) f(t)\|_{L^\infty_tL^2_x} \\ &\les & \sum_m
\int_0^1 m_k(\tau)
\min \big(2^m\tau^\f12, 2^{-m}\tau^{-\f12}\big)^\eps\, 
\big (\De_0 \|P_m f(t)\|_{L^\infty_t L^2_x} + \int_0^t \|P_m f\|_{L^2_x}\big )\\
& = & \sum_m 
2^{\eps m}\int_0^{2^{-2m}} m_k(\tau)
\tau^{\frac \eps 2}\,  \big (\De_0 \|P_m f(t)\|_{L^\infty_t L^2_x} + \int_0^t \|P_m f\|_{L^2_x}\big )\\ 
&+& \sum_m 2^{-\eps m}\int_{2^{-2m}}^1 m_k(\tau)\tau^{-\frac \eps 2}
\,\,   \big (\De_0 \|P_m f(t)\|_{L^\infty_t L^2_x} + \int_0^t \|P_m f\|_{L^2_x}\big )\\ &=&\sum_m
2^{\eps(m-k)}\int_0^{2^{2(k-m)}} \tilde m(\tau)\,\,  \big (\De_0 \|P_m f(t)\|_{L^\infty_t L^2_x} +
\int_0^t
\|P_m f\|_{L^2_x}\big )\\ &+& \sum_m 2^{\eps(k-m)}\int_{2^{2(k-m)}}^{2^{2k}} \hat m(\tau)\,\,
\big (\De_0
\|P_m f(t)\|_{L^\infty_t L^2_x} + \int_0^t \|P_m f\|_{L^2_x}\big )\\
 &\les & \sum_m 2^{-|k-m|}\big (\De_0 \|P_m f(t)\|_{L^\infty_t L^2_x} + \int_0^t \|P_m f\|_{L^2_x}\big )
\eeaa
Hence,
$$\sum_{k\ge 0}\|\int_0^t [\nab_L, P_k] \|_{L^\infty_t L^2_x}\les\sum_m \big(\De_0\c
\|P_m f(t)\|_{L^\infty_t L^2_x}+\int_0^t \|P_m f\|_{L^2_x}\big )$$
as desired.
\end{proof}
\begin{remark} In the argument above we have used  
  the following  simple bounds for $\tilde m(\tau) = \tau^{\frac \eps 2} m(\tau)$ and 
$\hat m(\tau) =\tau^{-\frac \eps 2} m(\tau)$:
\beaa
\int_0^\infty \tilde m(\tau), \quad \int_0^\infty \hat m(\tau) \les 1,\quad
\int_0^a \tilde m(\tau) \les a, \quad \int_A^\infty \hat m(\tau) \les A^{-1}
\eeaa 
which hold for all sufficiently small $a$ and all sufficiently large $A$.
Proposition \ref{prop:Sharp-Comm} now easily follows.
\end{remark}

\begin{proof}{\bf of Lemma \ref{le:main-comm}}:\quad
We prove lemma \ref{le:main-comm} by a bootstrap argument. More precisely we assume
that for any $\tau\approx 2^{-2m_0}$ and any $g$,
\be{eq:Ltau-boundM}
\|\L(\tau ) g\|_{L^\infty_t L^2_x} \les M \sum_{k\ge 0} 2^{-\eps |m_0-k|}
\|P_k g\|_{L^\infty_t L^2_x}
\end{equation}
with some positive sufficiently large constant $M$ and a fixed positive constant $\epsilon >0$,
 both independent on the function\footnote{Clearly this estimate holds true
with a constant   $M$ which  might depend on more
derivatives of the background metric.} $g$. 
Also, since the estimate above holds for 
all positive  values of $m_0$ we deduce,
\be{eq:Ltau-bound-trivialM}
\sup_{0\le \tau\le 1}\|\L(\tau ) g\|_{L^\infty_t L^2_x} \les M \|g\|_{\BB^0}
\end{equation}
We shall show that the bounds \eqref{eq:Ltau-boundM} and \eqref{eq:Ltau-bound-trivialM}
implies the stronger estimate 
\be{eq:Ltau-bound'}
\|\L(\tau ) g\|_{L^\infty_t L^2_x} \les \sum_{k\ge 0} 2^{-\eps |m_0-k|}\big(\De_0 M\c
\|P_k g\|_{L^2_t L^2_x} + 
\int_0^t \|P_k g\|_{L^2_t L^2_x}\big)
\end{equation}
from which the desired estimates \eqref{eq:Ltau-bound}, \ref{eq:Ltau-bound-trivial} follow.

We shall often use the following heat flow estimate similar to  \eqref{eq:altern-bound-U}
\be{eq:alternat-bound}
\|P_k U(\tau) f\|_{L^2_x}\les (1+ 2^{2k} \tau)^{-2}\|P_k  f\|_{L^2_x}
\end{equation}
We derive 
$$
[\nab_L, U(\tau)] f = \int_0^\tau U(\tau-\tau') [\nab_L, \Delta] U(\tau') f\, ds
$$
Recall, see proposition \ref{prop:comm-nabL-symbolic}, the  following
commutator
formula for scalars,
 \beaa
[\nab_L, \Delta ] f&=& -\trch \lap f+\Gd\c\nab^2 f+\nab\trch\c\nab f +(\Gd+\frac{1}{r})\c\Gd\c\nab
f \\
&=&-\frac{2}{r} \lap f+\Gd\c\nab^2 f+\nab\trch\c\nab f +(\Gd+\frac{1}{r})\c\Gd\c\nab f
\eeaa
Therefore,
\bea\L(\tau) f&=&I_1+I_2+I_3+I_4\label{eq:LdecompI1-I4}\\
I_1&=&\int_0^t ds \int_0^\tau d\tau U(\tau-\tau') \Big (\Gd\c\nab^2 U(\tau') f\Big )\,\nn\\
I_2&=&\int_0^t ds \int_0^\tau d\tau U(\tau-\tau') \Big (\nab \Gd\c\nab U(\tau') f\Big )\,\nn \\
I_3&=&-\frac{2}{r}\int_0^t ds \int_0^\tau d\tau U(\tau-\tau')\lap  U(\tau') f\nn\\
I_4&=&\int_0^t  ds\int_0^\tau d\tau U(\tau-\tau')\Big((\Gd+\frac{1}{r})\c\Gd\c\nab U(\tau') f\Big)\nn
\eea
{\sl  Estimate for  $I_2$.}\quad   Using the dual strong scalar
 Bernstein  inequality for $U(\tau-\tau')$ followed by H\"older and the heat flow
 estimate \eqref{eq:alternat-bound} , we infer 
\beaa
\|I_2(\tau)\|_{L^\infty_t L^2_x}
&\les &\int_0^\tau (\tau-\tau')^{-\f12}\|\nab\Gd\c\nab U(\tau') f\|_{L^1_t L^1_x}\\
&\les& \|\nab \Gd\|_{L^2_t L^2_x}\c
\int_0^\tau (\tau-\tau')^{-\f12}\|\nab U(\tau') f\|_{L^2_t L^2_x}
\\ &\les& \De_0 \sum_m \int_0^\tau (\tau-\tau')^{-\f12}
\|\nab P_m U(\tau') f\|_{L^2_t L^2_x}\\ 
&\les & \De_0\sum_m 2^m \|P_m f\|_{L^2_t L^2_x}\int_0^\tau 
(1+ 2^{2m} \tau')^{-2} (\tau-\tau')^{-\f12}d\tau'
\eeaa
Hence,
\be{eq:II-bound}
\|I_2(\tau)\|_{ L_t^\infty L^2_x}\les \De_0 \sum_m 2^{-|m-m_0|}\|P_m f\|_{L^2_t L^2_x}\quad
\mbox{for}\quad \tau\approx 2^{-2 m_0}
\end{equation}
The last statement follows from the following:
\begin{lemma}
For $\tau\approx 2^{-2 m_0}$ we have,
$$
J_m=2^m\int_0^\tau 
 (1+ 2^{2m} \tau')^{-2} (\tau-\tau')^{-\f12}d\tau'\les2^{-|m-m_0|}
$$
\end{lemma}
\begin{proof}:\quad  Indeed  for $\tau\approx 2^{-2m_0}\le  2^{-2m}$
\beaa
J_m&\les& 2^m\int_0^\tau 
  (\tau-\tau')^{-\f12}d\tau'\les 2^m\tau^{\f12}\les 2^{m-m_0}=2^{-|m-m_0|}
\eeaa
For $\tau\approx 2^{-2m_0}\ge  2^{-2m}$,
\beaa
J_m&\les& 2^m\int_0^{\tau/2}(1+ 2^{2m} \tau')^{-2} (\tau-\tau')^{-\f12}d\tau'
+2^m\int_{\tau/2}^\tau (1+ 2^{2m} \tau')^{-2} (\tau-\tau')^{-\f12}d\tau'\\
&\les&2^m\tau^{-1/2}\int_0^{\tau/2} (1+ 2^{2m} \tau')^{-2}d\tau'+2^m(1+ 2^{2m} \tau)^{-2}\tau^{1/2}
\les 2^{-m}\tau^{-1/2}\les 2^{m_0-m}\\
&=&2^{-|m-m_0|}
\eeaa
\end{proof}

{\sl  Estimate for  $I_3$.}\quad
We rewrite $I_3$ as follows\footnote{We neglect the factor $-\frac{2}{r}$
which plays no role in the estimates.},
\beaa
I_3&=&\int_0^t ds \int_0^\tau d\tau' U(\tau-\tau')\lap  U(\tau') f=
\int_0^t ds \int_0^\tau d\tau' U(\tau-\tau') U(\tau') \lap f\\
&=&\tau \int_0^t ds   \lap U(\tau) f
\eeaa
Hence, using \eqref{eq:alternat-bound} and $\tau\approx 2^{-2m_0}$,
\bea
\|I_3(\tau)\|_{L_t^\infty L_x^2}&\les&\tau\|  \lap U(\tau) f\|_{L_t^1L_x^2}\les 
\tau \sum_m\|  \lap U(\tau) P_mf\|_{L_t^1L_x^2}\nn\\&\les &
\sum_m \tau\c2^{2m} \big(1+\tau\c 2^{2m}\big)^{-2}\|P_mf\|_{L_t^1L_x^2}\nn\\
&\les&\sum _m 2^{2(m-m_0)} \big(1+ 2^{2(m-m_0)}\big)^{-2}\|P_mf\|_{L_t^1L_x^2}\nn\\
&\les&\sum_m 2^{-|m-m_0|}\|P_mf\|_{L_t^1L_x^2}\label{eq:I3-bound}.
\eea

{\sl Estimates for $I_4$.}\quad We proceed precisely
as for $I_2$, by noticing that $\|\Gd\|_{L_t^4L_x^4}\les \De_0$
in view of \eqref{eq:BA1BA2} and lemma \ref{lem:Basic-est},
\bea
\|I_4(\tau)\|_{L^\infty_t L^2_x}
&\les &\int_0^\tau (\tau-\tau')^{-\f12}\|(\Gd+\frac{1}{r})\c\Gd\c\nab U(\tau') f\|_{L^1_t L^1_x}\nn\\
&\les& \| (\Gd+\frac{1}{r})\c\Gd\|_{L^2_t L^2_x}\c
\int_0^\tau (\tau-\tau')^{-\f12}\|\nab U(\tau') f\|_{L^2_t L^2_x}
\nn\\ &\les& \De_0 \sum_m \int_0^\tau (\tau-\tau')^{-\f12}
\|\nab P_m U(\tau') f\|_{L^2_t L^2_x}\nn\\ 
&\les & \De_0\sum_m 2^m \|P_m f\|_{L^2_t L^2_x}\int_0^\tau 
(1+ 2^{2m} \tau')^{-2} (\tau-\tau')^{-\f12}d\tau'\\
&\les& \De_0\c\sum_m 2^{-|m-m_0|}\|P_mf\|_{L_t^1L_x^2}\label{eq:I4-bound}
\eea
for $\tau\approx2^{-2m_0}$.

{\sl Estimates for $I_1$}\quad This is  the most delicate term; indeed it is because of this term
that we need to make the  bootstrap  assumption \eqref{eq:Ltau-boundM}.   The difficulty
stems from the fact that we need to use the trace norm assumption  $\|\Gd\|_{L_x^\infty L_t^2}\les \De_0$
which means we have to bring the integration $\int_0^t$, along null geodesics, 
in front of $U(\tau-\tau')$ in the formula for $I_1$, see
\eqref{eq:LdecompI1-I4}. This brings in the commutator between $\int_0^t$
and $U(\tau-\tau')$ which we shall treat  according to the formula:
\be{eq:commint-U(tau)}
[\int_0^t, U(\tau-\tau')] g = \int_0^t [\nab_L, U(\tau-\tau') ]\int_0^{t'} g =
\L(\tau-\tau') \int_{0}^{t} g.
\end{equation}
Thus,
\beaa
I_1  &=& \int_0^\tau U(\tau-\tau') 
\Big (\int_0^t \Gd\c \nab^2 U(\tau') f \Big )\, d\tau' + 
\int_0^\tau \L(\tau-\tau') \Big (\int_0^t \Gd\c \nab^2 U(\tau') f \Big )\, d\tau' \\
&=& I_{11} + I_{12 }
\eeaa
Using Lemma \ref{lem:Reverse} and the dual strong scalar Bernstein inequality,
we estimate 
\beaa
\|I_{11}(\tau)\|_{L^\infty_t L^2_x}& \les &\int_0^\tau \|U(\tau-\tau')\int_0^t \div (\Gd\c \nab U(\tau')
f)\|_{L^2_x}\\& + &\int_0^\tau \|U(\tau-\tau')\int_0^t \nab\Gd\c \nab U(\tau') f\|_{L^2_x}\\  &\les &
\int_0^\tau (\tau-\tau')^{-\f12}d\tau'\big( \|\Gd\c\nab U(\tau') f\|_{L^2_x L^1_t} + 
 \|\nab \Gd\c\nab U(\tau') f\|_{L^1_x L^1_t} \big)
\\
&\les &  \Big (\|\Gd\|_{L^\infty_x L^2_t}  + \|\nab\Gd\|_{L^2_t L^2_x}\Big )
\int_0^\tau (\tau-\tau')^{-\f12}\|\nab U(\tau') f\|_{L^2_t L^2_x}\\
&\les&\De_0\int_0^\tau (\tau-\tau')^{-\f12}\|\nab U(\tau') f\|_{L^2_t L^2_x}
\eeaa
Thus, continuing  exactly as for $I_2$,
\bea
\|I_{11}(\tau)\|_{L_t^\infty L_x^2}\les  \De_0 \sum_m 2^{-|m-m_0|} \|P_m f\|_{L^2_t L^2_x},\quad
\mbox{for}\quad \tau\approx 2^{-2m_0}
\label{eq:I-tau-2m}
\eea
We now estimate $I_{12}$
with the help
of the bootstrap assumption \eqref{eq:Ltau-boundM},\eqref{eq:Ltau-bound-trivialM}
According to these we have, for any smooth function $g$, the following bounds,
\begin{align*}
&\|\L(\tau') g\|_{L^\infty_t L^2_x} \le 2M \sum_m 2^{-\eps |m-m_0|} 
\|P_m g\|_{L^\infty_t L^2_x},\quad\tau'\in [\frac \tau 2,\tau], \,\,\,\,
\tau\approx 2^{-2m_0}
\\
&\sup_{0\le\tau'\le 1}\|\L(\tau') g\|_{L^\infty_t L^2_x} \le M \|g\|_{B^0}, 
\end{align*}
Therefore,

\beaa
\|I_{12}\|_{L^\infty_t L^2_x} &\les & 
\|\int_0^\tau \L(\tau-\tau') \Big (\int_0^t 
\Gd\c \nab^2 U(\tau') f \Big )\, ds\|_{L^\infty_t  L^2_x}
\\ &\les & \int_0^{\frac \tau 2} \|\L(\tau-\tau') \Big (\int_0^t 
\Gd\c \nab^2 U(\tau') f \Big )\, ds\|_{L^\infty_t  L^2_x} \\ &+& 
\int_{\frac \tau 2}^\tau \|\L(\tau-\tau') \Big (\int_0^t 
\Gd\c \nab^2 U(\tau') f \Big )\, ds\|_{L^\infty_t L^2_x}\\ &\les &
2M\int_0^{\frac \tau 2}
\sum_m 2^{-\eps |m-m_0|}
\|P_m\int_0^t \Gd\c\nab^2 U(\tau') f  \|_{L^\infty_t L^2_x}\\
&+& M \int_{\frac \tau 2}^\tau
\|\int_0^t \Gd\c\nab^2 U(\tau') f  \|_{\BB^0}=\I_1+  \I_2
\eeaa
Using
the product  Besov estimate \eqref{eq:mainBesov-scalar}  and the estimates \eqref{eq:compete1}, \eqref{eq:compete2} we obtain
that
\beaa
\int_{\frac \tau 2}^\tau \|\int_0^t \Gd\c\nab^2 U(\tau') f  \|_{\BB^0}\,d\tau' &\les &
\De_0\int_{\frac \tau 2}^\tau \|\nab^2 U(\tau') f  \|_{\PP^0}\,d\tau'\\
&\les & \De_0 \tau\sum_{k\ge 0} \min (2^{2k}, 2^{-2k}\tau^{-2}) \|P_k f\|_{L^2_t L^2_x}\\
&\les & \De_0 \sum_{k\ge 0} 2^{-2|k-m_0|} \|P_k f\|_{L^2_t L^2_x}
\eeaa
We conclude that
\be{eq:bound-I12}
\I_2\les \De_0 M\sum_{k\ge 0} 2^{-2|k-m_0|} \|P_k f\|_{L^2_t L^2_x}
\end{equation}
To estimate $\I_1$ we observe that according to the estimate 
\eqref{eq:mainlemma2-dyadic} of
the remark \ref{rem:mainBesov-dyadic} 
\beaa
\|P_m\int_0^t \Gd\c\nab^2 U(\tau') f  \|_{L^\infty_t L^2_x}&\les& \Big (\NN_1[\Gd]
+ \|\Gd\|_{L^\infty_x L^2_t}\Big ) \c\sum_{k\ge 0} 2^{-\si|m-k|}\|P_k\nab^2 U(\tau') f\|_{L^2_t
L^2_x}\\ &\les & 
\De_0 \sum_{k\ge 0} 2^{-\si|m-k|}\|P_k\nab^2 U(\tau') f\|_{L^2_t L^2_x}
\eeaa
for some positive constant $\si>\eps$.
Therefore, with the help
of  lemma \ref{le:last-homogBesov} below,
\bea
\I_{1} &\les &
2\De_0 M \sum_m 2^{-\eps |m-m_0|}\sum_{k\ge 0} 2^{-\si|m-k|}
\int_0^{\frac \tau 2}\|P_k \nab^2 U(\tau') f  \|_{L^2_t L^2_x}\nn\\ &\les &
\De_0 M \c\sum_m 2^{-\eps |m-m_0|}\int_0^{\frac \tau 2}
\|P_m \nab^2 U(\tau') f  \|_{L^2_t L^2_x}=\De_0 M\c \J \nn\\
&\les&  \De_0 M\c \sum_{k\ge 0}  2^{-\eps |k-m_0|} 
\|P_k f\|_{L^2_t L^2_x}\label{eq:boundI111}
\eea
 Thus,
\beaa
\|I_{12}(\tau)\|_{L_t^\infty L_x^2}&\les &M \De_0 \sum_m 2^{-|m-m_0|} \|P_m f\|_{L^2_t L^2_x},\quad
\mbox{for}\quad \tau\approx 2^{-2m_0},
\label{eq:I-tau-2m-2}
\eeaa
which ends the proof of the improved estimate
\eqref{eq:Ltau-bound'}  and therefore also
of lemma \ref{le:main-comm}.
\end{proof}
It only remains to prove the following,
\begin{lemma} For $\tau\approx 2^{-2m_0}$ the integral $\J=\J^{+}+\J^{-}$
defined by 
\beaa
\J^{+}&=&\sum_{m\ge m_0} 2^{-\eps |m-m_0|}\int_0^{\frac \tau 2}
\|P_m \nab^2 U(\tau') f  \|_{L^2_t L^2_x}\\
\J^{-}&=&\sum_{0\le m<m_0} 2^{-\eps |m-m_0|}\int_0^{\frac \tau 2}
\|P_m \nab^2 U(\tau') f  \|_{L^2_t L^2_x}
\eeaa
 verifies the estimate,
\beaa
\J\les   \sum_{k\ge 0}  2^{-\eps |k-m_0|} 
\|P_k f\|_{L^2_t L^2_x}
\eeaa

\label{le:last-homogBesov}
\end{lemma}
\begin{proof}:\quad

{\bf 1)}\,\,\,{\sl  Estimates for $\J_{+}$.}\quad
We are in the case  $m\ge m_0$ ( or $2^{-2m}\le\tau$)
We  decompose further
\beaa
\|P_m \nab^2 U(\tau') f  \|_{L^2_t L^2_x} &\les &\sum_{k\ge 0} \|P_m \nab^2 P_k U(\tau') f  \|_{L^2_t L^2_x}\\ &\les & \sum_{k\le m_0}
 \|P_m \nab^2 P_k U(\tau') f  \|_{L^2_t L^2_x} + \sum_{k \ge m}
 \|P_m \nab^2 P_k U(\tau') f  \|_{L^2_t L^2_x} \\ &+& 
 \sum_{m_0<k<m} \|P_m \nab^2 P_k U(\tau') f  \|_{L^2_t L^2_x}
 =J_1 + J_2 + J_3 
\eeaa
Let $\Jp_1,\Jp_2,\Jp_3$ be the corresponding contributions to $\J$.
The term $\J_3$   is easiest to  estimate,
$$
J_3=\sum_{m_0<k<m} \|P_m \nab^2 P_k U(\tau') f  \|_{L^2_t L^2_x}\les 
\sum_{m_0<k<m} \|\nab^2 P_k U(\tau') f  \|_{L^2_t L^2_x}
$$
Using lemma \ref{lem:Besov-heat} we obtain 
$$
\int_0^{\frac \tau 2} J_3\les 
\sum_{m_0<k<m} \int_0^{\frac \tau 2}\|\Delta U(\tau') P_k  f  \|_{L^2_t L^2_x}\les 
\sum_{m_0<k<m} \|P_k f\|_{L^2_t L^2_x}
$$
Therefore,
$$
\Jp_3= \sum_{m>m_0} 2^{-\eps |m-m_0|} \sum_{m_0<k<m} 
\|P_k f\|_{L_t^2L_x^2}\les  \sum_{k\ge 0}  2^{-\eps |k-m_0|} 
\|P_k f\|_{L^2_t L^2_x}
$$
To estimate $\Jp_1$ we proceed as follows. Applying the dual finite band property 
followed by the finite band property  of the $P_k's$, we obtain
\beaa
J_1 \les 2^m \sum_{k\ge m}\|\nab P_k U(\tau') f\|_{L^2_t L^2_x} &\les & 
\sum_{m\ge k} 2^{m+k} \|P_k U(\tau') f\|_{L_t^2L_x^2} \\&\les & 
\sum_{k\ge m}2^{m-k} \| \Delta U(\tau') P_k f\|_{L^2_t L^2_x}
\eeaa
Once again, lemma \ref{lem:Besov-heat}  gives 
$$
\int_0^{\frac \tau 2} J_1 \les \sum_{k\ge m}2^{m-k} \|P_k f\|_{L^2_t L^2_x}
$$
Thus,
\beaa
\Jp= \sum_{m>m_0} 2^{-\eps |m-m_0|} \sum_{k\ge m} 2^{-|m-k|}
\|P_k f\|_{L^2_t L^2_x}\les  \De_0 M \sum_{k\ge 0}  2^{-\eps |k-m_0|} 
\|P_k f\|_{L^2_t L^2_x}.
\eeaa
Finally, we have for $J_2$.
\beaa
J_2 \les \sum_{k<m_0} \|\nab^2 P_k U(\tau') f\|_{L_t^2L_x^2}\les 
 \sum_{k<m_0} 2^{2k} \|P_k f\|_{L_t^2L_x^2}
\eeaa
Therefore,
\beaa
\int_0^{\frac \tau 2} J_2 \les \sum_{k<m_0} 2^{2k} \tau 
\|P_k f\|_{L^2_t L^2_x} \les  \sum_{k<m_0} 2^{-2|k-m_0|} 
\|P_k f\|_{L^2_t L^2_x}.
\eeaa
and we infer that
$$
\Jp_2=\sum_{m>m_0} 2^{-\eps |m-m_0|} \sum_{k< m_0} 2^{-2|m_0-k|}
\|P_k f\|_{L^2_t L^2_x}\les \sum_{k\ge 0}  2^{-\eps |k-m_0|} 
\|P_k f\|_{L^2_t L^2_x}.
$$

{\bf 2)}\quad {\sl Estimate for $\Jm$.}\quad
We are in the case  $m\le m_0$ ($2^{-2m}\ge \tau$)
As before we decompose further as follows:
\beaa
\|P_m \nab^2 U(\tau') f  \|_{L^2_t L^2_x} &\les &\sum_{k\ge 0} \|P_m \nab^2 P_k U(\tau') f  \|_{L^2_t L^2_x}\\ &\les & \sum_{k\le m}
 \|P_m \nab^2 P_k U(\tau') f  \|_{L^2_t L^2_x} + \sum_{k \ge m_0}
 \|P_m \nab^2 P_k U(\tau') f  \|_{L^2_t L^2_x} \\ &+& 
 \sum_{m<k<m_0} \|P_m \nab^2 P_k U(\tau') f  \|_{L^2_t L^2_x}
 =N_1 + N_2 + N_3 
\eeaa
and denote by $\Jm_1, \Jm_2, \Jm_3$ the corresponding contributions to $\Jm$.
The estimates for $N_3$ and $\Jm_3$  are identical to that for $J_3, \Jp_3$.

To estimate  $N_1$ we proceed as follows.
\beaa
N_1 \les \sum_{k\le m}\|\nab^2 P_k U(\tau') f\|_{L^2_t L^2_x} \les  
\sum_{k\le m} 2^{2k} \|P_k f\|_{L^2_t L^2_x} 
\eeaa
Integrating in $s$ and using that $\tau\approx 2^{-2m_0}$ we obtain the bound
\beaa
\Jm_1= \sum_{m\le m_0} 2^{-\eps |m-m_0|} \sum_{k\le m} 2^{-2|m_0-k|}
\|P_k f\|_{L^2_t L^2_x}\les \sum_{k\ge 0}  2^{-2|k-m_0|} 
\|P_k f\|_{L^2_t L^2_x}
\eeaa
Finally, we have for $N_2$.
\beaa
N_2 \les \sum_{k\ge m_0} \|P_m\nab^2 P_k U(\tau') f\|_{L^2_t L^2_x}\les 
 \sum_{k\ge m_0} 2^{m-k} \|\Delta U(\tau') P_k f\|_{L^2_t L^2_x}
\eeaa
Integrating in $s$, using lemma \ref{lem:Besov-heat}, we obtain the estimate
\beaa
\Jm_2 =\sum_{m\le m_0} 2^{-\eps |m-m_0|} \sum_{k\ge m_0} 2^{-|m-k|}
\|P_k f\|_{L^2_t L^2_x}&\les& \sum_{k\ge 0}  2^{-\eps |k-m_0|} 
\|P_k f\|_{L^2_t L^2_x}.
\eeaa
Combining the estimates involving $\Jp_1, \Jp_2, \Jp_3$ and
 $\Jm_1, \Jm_2, \J_3$  we 
conclude that
$$
\J\les  \sum_{k\ge 0}  2^{-\eps |k-m_0|} 
\|P_k f\|_{L^2_t L^2_x}.
$$
as desired.
\end{proof}
\section{Commutator estimates}\label{sec:Commutator}
This section  we prove the  commutator lemma  \ref{lem:Comm-nabL}  
which we recall below.
\begin{proposition}\label{prop:Commut-nabL}
For any smooth $S$-tangent tensor field $F$ and an arbitrary \mbox {$1\le q<2$}
the following estimate holds true 
\bea
\|[P_k,\nab_L] f\|_{L^q_t L^2_x} + 2^{-k} \|\nab [P_k,\nab_L] f\|_{L^q_t L^2_x} 
\les  2^{-\frac k2+}\,\N_1[F]\label{eq:comm-estim-q-nonsh}
\eea
In addition,
\bea
\|[P_k,\nab_L] f\|_{L^1_t L^2_x} + 2^{-k} \|\nab [P_k,\nab_L] f\|_{L^1_t L^2_x} 
\les  2^{-k+}\,\N_1[F] \label{eq:comm-estim-1-nonsh}
\eea
\end{proposition}
\begin{proof}:\quad
Using the definition of the LP projection $P_k$ we obtain
\bea
[P_k, \nab_L]  F& =& \int_0^\infty m_k (\tau) [U(\tau) , \nab_L]  F\, d\tau\\
\Phi (\tau)&=& [U(\tau), \nab_L] F
\eea
To calculate $[U(\tau), \nab_L] F$ we recall the commutator formula 
for $[\ddd_L,\lap]$ from proposition \ref{prop:comm-nabL-symbolic}
\beaa
[\ddd_L,\lap] F&=&-\trch\lap F+\Gd\c\nab^2 F+\nab \Gd\c\nab F+(\Gd+\frac{1}{r})\c\Gd\c\nab F\\
&+&\b \c \nab F +\nab\big(\b\c F+(\Gd+\frac{1}{r} )\c \Gd\c F\big)
\eeaa
which we rewrite in the simplified, symbolic,  form,
\beaa
[\ddd_L,\lap] F&=&\nab\bigg((\frac{1}{r}+\Gd)\c \nab F +\big(\b+ (\Gd+\frac{1}{r})\c \Gd\big)\c F \bigg)\\
&+& \bigg(\nab\Gd+\big(\b+ (\Gd+\frac{1}{r})\c \Gd\big) \bigg)  \c \nab F
\eeaa
Observe that the terms $ (\Gd+\frac{1}{r})\c \Gd$ and $\nab A$  have same 
or   better estimates than $\b$. Therefore we
can discard them   and simplify,
\bea
[\ddd_L,\lap] F&=&\nab\bigg((\frac{1}{r}+\Gd)\c \nab F +\b\c F \bigg)+\b\c \nab F
\eea
\def\Gdh{\hat{\Gd}}
Consequently, setting $\Gdh=\Gd+\frac{1}{r}$,

\beaa
\Phi(\tau )&=&\int_0^\tau U(\tau-\tau') [\Delta, \nab_L] U(\tau') F\\
&=& \int_0^\tau U(\tau-\tau') \nab \Big (\Gdh\c\nab  U(\tau') F + 
\b  U(\tau') F\Big )\\ &+&\int_0^\tau U(\tau-\tau') \b\c \nab U(\tau')  F\\ &=& \Phi_1(\tau) + \Phi_2(\tau)
\eeaa
In what follows we shall rely on the  following heat flow estimates:
\begin{lemma}
\label{le:funny-heatflow-estim}
For  any $2\le p<\infty$,
\bea
\|U(\tau) G\|_{L^\infty_x} &\les &\big (1 +
\|K\|_{L^2_x}^{\frac 1{p-1}}\big )  \|U(\tau'/2)G\|_{H^1_x},\label{eq:heat-inf}\\
\|U(\tau) G\|_{L^\infty_x} &\les &\langle \tau^{-\frac 1r}\rangle \big (1 +
\|K\|_{L^2_x}^{\frac 2{r(p-1)}}\big )  \|G\|_{L^r_x},\qquad 2\le \forall r\le \infty\label{eq:heat-L4}\\
\|\nab^2 U(\tau) G\|_{L^2_x} &\les & (\tau^{-\f12}  
+ \|K\|_{L^2_x}^{\frac p{p-1}}\big ) \|U(\tau'/2)G\|_{H^1_x}.\label{eq:heat-nab2}
\eea
\end{lemma}
\begin{proof}:\quad
The estimates  follow easily, by interpolation,  from the  heat flow estimates 
of  proposition \ref{prop:heat-flow-estimates}combined
with  the Bochner
and
$L^\infty$ inequalities  \eqref{eq:conseq-Bochner-ineq}, \eqref{eq:LinftyL2tensor}.
\end{proof}

Our proof will now proceed as follows:

{\bf Step 1}: \quad We shall first establish the following  estimates for $\Phi_1$:
\beaa
\|\Phi_1(\tau)\|_{L^q_t L^2_x} &\les& \max (\tau^{\frac 14},\, \tau^{\frac 12})
\N_1[F]\\
 \|\Phi_1(\tau)\|_{L^1_t L^2_x}& \les& \max (\tau^{\f12-},\, \tau^\f12) \N_1[F]
\eeaa
for any $1\le q<2$ close to $2$ and $\tau>0$.

{\bf Step 2}:\quad We show that $\Phi_2$ verifies the  weak estimate,
$$ \|\Phi_2(\tau)\|_{L^q_t L^2_x} \les \max (\tau^{\frac 1q-\f12-},\, \tau^{\f12 +})\,
 \N_1[F]
$$
for any $1\le q<2$ close to $2$.

{\bf Step 3}:\quad  Using the results above we conclude that the estimate
\eqref{eq:comm-estim-1-nonsh} holds true.

{\bf Step 4}:\quad 
With the help of \eqref{eq:comm-estim-1-nonsh} we can improve the estimate
for $\Phi_2$ and derive the desired estimate \eqref{eq:comm-estim-q-nonsh}.

\subsection{$L^q_t L^2_x$ estimates for $\Phi_1$}\quad
\begin{lemma}\label{lem:Phi1-straight}
For all $q<2$ sufficiently close to $q=2$ and any $\tau\ge 0$,
$$
\|\Phi_1(\tau)\|_{L^q_t L^2_x} \les \max (\tau^{\frac 14},\, \tau^{\frac 12})
\N_1[F]
$$
\end{lemma}
\begin{proof}:\quad 
We start by estimating the expression  
$$
\phi_1(\tau') =\Gdh\c \nab  U(\tau') F + \b \c U(\tau') F
$$
Applying H\"older inequality followed by Gagliardo-Nirenberg \eqref{eq:GNirenberg} and
the estimates \eqref{eq:heat-inf}-\eqref{eq:heat-nab2}, we obtain
\beaa
\|\phi_1(\tau')\|_{ L^2_x}&\les & \|\Gdh \c \nab U(\tau') F\|_{ L^2_x} +
 \|\b \c U(\tau') F\|_{L^2_x}
\\ &\les& \|\Gdh\|_{L^4_x}\c \|\nab U(\tau') F\|_{ L^4_x} + 
\|\b\|_{L^2_x} \c \|U(\tau') F\|_{ L^\infty_x} \\ &\les & 
\|\Gdh\|_{L^4_x}\c \|\nab^2 U(\tau') F\|_{ L^2_x}^{\f12} \c
\|\nab U(\tau') F\|_{ L^2_x}^{\f12} + 
\|\b\|_{L^2_x}\c 
\|U(\tau') F\|_{ L^\infty_x} \\ &\les &
\|\Gdh\|_{L^4_x}\big ({\tau'}^{-\frac 14} + \|K\|_{L^2_x}^{\frac p{2(p-1)}}\big )
\|U(\tau'/2)F\|_{H^1_x} \\ &+& 
\langle {\tau'}^{-\frac 14}\rangle\|\b\|_{L^2_x} 
\big (1+\|K\|_{L^2_x}^{\frac 1{2(p-1)}}\big ) \|F\|_{L^4_x}
\eeaa
\begin{remark}\label{rem:take-p-infty}
To simplify the  exposition we will take $p=\infty$ 
in the above estimate. Properly speaking this is not
acceptable since the estimates of lemma \ref{le:funny-heatflow-estim}
do not hold true for $p=\infty$. The correct proof requires the    choice
of an appropriate large  value of $p<\infty$ dependent on the exponent
$q<2$. As long as we stay away from the critical exponent $q=2$ one
can easily correct  the slightly idealized  exposition below. We note in fact
that  the restriction to $q<2$  is due to the presence of the Gauss curvature terms $\|K\|_{L^2_x}$,
which were generated by the B\"ochner identity and inequality for tensors. 
Many of the technical complications below are due to the presence of these terms.

\end{remark}
Hence,
\beaa
\|\phi_1(\tau')\|_{L^2_x}&\les & 
\|\Gdh\|_{L^4_x}\big ({\tau'}^{-\frac 14} + \|K\|_{L^2_x}^{\frac 1{2}}\big )
\|U(\tau'/2)F\|_{H^1_x} +
\langle {\tau'}^{-\frac 14}\rangle\|\b\|_{L^2_x}  \|F\|_{L^4_x}
 \eeaa
 Fix an exponent $q,\, 1\le q <2$ and use the interpolated heat flow  estimate, see proposition
\ref{prop:heat-flow-estimates}:
 $$
 \|U(\tau'/2)F\|_{H^1_x}\les\langle{\tau'}^{-(\frac 34-\frac 1q)} \rangle
  \| F\|_{H^1_x}^{\frac 2q-\f12} \|F\|_{L^2_x}^{\frac 32 - \frac 2q}
 $$
 Therefore,
\beaa
\|\phi_1(\tau')\|_{L^2_x}&\les &\|\Gdh\|_{L^4_x}\Big ( {\tau'}^{-\frac 14}\| F\|_{H^1_x}  + 
 \langle{\tau'}^{-(\frac 34-\frac 1q)}\rangle\|K\|_{L^2_x} ^{\frac 12} \| F\|_{H^1_x}^{\frac 2q-\f12} \|F\|_{L^2_x}^{\frac 32 - \frac 2q}\Big ) 
\\ & +& \langle{\tau'}^{-\frac 14}\rangle
 \|\b\|_{L^2_x}  \|F\|_{L^4_x} 
\eeaa
We now take  the  $L^q_t$-norm, use  the  assumptions   {\bf A1}, {\bf A2}, {\bf K1}
 $$
\|\Gdh\|_{L^4_x L^\infty_t},  \|\b\|_{L^2_t L^2_x},  \|K\|_{L^2_t L^2_x}
\les 1,
$$
apply  the Gagliardo-Nirenberg estimate \eqref{eq:GNirenberg}, and the inequality 
$\|F\|_{L^{\frac {2q}{2-q}}_t L^4_x} \les \N_1[F]$ to 
 obtain 
\beaa
\|\phi_1(\tau')\|_{L^q_t L^2_x} &\les & 
\Big ({\tau'}^{-\frac 14} \|F\|_{L^q_t H^1_x} + 
 \langle{\tau'}^{-(\frac 34-\frac 1q)}\rangle
\|K\|_{L^2_t L^2_x}^\f12 \|F\|_{L^2_t H^1_x}^{\frac 2q-\f12} 
\|F\|_{L^\infty_t L^2_x}^{\frac 32 - \frac 2q}\Big ) \\
& + &\langle{\tau'}^{-\frac 14}\rangle
  \|F\|_{L^{\frac {2q}{2-q}}_t L^4_x}
 \les \langle{\tau'}^{-\frac 14}\rangle \N_1[F]
\eeaa

 Returning to the estimates for $\Phi_1(\tau)$ we derive
 \bea
 \|\Phi_1(\tau) \|_{L^q_t L^2_x} &\les & \|\int_0^\tau U(\tau-\tau') 
 \nab \phi_1 (\tau') \|_{L^q_t L^2_x} \nn\\ &\les &  
  \int_0^\tau (\tau-\tau')^{-\frac 12} \|\phi_1 (\tau')\|_{L^q_t L^2_x}\nn\\
  &\les & \int_0^\tau (\tau-\tau')^{-\frac 12} 
  \langle{\tau'}^{-\frac 14}\rangle\, \N_1[F]\nn\\
  &\les & \max (\tau^{\frac 14}, \, \tau^\f12 )\, \N_1[F]\label{eq:Phi1-tau}
 \eea
 \end{proof}
 \subsection{$L^1_t L^2_x$ estimate for $\Phi_1$}\quad
 
 In this section we derive an improved estimate for $\Phi_1(\tau)$ in 
 the $L^1_t L^2_x$-norm.
 \begin{lemma}\label{lem:Phi1-twist}
 For any $\tau\ge 0$,
 $$
 \|\Phi_1(\tau)\|_{L^1_t L^2_x} \les \max (\tau^{\f12-},\, \tau^\f12) \N_1[F]
 $$
 \end{lemma}
 \begin{proof}:\quad 
 Fix a sufficiently large exponent $r<\infty$. Then applying H\"older inequality
 followed by Gagliardo-Nirenberg estimate and heat flow
estimates \eqref{eq:heat-L4}, \eqref{eq:heat-nab2}, we obtain
 \beaa
 \|\phi_1(\tau')\|_{ L^2_x}&\les & \|\Gdh\c \nab U(\tau') F\|_{ L^2_x} +
 \|\b \c U(\tau') F\|_{L^2_x}
\\ &\les& \|\Gdh\|_{L^r_x} \|\nab U(\tau') F\|_{ L^{\frac {2r}{r-2}}_x} + 
\|\b\|_{L^2_x}  \|U(\tau') F\|_{ L^\infty_x} \\ &\les & 
\|\Gdh\|_{L^r_x} \|\nab^2 U(\tau') F\|_{ L^2_x}^{\frac 2r} 
\|\nab U(\tau') F\|_{ L^2_x}^{1-\frac 2r} + 
\|\b\|_{L^2_x} \c
\|U(\tau') F\|_{ L^\infty_x} \\ &\les &
\|\Gdh\|_{L^r_x}\big ({\tau'}^{-\frac 1r} + \|K\|_{L^2_x}^{\frac {2p}{r(p-1)}}\big )
\|U(\tau'/2)F\|_{H^1_x} \\ &+& 
\langle {\tau'}^{-\frac 1r}\rangle\|\b\|_{L^2_x} 
\big (1+\|K\|_{L^2_x}^{\frac 2{r(p-1)}}\big ) \|F\|_{L^r_x}
 \eeaa
Once again we set $p=\infty$, see remark \ref{rem:take-p-infty}, 
\be{eq:estim-phi1-x}
 \|\phi_1(\tau')\|_{ L^2_x}\les  \big({\tau'}^{-\frac 1r}+\|K\|_{L_x^2}^{\frac{2}{r}}\big) \|\Gdh\|_{L^r_x}
\|F\|_{H^1_x} + 
\langle {\tau'}^{-\frac 1r}\rangle\|\b\|_{L^2_x} 
 \|F\|_{L^r_x}
\end{equation}
Taking $L^1_t$-norm and using our assumptions $  {\bf A1}, {\bf A2}, {\bf K1}$
we derive
\bea
\|\phi_1(\tau')\|_{ L^1_t L^2_x}&\les & {\tau'}^{-\frac 1r} \|\Gdh\|_{L^2_t L^r_x}
\|F\|_{L^2_t H^1_x} + 
\langle {\tau'}^{-\frac 1r}\rangle\|\b\|_{L^2_t L^2_x} 
 \|F\|_{L^2_t L^r_x}\nn\\
&+&\|\Gdh\|_{L_t^{\frac{2r}{r-2}} L_x^r}\c \|K\|_{L_t^2L_x^2}^{\frac{2}{r}} \|F\|_{L_t^2 H^1_x}
  \les  
\langle {\tau'}^{-\frac 1r}\rangle \N_1[F].\label{eq:finally1}
\eea
for an arbitrarily large $r$.

To derive the last estimate we have used lemma \ref{lem:no-end} below.

 Finally  we recall that $\Phi_1(\tau) = \int_0^\tau U(\tau-\tau') 
\nab \phi_1(\tau')$  and proceeding as before, with $r$ arbitrarily large,  we finish the proof.
 \end{proof}
\begin{lemma}\label{lem:no-end}
For any $2\le p<\infty$  and $2\le q< \frac{2p}{p-4} $ 
$$
\|F\|_{L^q_t L^p_x} \les \N_1[F]
$$
\end{lemma}
\begin{proof}:\quad The proof follows immediately
by interpolating between  the following:  
$$
\|F\|_{L^2_t L^r_x}\les \N_1[F],\qquad 
\|F\|_{ L^\infty_t L^4_x} \les \N_1[F].
$$
\end{proof}
\subsection{$L^q_t L^2_x$ estimates for $\Phi_2 $}\quad
 \begin{lemma}\label{lem:Phi2-straight}
 For any $1\le q<2$ and all $\tau \ge 0$,
 $$
 \|\Phi_2(\tau)\|_{L^q_t L^2_x} \les \max (\tau^{\frac 1q-\f12-},\, \tau^{\f12 +})\,
 \N_1[F]
 $$
 \end{lemma}
 \begin{proof}:\quad
We start by estimating the function  
$$
\phi_2 (\tau')= \nab \Gd\c \nab U(\tau') F + \b  \c \nab U(\tau') F
$$
Fix $r>1$ sufficiently close to $r=1$ and apply H\"older inequality followed by the 
 Gagliardo-Nirenberg estimate, 
\beaa
\|\phi_2(\tau')\|_{L^r_x}&\les & \Big (\|\nab\Gd\|_{L^2_x} +
\|\b\|_{L^2_x} \Big ) 
\|\nab U(\tau') F\|_{L^{\frac {2r}{2-r}}_x} \\&\les & 
\Big (\|\nab\Gd\|_{L^2_x} +\|\b\|_{L^2_x} \Big ) 
\|\nab ^2 U(\tau') F\|_{L^2_x}^{2-\frac 2r}
\|\nab U(\tau') F\|_{L^2_x}^{\frac 2r-1}.
\eeaa
We now use estimate \eqref{eq:heat-nab2} to obtain
\beaa
\|\phi_2(\tau')\|_{ L^r_x}\les  
 \Big (\|\nab\Gd\|_{L^2_x} +\|\b\|_{L^2_x}\Big )
\Big (\tau^{-\frac {r-1}r}+ 
 \|K\|_{L^2_x} ^{2\frac {p(r-1)}{r(p-1)}}\Big )
\| U(\tau'/2) F\|_{H^1_x}
\eeaa
Proceeding in the spirit of remark \ref{rem:take-p-infty} we  set $r=1$ 
on the right hand side
 and  replace the above inequality with, 
\be{eq:phi2-inter}
\|\phi_2(\tau')\|_{L^r_x} \les \Big (\|\nab\Gd\|_{L^2_x} +\|\b\|_{L^2_x}
 \Big )\|\nab U(\tau'/2)F\|_{L^2_x}
\end{equation}
For a given exponent $q$ we use the standard heat flow estimates
 $\|\nab U(\tau'/2) F\|_{L^2_x} \les {\tau'}^{-\f12} \|F\|_{L^2_x}$  and  
 $\|\nab U(\tau'/2) F\|_{L^2_x} \les \|\nab F\|_{L^2_x}$ to  write
$$
\|\nab U(\tau'/2) F\|_{L^2_x} \les {\tau'}^{-(1-\frac 1q)}
\|\nab F\|_{L^2_x}^{\frac 2q -1} 
\|F\|_{L^2_x}^{2-\frac 2q}
$$
Now taking the $L^q_t $-norm and using, 
$
\,\|\nab\Gd\|_{L^2_t L^2_x} ,\, 
 \|\b\|_{L^2_t L^2_x}\les \De_0\les 1, \qquad 
$
we derive with the help of Gagliardo-Nirenberg estimate that
 \beaa
 \|\phi_2(\tau') \|_{L^q_t L^r_x} &\les& \, {\tau'}^{-(1-\frac 1q)}
 \|\nab F\|_{L^2_t L^2_x}^{\frac 2q -1} 
\|F\|_{L^\infty_t L^2_x}^{2-\frac 2q}+ \|\nab F\|_{L^q_t L^2_x}\\
 &\les & \, \langle {\tau'}^{-(1-\frac 1q)}\rangle\, \N_1[F]
 \eeaa
Returning to the estimates for $\Phi_2(\tau)$ and applying the  dual weak Bernstein 
 inequality $L^2\to L^r$, see \eqref{eq:dual-heat-GN},  for some $r>1$ sufficiently close to $r=1$,  we
obtain
 \bea
 \|\Phi_2(\tau')\|_{L^q_t L^2_x}&\les& \|\int_0^\tau U(\tau-\tau') 
 \phi_2(\tau')\|_{L^q_t L^2_x}\nn\\ &\les &
  \int_0^\tau (\tau-\tau')^{-(\frac 1r-\f12 )}\| \phi_2(\tau')\|_{L^q_t L^r_x}\nn\\
  &\les & \int_0^\tau (\tau-\tau')^{-(\frac 1r-\f12 )}
  \langle {\tau'}\rangle^{-(1-\frac 1q)}\, \N_1[F]\nn\\
  &\les & \max (\tau^{\frac 12+\frac 1{q}-\frac 1r},\, \tau^{\frac 32-\frac 1r})\,
  \N_1[F]\nn \\&\les & 
   \max (\tau^{\frac 1{q}-\frac 12},\, \tau^{\frac 12+})\,
  \N_1[F]
  \label{eq:Phi2-tau}
 \eea
 \end{proof}
We can now finish the  proof of proposition \ref{prop:Commut-nabL} for the case 
 of $q=1$. 
Recall  that
$$
[P_k,\nab_L] F = \int_0^\infty m_k (\tau) \Phi(\tau)\, d\tau,
$$
where $\Phi(\tau)=\Phi_1(\tau) +\Phi_2(\tau)$. 
 Combining the results of lemmas \ref{lem:Phi1-twist} and \ref{lem:Phi2-straight} 
 we derive
 \bea
 \|[P_k, \nab_L] F\|_{L^1_t L^2_x} &\les &  \N_1[F]\c\int_0^\infty m_k (\tau) 
  \max (\tau^{\frac 12-}, \, \tau^{\f12+} ) \,\nn\\ &\les & 
 \N_1[F] \c\big(\int_0^1 m_k (\tau) \,\tau^{\frac 12-}\, + 
 \int_1^\infty m_k(\tau) \c\,\tau^{\f12 +}\big)\nn\\ &\les & \N_1[F]\c
\big( 2^{-k+}\int_0^{2^{2k}} m(\tau) \,\tau^{\frac 12-}\,  + 
 2^{-k-}\int_{2^{2k}}^\infty m(\tau) \,\tau^{\f12 +}\big)\nn
\\&\les &
 2^{-k+}\,\c \N_1[F]\label{eq:PknabL-comm}
 \eea
 To obtain the estimate for the term $\|\nab [P_k, \nab_L] F \|_{L_t^1L_x^2}$ we need to control
 $\nab\Phi_1(\tau)$ and $\nab\Phi_2(\tau)$.
 Observe that the tensor-fields $\Phi_j(\tau),\, j=1,2$ are the solutions 
 of the respective heat equations
 \begin{align*}
& \pr_\tau \Phi_1 - \Delta \Phi_1 = \nab \phi_1,\qquad \Phi_1|_{\tau=0}=0,\\
&\pr_\tau \Phi_2 - \Delta \Phi_2 = \phi_2,\qquad \Phi_2|_{\tau=0}=0
 \end{align*}
 The standard heat flow estimates imply that
 \begin{align}
 &\|\nab \Phi_1\|_{L^2_\tau L^2_x} \les \|\phi_1\|_{L^2_\tau L^2_x},\label{eq:st-heat-flow-1}\\
&\|\nab \Phi_2\|_{L^2_\tau L^2_x} \les \|\phi_2\|_{L^{\frac {2p}{3p-2}}_\tau L^p_x}, \qquad 1\le p\le 2
\label{eq:st-heat-flow-2}
 \end{align}
where $\|\, \|_{L^2_\tau }$ refer to  the $L^2$ norm on the interval $[0,\tau]$.
The proof of \ref{eq:st-heat-flow-1} follows by multiplying the  heat equation
for $\Phi_1$ by $\Phi_1$ and integrating by parts. Estimate \eqref{eq:st-heat-flow-2}
follows in the same  manner, replacing the integration by parts
by an appropriate use of H\"older and interpolation.

Taking the $L_\tau^2$ norm of the equations \eqref{eq:estim-phi1-x}
and \eqref{eq:phi2-inter} followed by the $L_t^1$ norm we derive, for any
$p>1$ sufficiently close to $1$ and all sufficiently large $r$,
\beaa
\|\phi_1\|_{ L_t^1 L^2_\tau L^2_x}&\les&\max(\tau^{\f12}, \tau^{\f12-\frac{1}{r}})\c\N_1[F]\\
\|\phi_2\|_{ L_t^1 L^{\frac {2p}{3p-2}}_\tau L^p_x}&\les& \tau^{\frac{3p-2}{2p}}\c\N_1[F]
\eeaa
Therefore,
\beaa
\|\nab\Phi_1\|_{ L_t^1 L^2_\tau L^2_x}&\les&\max(\tau^{\f12}, \tau^{\f12-})\c\N_1[F]\\
\|\nab\Phi_2\|_{ L_t^1 L^2_\tau L^2_x}&\les& \max(\tau^{\f12}, \tau^{\f12-})\c\N_1[F]
\eeaa
One  can, by a standard argument, convert  the above into the following
weighted estimates,
\beaa
\|\min(\tau^{-\f12-}, \tau^{-\f12+})\nab\Phi_1\|_{ L_t^1 L^2_\tau L^2_x}&\les&\N_1[F]\\
\|\min(\tau^{-\f12-}, \tau^{-\f12+})\nab\Phi_2\|_{ L_t^1 L^2_\tau L^2_x}&\les& \N_1[F]
\eeaa
 Thus
 \beaa
 \|\nab [P_k, \nab_L] F\|_{L^1_t L^2_x} &\les & 
 \|\int_0^\infty m_k (\tau ) \nab \Phi(\tau)\|_{L^1_t L^2_x}\\
 &\les &  \| \max (\tau^{\f12-},\, \tau^{\f12 +})\, m_k (\tau )\|_{L^2_\tau} 
  \|\min (\tau^{-\f12+},\, \tau^{-\f12 -})\,\nab \Phi_j\|_{L^1_t L^2_\tau L^2_x}\\
  &\les & 2^{k} 2^{-k+}\,\N_1[F]
 \eeaa
To finish the proof of proposition \ref{prop:Commut-nabL} for 
an arbitrary $q<2$ we first need to obtain improved estimates 
for $\Phi_2$.

\subsection{Improved estimates for $\Phi_2$}\quad

\begin{lemma}\label{lem:Phi2-twist}
For all $1\le q<2$ and any $\tau\ge 0$
$$
\|\Phi_2(\tau)\|_{L^q_t L^2_x} \les 
\max (\tau^{\frac 14-},\, \tau^{\frac 12})\, N_1[F]
$$
\end{lemma}
\begin{proof}:\quad 
We start by proceeding in the same way as in the estimates  
for $\Phi_2(\tau)$ of lemma \ref{lem:Phi2-straight} 
until equation \eqref{eq:phi2-inter},
$$
\|\phi_2(\tau')\|_{L^r_x} \les \Big (\|\nab\Gd\|_{L^2_x} +\|\b\|_{L^2_x}
 \Big )\|\nab U(\tau'/2)F\|_{L^2_x}
 $$
 Taking the $L^q_t$-norm,
 \beaa
 \|\phi_2(\tau')\|_{L^q_t L^r_x}\les 
 \|\nab U(\tau'/2)F\|_{L^{\frac {2q}{2-q}}_t L^2_x} \les  
 \|\nab U(\tau'/2)F\|_{L^{\infty}_t L^2_x} 
 \eeaa
We now estimate the term $ \|\nab U(\tau'/2)F\|_{L^\infty_t L^2_x} $.  Squaring,  and integrating
by parts in $t$ and commuting   $\ddd_L$ with  $\nab $  and then with $U(\tau'/2)$ and
then using the heat flow  estimates    we derive,
\beaa
 \|\nab U(\tau'/2)F\|^2_{L^\infty_t L^2_x}&\les & 
 \int_{\cal H} \nab_L  \nab U(\tau'/2)F\c \nab U(\tau'/2)F + (\N_1[F])^2\\
 &\les &  \int_{\cal H} \nab \nab_L   U(\tau'/2)F\c \nab U(\tau'/2)F\\ & +&
 \|\,[\nab, \nab_L] U(\tau'/2)F\|_{L^1_t L^2_x} 
 \|\nab U(\tau'/2)F\|_{L^\infty_t L^2_x}  + 
  (\N_1[F])^2\\
 &\les &  \int_{\cal H} \nab    U(\tau'/2) \nab_L F\c \nab U(\tau'/2)F \\ &+&
 \|\,[\nab_L, U(\tau'/2)] F\|_{L^1_t L^2_x} \|\Delta  U(\tau'/2)F\|_{L^\infty_t L^2_x}
\\ &+&  \|\,[\nab, \nab_L] U(\tau'/2)F\|_{L^1_t L^2_x} 
 \|\nab U(\tau'/2)F\|_{L^\infty_t L^2_x}  + 
  (\N_1[F])^2\\ &\les & 
  \langle{\tau'}^{-\f12}\rangle (\N_1[F])^2 + \|\,[\nab_L, U(\tau'/2)] F\|_{L^1_t L^2_x} \|\Delta  U(\tau'/2)F\|_{L^\infty_t L^2_x}
\\ &+&  \|\,[\nab, \nab_L] U(\tau'/2)F\|_{L^1_t L^2_x} 
 \|\nab U(\tau'/2)F\|_{L^\infty_t L^2_x} 
\eeaa
Thus, since $$\|\lap U(\tau'/2) F\|_{L_t^\infty L_x^2}\les(\tau')^{-1}\| F\|_{L_t^\infty L_x^2}\les
(\tau')^{-1}\NN_1[F],$$
\bea
\|\nab U(\tau'/2)F\|_{L^\infty_t L^2_x}&\les & 
  \langle{\tau'}^{-\frac 14}\rangle \N_1[F]+ {\tau'}^{-\frac 12}
\|\,[\nab_L, U(\tau'/2)] F\|^\f12_{L^1_t L^2_x} \c(\N_1[F])^\f12\nn\\ &+&  
  \|\,[\nab, \nab_L] U(\tau'/2)F\|_{L^1_t L^2_x} \label{eq:tau'-nabU}
\eea
The results of lemmas \ref{lem:Phi1-twist} and \ref{lem:Phi2-twist}
imply that
$$
\|\,[\nab_L, U(\tau'/2)] F\|_{L^1_t L^2_x} \les \max ({\tau'}^{\frac 12-},\,{\tau'}^{\f12+})\, \N_1[F]
$$
On the other hand, observe that we  can write  the commutator  formula 
for 
 $[\nab,\ddd_L]$ in the form $[\nab,\ddd_L]=\Gdh\c\nab+\b\c$.
Therefore remembering the definition of $\phi_1(\tau)$,
\beaa
[\nab, \nab_L] U(\tau'/2)F\approx\phi_1(\tau'/2)
\eeaa
Thus, in view of \ref{eq:finally1},
\beaa
 \|\,[\nab, \nab_L] U(\tau'/2)F\|_{L^1_t L^2_x} \les  \langle (\tau')^{0-}\rangle\c 
\N_1[F]
\eeaa
Returning to \eqref{eq:tau'-nabU} we infer that
$$
\|\nab U(\tau'/2) F\|_{L^\infty_t L^2_x}\les \langle{\tau'}^{-\frac 14-}\rangle\, 
\N_1[F].
$$
It then follows that 
$$
\|\phi_2(\tau')\|_{L^q_t L^r_x}
 \les  
 \|\nab U(\tau'/2)F\|_{L^{\infty}_t L^2_x} \les
 \langle{\tau'}^{-\frac 14-}\rangle\,\c 
\N_1[F].
$$
Proceeding as in lemma \ref{lem:Phi2-straight} we obtain for 
$r>1$ sufficiently close to $r=1$,
\beaa
\|\Phi_2(\tau)\|_{L^q_t L^2_x} &\les &\|\int_0^\tau U(\tau-\tau') 
\phi_2(\tau')\|_{L^q_t L^2_x}\\ &\les & \int_0^\tau (\tau-\tau')^{-(\frac 1r-\f12)}
\|\phi_2(\tau')\|_{L^q_t L^r_x}\\ &\les & \int_0^\tau (\tau-\tau')^{-(\frac 1r-\f12)}
\langle{\tau'}^{-\frac 14-}\rangle\, \N_1[F]\\ &\les &
\max (\tau^{\frac 54-\frac 1r-},\, \tau^{\frac 32-\frac 1r})\, N_1[F]\\
&\les & \max (\tau^{\frac 14-},\, \tau^{\frac 12})\, N_1[F]
\eeaa
\end{proof}
We can now  finish the proof of proposition \ref{prop:Commut-nabL}.

 Combining the results of lemmas \ref{lem:Phi1-straight} and \ref{lem:Phi2-twist} 
 we derive
 \beaa
 \|[P_k, \nab_L] F\|_{L^q_t L^2_x} &\les & \int_0^\infty m_k (\tau) 
  \max (\tau^{\frac 14-}, \, \tau^{\f12} ) \, \N_1[F]\nn\\ &\les & 
 \int_0^1 m_k (\tau) \,\tau^{\frac 14-}\, \N_1[F] + 
 \int_1^\infty m_k(\tau) \,\tau^{\f12 +}\N_1[F]\nn\\ &\les & 
 2^{-\frac k2+}\int_0^{2^{2k}} m(\tau) \,\tau^{\frac 14-}\, \N_1[F] + 
 2^{-k}\int_{2^{2k}}^\infty m(\tau) \,\tau^{\f12 }\N_1[F]\nn\\&\les &
 2^{-\frac k2+}\,\N_1[F]
 \eeaa
 Arguing as in the case of the $L^1_t L^2_x$ estimates can we also obtain 
 that 
 $$
 \|\nab [P_k, \nab_L] F\|_{L^q_t L^2_x} \les  2^k 2^{-\frac k2+}\,\N_1[F]
 $$
 \end{proof}

\end{document}